\documentclass[12pt]{article}
\usepackage{amsfonts}
\usepackage{full page,
amssymb, amscd, amsmath, graphicx, 
}
\newtheorem{thm}{Theorem}[section]

\newtheorem{prop}[thm]{Proposition}
\newtheorem{cor}[thm]{Corollary}
\newtheorem{lemma}[thm]{Lemma}
\newtheorem{rema}[thm]{Remark}

\newtheorem{assum}[thm]{Assumption}

\newtheorem{data}[thm]{Data}

\newcommand{\halmos}{\rule{1ex}{1.4ex}}

\newcommand{\nn}{\nonumber \\}

 \newcommand{\res}{\mbox{\rm Res}}

\renewcommand{\hom}{\mbox{\rm Hom}}

 \newcommand{\pf}{{\it Proof.}\hspace{2ex}}
 
 \newcommand{\epfv}{\hspace*{\fill}\mbox{$\halmos$}\vspace{1em}}

\newcommand{\wt}{\mbox{\rm wt}\,}
\newcommand{\swt}{\mbox{\rm {\scriptsize wt}}\,}

\newcommand{\A}{\mathcal{A}}

\newcommand{\C}{\mathbb{C}}
\newcommand{\Z}{\mathbb{Z}}
\newcommand{\R}{\mathbb{R}}

\newcommand{\N}{\mathbb{N}}

\newcommand{\one}{\mathbf{1}}

\title{ {\bf A construction of lower-bounded 
generalized twisted modules for a grading-restricted vertex (super)algebra} }
\date{}
\author{Yi-Zhi Huang}

\begin{document}

\bibliographystyle{alpha}
\maketitle
\begin{abstract}
We give a general, direct and explicit 
construction of lower-bounded generalized twisted modules satisfying a universal property 
for a grading-restricted 
vertex (super)algebra $V$ associated to an automorphism $g$ of $V$. 
In particular, when $g$ is the identity, we obtain lower-bounded 
generalized $V$-modules satisfying a universal property. 
Let $W$ be a lower-bounded graded vector space
equipped with a set of ``generating twisted fields'' and a set of
``generator twist fields'' satisfying a weak commutativity for generating twisted fields,
a generalized weak commutativity for one generating twisted field 
and one generator twist field and 
some other properties that are relatively easy to verify. We first prove the convergence and commutativity 
of products of an  arbitrary number of generating twisted fields, one twist generator 
field and an arbitrary number of generating fields for $V$. Then using the convergence and 
commutativity, we define a twisted vertex 
operator map for $W$ and prove that  $W$ equipped with this twisted vertex operator map
is a lower-bounded generalized $g$-twisted $V$-module. Using this result, we
give an explicit construction of lower-bounded generalized $g$-twisted $V$-modules
satisfying a universal property 
starting from vector spaces graded by weights, $\Z_{2}$-fermion numbers and 
$g$-weights (eigenvalues of $g$) and real numbers corresponding to the lower 
bounds of the weights of the modules to be constructed. In particular, every lower-bounded generalized 
$g$-twisted $V$-module (every lower-bounded generalized $V$-module when $g$ is the identity)
 is a quotient of such a universal lower-bounded generalized $g$-twisted $V$-module
(a universal lower-bounded generalized $V$-module). 
\end{abstract}

\renewcommand{\theequation}{\thesection.\arabic{equation}}
\renewcommand{\thethm}{\thesection.\arabic{thm}}
\setcounter{equation}{0}
\setcounter{thm}{0}
\section{Introduction}

In the representation theory of associative algebras and Lie algebras, modules
satisfying universal properties (for examples, free modules, Verma modules and so on) 
in suitable categories of modules play a fundamental role. Modules in these categories 
are quotients of these biggest or universal modules and therefore can be studied using 
these modules whose structures are relatively simple to understand. 

In the representation theory of vertex (operator) (super)algebras and conformal field theory, 
finding a construction of modules 
satisfying universal properties in suitable categories of modules has been a long-standing problem. 
Finding such a construction will provide us with a powerful tool and will allow us to use the powerful 
homological algebra techniques (for example, the construction and applications of 
resolutions of modules) for the study of modules for 
vertex (operator) (super)algebras. To study the fixed-point subalgebra 
of a vertex (operator) (super)algebra under a group of automorphisms,
we have to construct and study twisted modules. 
It is also a long-standing problem in the case that 
the automorphism is of finite order to find such a construction of twisted modules satisfying 
universal properties in suitable categories of twisted modules. 

In this paper, we give a general, direct and explicit 
construction of lower-bounded generalized
twisted modules satisfying a universal property for a grading-restricted vertex (super)algebra
$V$ associated to an automorphism $g$ of $V$.  
In particular, in the case that $g$ is the identity, our construction 
give a general, direct and explicit construction of lower-bounded generalized
$V$-modules.  Our construction is for an arbitrary 
automorphism of the algebra. In particular,  in the case that the 
automorphism of the algebra is of infinite order and does not act on the algebra
semisimply, our construction gives  lower-bounded generalized twisted modules whose twisted 
vertex operators in general involve logarithm of the variable.

Twisted modules associated to automorphisms of finite order
of a vertex operator algebra were introduced by
Frenkel, Lepowsky and Meurman in their construction \cite{FLM1} \cite{FLM2}
\cite{FLM}
of the moonshine module vertex operator algebra $V^{\natural}$.  
Twisted module associated a general automorphism $g$ of 
 a vertex operator algebra $V$ were introduced by
the author in \cite{H-log-twisted-mod}. One of the main conjecture in the 
representation theory of vertex operator algebras
is that for a suitable vertex operator algebra $V$ and  a finite group 
$G$ of automorphisms of $V$,
the category of $g$-twisted $V$-modules for 
all $g\in G$ has a natural structure of $G$-crossed braided 
tensor category satisfying additional properties (see \cite{H-problems}). 
This conjecture follows from another stronger conjecture stating that 
twisted intertwining operators  (see \cite{H-twisted-int}) among $g$-twisted $V$-modules 
for $g\in G$ satisfy associativity, commutativity 
and modular invariance property (see also \cite{H-problems}). 
The second conjecture corresponds to a construction of orbifold conformal field theories
 and its solution will certainly depend on a deep 
understanding of twisted $V$-modules.

Twisted modules for
vertex (operator) (supper)algebras have been constructed and studied in many papers
(see for example, \cite{Le1}, \cite{FLM2}, 
\cite{Le2}, \cite{FLM}, \cite{D}, \cite{DL},
\cite{DonLM1}, \cite{DonLM2},  \cite{Li}, \cite{BDM}, \cite{DoyLM1},
\cite{DoyLM2}, \cite{BHL}, \cite{H-log-twisted-mod},
\cite{B}, \cite{Y}, and the references in these papers). But these constructions and studies are 
for special classes of vertex operator algebras and/or special classes of 
automorphisms. To study twisted modules, twisted intertwining operators
and the category of twisted modules, we need a general construction
of twisted modules. In principle, twisted modules can be constructed 
using the the functors constructed in \cite{HY} from categories of modules for 
the associative algebras introduced in \cite{DonLM1}
(for automorphisms of $V$ of finite orders) and in \cite{HY} (for 
general automorphisms) to suitable categories of twisted modules. But this indirect approach
is very difficult to use in general because the  abstract 
functors in \cite{HY} from the categories of modules for the associative algebras to the 
categories of suitable twisted modules are not equivalence of categories. 
It is therefore important 
to have a general, direct and explicit construction of suitable twisted modules satisfying universal properties. 
As we mentioned above, finding such a construction is a long-standing problem 
even in the case that the automorphism is of finite order or is even the identity. 
We solve this problem in this paper in the category of lower-bounded generalized 
$g$-twisted $V$-module for a general automorphism $g$ of a general grading-restricted 
vertex (super)algebra $V$. 

The approach used in our construction is the one that 
the author developed for the first construction 
of grading-restricted vertex algebras in 
\cite{H-2-const}. But the construction of lower-bounded generalized 
twisted modules in this paper, especially of those 
twisted modules whose twisted vertex operators involving the logarithm of the variable, 
is much more difficult than the one in \cite{H-2-const}, because the twisted vertex operators
are multivalued and because we do not have skew-symmetry for twisted modules 
(even for modules).
Besides twisted vertex operators, one 
crucial ingredient in the construction in this paper is
the twist vertex operators introduced and studied in \cite{H-twist-vo}.

Our construction is divided into two steps. We first prove a general construction theorem
which will be very useful also for the constructions of grading-restricted twisted modules,
twisted modules and other types of 
lower-bounded generalized twisted modules.
Let $W$ be a lower-bounded graded vector space
equipped with a set of ``generating twisted fields'' and a set of
``generator twist fields'' satisfying a weak commutativity for generating twisted fields,
a generalized weak commutativity for one generating twisted field 
and one generator twist field and 
some other properties that are relatively easy to verify.
We first prove the convergence and commutativity 
of products of an  arbitrary number of generating twisted fields, one twist generator 
field and an arbitrary number of generating fields for $V$.
Then using the convergence and 
commutativity, we define a twisted vertex 
operator map for $W$
and prove the construction theorem that  $W$ equipped with this twisted vertex operator map
is a lower-bounded generalized $g$-twisted $V$-module. If  $W$ is grading-restricted,
we obtain a grading-restricted twisted module and if in addition the operator 
$L_{W}(0)$ acts on $W$ semisimply, we obtain a twisted module. 
In the special case that $g=1_{V}$, we obtain
lower-bounded generalized modules, grading-restricted generalized 
modules and modules.

Then using this construction theorem, we
give an explicit construction of a lower-bounded generalized $g$-twisted $V$-module 
$\widehat{M}^{[g]}_{B}$ 
satisfying a universal property 
starting from a vector space $M$
graded by weights, $\Z_{2}$-fermion numbers and 
$g$-weights (eigenvalues of $g$) and a real number $B$ less than 
the real parts of all weights of homogeneous elements of $M$. 
The real number $B$ is in fact a lower bound of the weights
of $\widehat{M}^{[g]}_{B}$ and, 
roughly speaking, $M$ 
together with the algebra $V$ 
gives the generators of $\widehat{M}^{[g]}_{B}$.
In particular, every lower-bounded generalized 
$g$-twisted $V$-module (every lower-bounded generalized $V$-module when $g$ is the identity)
 is a quotient of such a universal lower-bounded generalized $g$-twisted $V$-module
(a universal lower-bounded generalized $V$-module). 

The construction and results obtained in this paper can be used to 
study a number of problems in the 
representation theory of vertex (operaor) (super)algebras. We shall 
discuss these apllications in future papers. We shall also construct 
and study examples of twisted modules for lattice, affine Lie and Virasoro vertex operator algebras
in future papers using the construction and results in this paper. 

The formulations and construction in the present paper are based on the formulations and 
results in 
\cite{H-twist-vo}. We refer the reader to \cite{H-twist-vo} for the basic 
definitions of grading-restricted vertex (super)algebra, 
generalized twisted modules and variants and twist vertex operator, 
conventions on formal and complex variables, and results on twist vertex operators
and their proofs.

This paper is organized as follows:
In Section 2, assuming that $V$ is generated by a set of fields $\{\phi^{i}(x)\}_{i\in I}$
and $g$ is an automorphism of $V$,
we introduce our data, a graded vector space $W$ with an action of $g$, 
a set $\{\phi_{W}^{i}(x)\}_{i\in I}$ of generating twisted fields 
and a set $\{\psi_{W}^{a}(x)\}_{a\in A}$ of generator twist fields
and two operators $L_{W}(0)$ and $L_{W}(-1)$, and our assumptions 
on these data, including, in particular, 
a weak commutativity for generating twisted fields and 
a generalized weak commutativity for one generating twisted field 
and one generator twist field. We also give a number of 
immediate consequences of these assumptions in this section. In Section 3,
we prove that the weak commutativity and generalized weak commutativity 
mentioned above are equivalent to the convergence and commutativity 
of products of an  arbitrary number of generating twisted fields, one twist generator 
field and an arbitrary number of generating fields for $V$. 
In Section 4, we define a twisted vertex operator map $Y^{g}_{W}$ for $W$ using 
convergence and commutativity above and  prove our construction theorem 
that $W$ equipped with $Y^{g}_{W}$
is a lower-bounded generalized $g$-twisted $V$-module. 
In Section 5, starting from a vector space $M$ graded by weights, $\Z_{2}$-fermion numbers
and $g$-weights (eigenvalues of $g$) and a real number $B$ less than 
the real parts of all weights of homogeneous elements of $M$ , we construct 
 a lower-bounded generalized $g$-twisted $V$-module $\widehat{M}^{[g]}_{B}$
and prove that it satisfies
a universal property. We also state in this section the consequence that 
every lower-bounded generalized 
$g$-twisted $V$-module
 is a quotient of a universal lower-bounded generalized $g$-twisted $V$-module.

\paragraph{Acknowledgments}
The author is grateful to Jason Saied for questions on
the construction of modules for grading-restricted vertex algebras
using the approach in \cite{H-2-const}.

\renewcommand{\theequation}{\thesection.\arabic{equation}}
\renewcommand{\thethm}{\thesection.\arabic{thm}}
\setcounter{equation}{0}
\setcounter{thm}{0}
\section{Generating twisted fields and generator twist fields}

In this section, we introduce the basic assumptions needed in our construction theorem
and state some immediate consequences. 

In the present paper, we fix a grading-restricted vertex superalgebra
$V$ and an automorphism $g$ of $V$. Then $V=\coprod_{\alpha\in P_{V}}V^{[\alpha]}$,
where $V^{[\alpha]}$ is the generalized eigenspace for $g$ with eigenvalue $e^{2\pi i\alpha}$
and $P_{V}$ is the subset of $\{\alpha\in \C\;|\; \Re(\alpha)\in [0, 1), 
e^{2 \pi i\alpha} \;\text{is an eigenvalue of}\;g\}$. By Lemma 
2.5 in \cite{H-twist-vo}, there exists an operator $\mathcal{L}_{g}$ with the
semisimple and nilpotent parts $\mathcal{S}_{g}$ and $\mathcal{N}_{g}$, respectively,
on $V$ such that $g=e^{2\pi i\mathcal{L}_{g}}=
e^{2\pi i(\mathcal{S}_{g}+\mathcal{N}_{g})}$
and by Proposition 2.6 in \cite{H-twist-vo}, both $e^{2\pi i\mathcal{S}_{g}}$
and $e^{2\pi i\mathcal{N}_{g}}$ are also automorphisms of $V$. 
Generalized eigenvectors for $g$ are eigenvectors for 
$e^{2\pi i\mathcal{S}_{g}}$ with the same eigenvalues and $\mathcal{N}_{g}$ is a 
derivation of $V$ (Proposition 2.6 in \cite{H-twist-vo}). 

We first state our assumptions on $V$.

\begin{assum}\label{algebra}
{\rm We assume that $V$ is generated by  
$\phi^{i}(x)=Y_{V}(\phi^{i}_{-1}\one, x)$ for $i\in I$, where
$\phi^{i}(x)$ or $\phi^{i}_{-1}\one$ for $i\in I$ are 
homogeneous with respect to weights and $\Z_{2}$-fermion numbers and
where 
$\phi^{i}_{-1}$ is the constant term of $\phi^{i}(x)$ and 
$\phi^{i}_{-1}\one=\lim_{x\to 0}\phi^{i}(x)\one$ (see \cite{H-2-const} for more details.)
For $i\in I$, $\phi^{i}_{-1}\one$ 
is a generalized eigenvector of $g$ with eigenvalue
$e^{2\pi i\alpha_{i}}$. 
We also assume that for $i\in I$, either  $\mathcal{N}_{g}\phi^{i}_{-1}\one=0$ 
or there exists 
$\mathcal{N}_{g}(i)\in I$ such that $\mathcal{N}_{g}\phi^{i}_{-1}\one
=\phi^{\mathcal{N}_{g}(i)}_{-1}\one$.  }
\end{assum}

We denote 
the weights and the $\Z_{2}$-fermion numbers 
of $\phi^{i}(x)$ or $\phi^{i}_{-1}\one$ for $i\in I$
by $\wt \phi^{i}$ and $|\phi^{i}|$, respectively. 

Since
$\phi^{i}(x)=Y_{V}(\phi^{i}_{-1}\one, x)$, we have 
$$e^{2\pi i\mathcal{S}_{g}}\phi^{i}(x)e^{-2\pi i\mathcal{S}_{g}}
=e^{2\pi i\alpha_{i}}\phi^{i}(x)$$
and
\begin{align*}
[\mathcal{N}_{g}, \phi^{i}(x)]&=[\mathcal{N}_{g}, Y_{V}(\phi^{i}_{-1}\one, x)]\nn
&=Y_{V}(\mathcal{N}_{g}\phi^{i}_{-1}\one, x).
\end{align*}
For convenience,
we shall use $\phi^{0}(x)$ to denote the $0$ vertex operator $Y_{V}(0, x)=0$, 
add $0$ to the index set $I$ and denote the index set with $0$ added still by $I$. 
Then for $i\in I$, there exists $\mathcal{N}_{g}(i)\in I$ such that 
$\mathcal{N}_{g}\phi^{i}_{-1}\one
=\phi^{\mathcal{N}_{g}(i)}_{-1}\one$, or equivalently, 
$[\mathcal{N}_{g}, \phi^{i}(z)]=\phi^{\mathcal{N}_{g}(i)}(z).$
Since 
$\mathcal{N}_{g}$ is nilpotent, there exists $K$ such that $\mathcal{N}_{g}^{K}(i)=0$.

Next we give the data needed in our construction theorem (Theorem \ref{const-thm})
in Section 4.

\renewcommand{\labelenumi}{$($\alph{enumi}$)$}
\begin{data}\label{data}
{\rm 
\begin{enumerate}

\item\label{data-a} Let 
$$W = \coprod_{n \in \C, s\in \Z_{2}, [\alpha]\in \C/\Z} W_{[n]}^{s; [\alpha]}
=\coprod_{n \in \C, s\in \Z_{2}, \alpha\in P_{W}} W_{[n]}^{s; [\alpha]}$$
 be a ${\C}\times \Z_{2}\times \C/\Z$-graded
vector space such that $W_{[n]}=\coprod_{s\in \Z_{2},  
\alpha\in P_{W}}W^{s; [\alpha]}_{[n]}=0$ 
when the real part of $n$ is sufficiently negative, where $P_{W}$ is the subset of 
the set $\{\alpha\in \C\;|\;\Re(\alpha)\in [0, 1)\}$ such that $W_{[n]}^{s; [\alpha]}\ne 0$
for $\alpha\in P_{W}$.


\item\label{data-b} Let 
\begin{align*}
\phi^{i}_{W}: W&\to x^{-\alpha^{i}}W((x))[\log  x]\nn
w&\mapsto \phi^{i}_{W}(x)w=\sum_{k\in \N}\sum_{n\in \alpha^{i}+\Z}
(\phi^{i}_{W})_{n, k}wx^{-n-1}
(\log x)^{k}
\end{align*}
for $i\in I$ be a set of linear maps called the {\it generating twisted field maps}.
Since $\phi^{i}_{W}(x)w\in x^{-\alpha^{i}}W((x))[\log  x]$, 
we must have $(\phi^{i}_{W})_{n, k}w=0$ when $n-\alpha^{i}$ is sufficiently 
negative and $k$ is sufficiently large. These linear maps correspond to
multivalued analytic maps with the preferred 
branch $\phi_{W}^{i; 0}$ and labeled branches $\phi_{W}^{i; p}$ for 
$p\in \Z$ from $\C^{\times}$
to $\hom(W, \overline{W})$. 

\item\label{data-c} Let 
\begin{align*}
\psi_{W}^{a}: V&\to \sum_{\alpha\in P_{V}}x^{-\alpha}W((x))[\log x]\nn
v&\mapsto \phi^{a}_{W}(x)v=\sum_{k\in \N}\sum_{\alpha\in P_{V},\; n\in \alpha+\Z}
(\psi^{a}_{W})_{n, k}vx^{-n-1}
(\log x)^{k}
\end{align*} 
for $a\in A$ be a set of linear maps called the {\it generator twist field maps}
such that $\phi^{a}_{W}(x)v\in x^{-\alpha}W((x))[\log x]$ for $\alpha\in P_{V}$ and 
$v\in V^{[\alpha]}$. 
Since $\phi^{a}_{W}(x)v\in x^{-\alpha}W((x))[\log  x]$ for $v\in V^{[\alpha]}$, 
we must have $(\psi^{a}_{W})_{n, k}v=0$ when $n-\alpha$ is sufficiently 
negative and $k$ is sufficiently large.
These linear maps corresponds to
multivalued analytic maps with preferred branch $\psi^{a; 0}$ and labeled branches
$\psi^{a; p}$ for $p\in \Z$ from $\C^{\times}$
to $\hom(V, \overline{W})$. 

\item\label{data-d} Let $L_{W}(0)$ and $L_{W}(-1)$ be operators on $W$. 

\item\label{data-e} An action of $g$ on $W$, denoted still by $g$, and 
an operator, still denoted by $\mathcal{L}_{g}$ and its
semisimple and nilpotent parts, still denoted by $\mathcal{S}_{g}$ and $\mathcal{N}_{g}$, 
respectively, on $W$ such that $g=e^{2\pi i\mathcal{L}_{g}}=
e^{2\pi i(\mathcal{S}_{g}+\mathcal{N}_{g})}$ on $W$.

\end{enumerate}}
\end{data}

These data are assumed to satisfy the following properties:

\renewcommand{\labelenumi}{\arabic{enumi}.}
\begin{assum}\label{basic-properties}
The space $W$, the generating twisted field maps  
 $\phi^{i}_{W}$ for $i\in I$,  the generator twist field maps
$\psi_{W}^{a}$ for $a\in A$, the operators $L_{W}(0)$, $L_{W}(-1)$,
$g$,  $\mathcal{L}_{g}$, $\mathcal{S}_{g}$ and $\mathcal{N}_{g}$ on $W$
in Data \ref{data} have the following properties:

\begin{enumerate}

\item\label{property-1} 
There exist semisimple and nilpotent operators 
$L_{W}(0)_{S}$ and $L_{W}(0)_{N}$ on $W$ such that
$L_{W}(0)=L_{W}(0)_{S}+L_{W}(0)_{N}$. 
For $i\in I$, 
$[L_{W}(0), \phi_{W}^{i}(x)]=z\frac{d}{dx}\phi_{W}^{i}(x)+(\wt \phi^{i})\phi_{W}^{i}(x)$.
For $a\in A$, there exists $(\wt \psi_{W}^{a})\in \C$ and, when $L_{W}(0)_{N}\psi_{W}^{a}(x)\ne 0$,
there exists $L_{W}(0)_{N}(a)\in A$ such that
$L_{W}(0)\psi_{W}^{a}(x)-\psi^{a}(x)L_{V}(0)
=x\frac{d}{dx}\psi_{W}^{a}(x)+(\wt \psi_{W}^{a})\psi_{W}^{a}(x)+\psi_{W}^{L_{W}(0)_{N}(a)}(x)$,
where $\psi_{W}^{L_{W}(0)_{N}(a)}(x)=0$ when $L_{W}(0)_{N}\psi_{W}^{a}(x)= 0$.

\item\label{property-2} For $i\in I$,
$[L_{W}(-1), \phi_{W}^{i}(x)]=\frac{d}{dx}\phi_{W}^{i}(x)$
and for $a\in A$,
$L_{W}(-1)\psi_{W}^{a}(x)-\psi_{W}^{a}(x)L_{V}(-1)=\frac{d}{dx}\psi_{W}^{a}(x)$.

\item\label{property-3} For $a\in A$, $\psi_{W}^{a}(x)\one\in W[[x]]$ and its constant terms
$\lim_{x\to 0}\psi_{W}^{a}(x)\one$ 
is homogeneous with respect to weights,  $\Z_{2}$-fermion number and $g$-weights.

\item\label{property-4} The vector space $W$ is spanned by elements of the form
$(\phi_{W}^{i_{1}})_{n_{1}, l_{1}}\cdots (\phi_{W}^{i_{k}})_{n_{k}, l_{k}}(\psi_{W}^{a})_{n, l}
v$ for $i_{1}, \dots, i_{k}\in I$, $a\in \A$
and $n_{1}\in \alpha^{i_{1}}+\Z, \dots, n_{k}\in \alpha^{i_{k}}+\Z$, 
$n\in \C$, $l_{1}, \dots, l_{k}, l\in \N$, $v\in V$.

\item\label{property-5} (i) For $i\in I$, $g\phi_{W}^{i;p+1}(z)g^{-1}=\phi_{W}^{i; p}(z)$.
(ii) For $i\in I$, $\phi_{W}^{i}(x)=
x^{-\mathcal{N}_{g}}(\phi_{W}^{i})_{0}(x)x^{\mathcal{N}_{g}}$ and for $a\in A$, 
$\psi_{W}^{a}(x)=(\psi_{W}^{a})_{0}(x)x^{-\mathcal{N}_{g}}$
where 
$(\phi_{W}^{i})_{0}(x)$ and $(\psi_{W}^{a})_{0}(x)$ are the constant terms 
of $\phi_{W}^{i}(x)$ and $\psi_{W}^{a}(x)$, respectively,  viewed
as power series of $\log x$ (with coefficients being series in powers of $x$). 
(iii)  For $i\in I$, $e^{2\pi i\mathcal{S}_{g}}
\phi_{W}^{i}(z)e^{-2\pi i\mathcal{S}_{g}}=e^{2\pi \alpha^{i}}\phi_{W}^{i}(z)$ and 
$[\mathcal{N}_{g}, \phi_{W}^{i}(z)]=\phi_{W}^{\mathcal{N}_{g}(i)}(z).$
(iv) For $a\in A$, 
there exists $\alpha^{a}\in P_{W}$ such that $(\psi_{W}^{a})_{n, 0}\one$ for 
$n\in -\N-1$ are generalized eigenvectors of $g$ 
with eigenvalue $e^{2\pi i \alpha^{a}}$. 

\item\label{property-6}  For $i, j\in I$, there exists $M_{ij}\in \Z_{+}$ such that
\begin{equation}\label{phi-weak-comm}
(x_{1}-x_{2})^{M_{ij}}\phi_{W}^{i}(x_{1})\phi_{W}^{j}(x_{2})
=(x_{1}-x_{2})^{M_{ij}}(-1)^{|\phi^{i}||\phi^{j}|}\phi_{W}^{j}(x_{2})\phi_{W}^{i}(x_{1}).
\end{equation}

\item\label{property-7}  For $i\in I$ and $a\in A$, there exists $M_{ia}\in \Z_{+}$ such that
\begin{align}\label{psi-gen-weak-comm}
&(x_{1}-x_{2})^{\alpha_{i}+M_{ia}}
(x_{1}-x_{2})^{\mathcal{N}_{g}}\phi_{W}^{i}(x_{1})(x_{1}-x_{2})^{-\mathcal{N}_{g}}
\psi_{W}^{a}(x_{2})\nn
&\quad =(-x_{2}+x_{1})^{\alpha_{i}+M_{ia}}
 (-1)^{|\phi^{i}||\psi^{a}|}
\psi_{W}^{a}(x_{2})
(-x_{2}+x_1)^{\mathcal{N}_{g}} \phi^{i}(x_{1})
(-x_{2}+x_1)^{-\mathcal{N}_{g}}.
\end{align}

\end{enumerate}
\end{assum}

For a multivalued analytic function $\phi(z)$ with a 
preferred branch of $z$ with domain $\C^{\times}$, we 
use $\phi^{p}(e^{a}z)$ to denote composition of the single-valued analytic 
function of $l_{p}(z)$ given by the branch $\phi^{p}(z)$ and the analytic map given by 
$l_{p}(z) \mapsto l_{p}(z)+a$. 

We have some  immediate consequences:

\renewcommand{\labelenumi}{\arabic{enumi}.}
\begin{prop}\label{properties}
The space $W$, the maps  $\phi_{W}^{i}$ for $i\in I$, $\psi_{W}^{a}$ for $a\in A$,
$L_{W}(-1)$ have the following
properties:

\begin{enumerate}
\setcounter{enumi}{7}

\item\label{property-8} For $c\in \C$, $i\in I$ and $a\in A$, 
$$e^{cL^{g}_{W}(0)}\phi_{W}^{i;p}(z)e^{-cL_{V}(0)}
=e^{c(\swt \phi^{i})}\phi_{W}^{i;p}(e^{c}z)$$
and
$$e^{cL^{g}_{W}(0)}\psi^{a;p}(z)e^{-cL_{V}(0)}=e^{c(\swt \psi_{W}^{a})}\psi^{a;p}(e^{c}z).$$

\item\label{property-9} For $i_{1}, \dots, i_{k}\in I$, $a\in A$
and $n_{1}, \dots, n_{k}\in \C$, $l_{1}, \dots, l_{k}, l\in \N$, $n\in \Z$ and $v\in V$,
we have 
\begin{align*}
L&_{W}(0)(\phi^{i_{1}}_{W})_{n_{1}, l_{1}}
\cdots (\phi^{i_{k}}_{W})_{n_{k}, l_{k}}(\psi_{W}^{a})_{n, l}v\nn
&=\sum_{j=1}^{k}(\phi^{i_{1}}_{W})_{n_{1}, l_{1}}
\cdots (\phi^{i_{j-1}}_{W})_{n_{j-1}, l_{j-1}} \cdot\nn
&\quad\quad\cdot\left((-n_{j}-1)
(\phi^{i_{j}}_{W})_{n_{j}, l_{j}}
+(l_{j}+1)(\phi^{i_{j}}_{W})_{n_{j}, l_{j}+1}
+(\wt \phi^{i_{j}})(\phi^{i_{j}}_{W})_{n_{j}, l_{j}}\right)\cdot\nn
&\quad\quad\cdot
(\phi^{i_{j+1}}_{W})_{n_{j+1}, l_{j+1}}
\cdots (\phi^{i_{k}}_{W})_{n_{k}, l_{k}}(\psi_{W}^{a})_{n, l}v\nn
&\quad +(\phi^{i_{1}}_{W})_{n_{1}, l_{1}}
\cdots (\phi^{i_{k}}_{W})_{n_{k}, l_{k}} \cdot\nn
&\quad\quad\cdot
\left((-n-1)(\psi_{W}^{a})_{n, l}+(l+1)
(\psi_{W}^{a})_{n, l+1}
+(\wt \psi_{W}^{a})(\psi_{W}^{a})_{n, l}
+(\psi_{W}^{L_{W}(0)_{N}(a)})_{n, l}\right)v\nn
&\quad +(\phi^{i_{1}}_{W})_{n_{1}, l_{1}}
\cdots (\phi^{i_{k}}_{W})_{n_{k}, l_{k}}
(\psi_{W}^{a})_{n, l}L_{V}(0)v,
\end{align*}
and 
\begin{align*}
L&_{W}(-1) (\phi_{W}^{i_{1}})_{n_{1}, l_{1}}\cdots (\phi_{W}^{i_{k}})_{n_{k}, l_{k}}
(\psi_{W}^{a})_{n, l}v\nn
&= \sum_{j=1}^{k}
(\phi_{W}^{i_{1}})_{n_{1}, l_{1}}\cdots (\phi^{i_{j-1}}_{W})_{n_{j-1}, l_{j-1}}\cdot\nn
&\quad\quad\cdot 
\left(-n_{j}(\phi^{i_{j}}_{W})_{n_{j}-1, l_{j}}
+(l_{j}+1)(\phi^{i_{j}}_{W})_{n_{j}-1, l_{j}+1}\right)
(\phi^{i_{j+1}}_{W})_{n_{j+1}, l_{j+1}}\cdots 
(\phi^{i_{k}}_{W})_{n_{k}, j_{k}}(\psi_{W}^{a})_{n, l}v\nn
&\quad +(\phi_{W}^{i_{1}})_{n_{1}, l_{1}}\cdots (\phi_{W}^{i_{k}})_{n_{k}, l_{k}}
\left(-n(\psi_{W}^{a})_{n-1, l}+(l+1)(\psi_{W}^{a})_{n-1, l+1}\right)v\nn
&\quad +(\phi_{W}^{i_{1}})_{n_{1}, l_{1}}\cdots (\phi_{W}^{i_{k}})_{n_{k}, l_{k}}
(\psi_{W}^{a})_{n, l}L_{V}(-1)v.
\end{align*}

\item\label{property-10} For $c\in \C$, $z\in \C^{\times}$ satisfying $|z|>|c|$, $i\in I$ and $a\in A$,
$e^{cL_{W}(-1)}\phi_{W}^{i; p}(z)e^{-cL_{W}(-1)}=\phi_{W}^{i;p}(z+c)$ and 
$e^{cL_{W}(-1)}\psi^{a;p}(z)e^{-cL_{V}(-1)}=\psi^{a;p}(z+c)$.

\item\label{property-11} The operator $L_{W}(-1)$ has weight $1$ and its adjoint $L_{W}(-1)'$ as an operator
on $W'$ has weight $-1$. In particular, $e^{zL_{W}(-1)'}w'\in W'$ for $z\in \C$ and $w'\in W'$.

\item \label{property-12} For $i\in I$, there exists $K, N\in \N$ such that 
$\phi_{W}^{i}(x)=\sum_{k=0}^{K}\sum_{n\in \alpha^{i}+N-\N}
(\phi_{W}^{i})_{n, k}x^{-n-1}(\log x)^{k}$ and for 
$i\in I$  and $w'\in W'$,  $\langle w', 
\phi_{W}^{i}(x)\cdot\rangle$ has only 
finitely many terms containing
$x^{-\alpha^{i}+n}$ for $n\in \Z_{+}$.  
For $a\in A$ and $v_{1}, \dots, v_{l}, v\in V$, 
$\psi^{a}_{W}(x)Y_{V}(v_{1}, x_{1})\cdots Y_{V}(v_{l}, x_{l})v$
is a polynomial in $\log x$ and for $a\in A$, $w'\in W'$ and $\alpha\in P_{V}$, 
$\langle w', 
\psi^{a}_{W}(x)\cdot\rangle$ as a formal series in $x$  
with coefficients in $(V^{[\alpha]})^{*}[[\log x]]$ is of the form 
$\sum_{n\in \alpha+N-\N}\lambda_{n}(\log x)x^{-n-1}$ for some $N\in \N$ and 
$\lambda_{n}\in (V^{[\alpha]})^{*}[[\log x]]$. 

\end{enumerate}
\end{prop}
\pf
These properties follow immediately from Assumption \ref{algebra}, Data \ref{data}
and Properties \ref{property-1}--\ref{property-7}
in Assumption \ref{basic-properties}.
\epfv

\begin{rema}
{\rm The twist fields $\psi_{W}^{a}$ might look mysterious 
to the reader but are in fact crucial in the present paper. 
The  reason  why they  are  important  for  the  construction  of  twisted modules
(and in particular, modules) is explained in the introduction in \cite{H-twist-vo}.   
We repeat the explanation here.  If we start with 
only generating twisted fields satisfying the 
weak commutativity or commutativity, we can still define a 
twisted vertex operator map but can only prove the 
commutativity and weak commutativity.  For  twisted  modules  (or even modules), 
the associativity is the main property to be verified and is not a consequence of
the weak commutativity or commutativity.  If we already have a twisted module, 
then the twist vertex operator map studied in \cite{H-twist-vo} 
changes the associativity of the twisted vertex operators to the 
commutativity involving twisted vertex operators, twist vertex operators and 
vertex operators for the algebra. 
Thus  in our construction, by  introducing  generator twist  fields  $\psi_{W}^{a}$ and  
assuming  that  the  commutativity  holds
for the generating twisted fields, generator twist fields and 
generating fields for the algebra $V$, we are able to 
prove the associativity and construct twisted modules  using  the  same  method  as  in  the  first  construction  of  grading-restricted vertex algebras in \cite{H-2-const}.}
\end{rema}

\renewcommand{\theequation}{\thesection.\arabic{equation}}
\renewcommand{\thethm}{\thesection.\arabic{thm}}
\setcounter{equation}{0}
\setcounter{thm}{0}
\section{Convergence and commutativity}

The construction theorem in the next
section uses the approach developed in \cite{H-2-const}. 
In this approach, the twisted vertex operators shall be defined and proved to 
be well defined using the 
correlation functions obtained from the product of the 
generating twisted fields $\phi^{i}_{W}(x)$, the generator twist fields
$\psi^{a}_{W}(x)$ and the vertex operators for $V$. 
We also need the commutativity of these correlation functions. 
Therefore to use this approach, we first have to prove the 
convergence and commutativity of these products. 
In this section, we prove these properties.

We first give the convergence and commutativity for $\phi_{W}^{i}$. In fact, we show
that the  convergence and commutativity for $\phi_{W}^{i}$ are equivalent to 
Property \ref{property-6} in Assumption \ref{basic-properties}. 

\begin{thm}\label{phi-locality}
Let $W = \coprod_{n \in \C, s\in \Z_{2}, \alpha\in \C/\Z} W_{[n]}^{s; [\alpha]}$ be
a ${\C}\times \Z_{2} \times\C/\Z$-graded
vector space as in Data \ref{data} and
let 
\begin{align*}
\phi^{i}_{W}: W&\to x^{-\alpha^{i}}W((x))[\log  x]\nn
w&\mapsto \phi^{i}_{W}(x)w
\end{align*}
for $i\in I$ be a set of linear maps as in Data \ref{data}. Assume that they satisfy the parts
for $\phi_{W}^{i}$ for $i\in I$ in Properties \ref{property-1}, \ref{property-2}, 
\ref{property-5} in Assumption \ref{basic-properties}. 
Then Property \ref{property-6} in Assumption \ref{basic-properties} is equivalent to 
the following properties:
\begin{enumerate}
\setcounter{enumi}{12}

\item\label{property-13}  For $w'\in W'$, $w\in W$ and $i_{1}, \dots, i_{k}\in I$, 
the series
$$\langle w', \phi_{W}^{i_{1};p}(z_{1})\cdots \phi_{W}^{i_{k};p}(z_{k})w\rangle$$
is absolutely convergent in the region $|z_{1}|>\cdots>|z_{k}|>0$. 
Moreover, there exists a 
multivalued analytic function of the form 
$$\sum_{n_{1}, \dots, n_{k}=0}^{N}f_{n_{1}\cdots n_{k+l}}(z_{1}, \dots,
z_{k})  
z_{1}^{-\alpha_{i_{1}}}\cdots z_{k}^{-\alpha_{i_{k}}}
(\log z_{1})^{n_{1}}\cdots (\log z_{k})^{n_{k}},$$
denoted by 
$$F(\langle w', \phi_{W}^{i_{1}}(z_{1})\cdots \phi_{W}^{i_{k}}(z_{k})w\rangle),$$
where $N\in \N$ and $f_{n_{1}\cdots n_{k}}(z_{1}, \dots,
z_{k})$ for  $n_{1}, \cdots, n_{k}=0, \dots, N$ are rational functions
of $z_{1}, \dots, z_{k}$ with the only possible poles $z_{i}=0$ for
$i=1, \dots, k$, $z_{i}-z_{j}=0$
for $i, j=1, \dots, k$, $i\ne j$, such that its sum is equal to the branch
\begin{align*}
F&^{p}(\langle w', \phi_{W}^{i_{1}}(z_{1})\cdots \phi_{W}^{i_{k}}(z_{k})w\rangle)\nn
&=\sum_{n_{1}, \dots, n_{k}=0}^{N}f_{n_{1}\cdots n_{k}}(z_{1}, \dots,
z_{k})e^{-\alpha_{i_{1}}l_{p}(z_{1})}\cdots e^{-\alpha_{i_{k}}l_{p}(z_{k})}
(l_{p}(z_{1}))^{n_{1}}\cdots (l_{p}(z_{k}))^{n_{k}},
\end{align*}
of $F(\langle w', \phi_{W}^{i_{1}}(z_{1})\cdots \phi_{W}^{i_{k}}(z_{k})w\rangle)$ 
in the region 
given by $|z_{1}|>\cdots>|z_{k}|>0$. 
In addition,  the orders of the pole $z_{i}=0$ of the  rational functions 
$f_{n_{1}\cdots n_{k}}(z_{1}, \dots, z_{k})$ have a lower bound 
independent of $\phi^{i_{l}}$ for $i\ne i$ and $w'$;
the orders of the pole $z_{i}=z_{j}$ of the rational functions 
$f_{n_{1}\cdots n_{k}}(z_{1}, \dots, z_{k})$ 
have a lower bound independent of $\phi^{i_{l}}$ for $l\ne i, j$, 
$w$ and $w'$. 

\item\label{property-14} For $w\in W$,  $w'\in W'$, $i_{1}, i_{2}, i\in I$,
\begin{equation}\label{phi-w-i-comm}
F^{p}(\langle w', \phi_{W}^{i_{1}}(z_{1})\phi_{W}^{i_{2}}(z_{2})w\rangle)
=(-1)^{|\phi^{i_{1}}||\phi^{i_{2}}|} F^{p}(\langle v', \phi_{W}^{i_{2}}(z_{2})
\phi_{W}^{i_{1}}(z_{1})w\rangle).
\end{equation}

\end{enumerate}
\end{thm}
\pf
The proof that Properties \ref{property-13} and \ref{property-14} implies Property 
\ref{property-6} in Assumption \ref{basic-properties} is completely the same as the 
proof of Proposition 3.7 in \cite{H-twist-vo}. 
The proof that Property \ref{property-6} in Assumption \ref{basic-properties} implies 
Property \ref{property-13}  is completely the same as 
the proof of Theorem 3.10 in \cite{H-twist-vo}.
We refer the reader to those proofs in \cite{H-twist-vo}. 

Propertiy \ref{property-14} follows immediately from
Property \ref{property-13} in the case $k=2$ and (\ref{phi-weak-comm}).
Thus Property \ref{property-6} in Assumption \ref{basic-properties} also implies 
Property \ref{property-14}.
\epfv

The  result below is the most general convergence result and the commutativity 
involving generator twist fields. We in fact
prove that these convergence and commutativity properties are equivalent to 
Property \ref{property-7} in Assumption \ref{basic-properties} and thus by 
Theorem \ref{phi-locality}, together with the convergence and commutativity properties
in Theorem \ref{phi-locality},
are equivalent to Properties \ref{property-6} and \ref{property-7}
in Assumption \ref{basic-properties}.

\begin{thm}\label{psi-gen-locality}
Let $W = \coprod_{n \in \C, s\in \Z_{2}, \alpha\in \C/\Z} W_{[n]}^{s; [\alpha]}$ be
a ${\C}\times \Z_{2} \times \C/\Z$-graded
vector space as in Data \ref{data} and
let 
\begin{align*}
\phi^{i}_{W}: W&\to x^{-\alpha^{i}}W((x))[\log  x]\nn
w&\mapsto \phi^{i}_{W}(x)w
\end{align*}
for $i\in I$ and 
\begin{align*}
\psi_{W}^{a}: V&\to \sum_{\alpha\in P_{V}}x^{-\alpha}W((x))[\log x]\nn
v&\mapsto \phi^{a}_{W}(x)v
\end{align*} 
for $a\in A$ be  linear maps as in Data \ref{data}. Assume that they satisfy Properties 
\ref{property-1}--\ref{property-6} in 
Assumption \ref{basic-properties}. 
Then Property \ref{property-7} is equivalent to the following 
properties:
\begin{enumerate}
\setcounter{enumi}{14}

\item\label{property-15}  For $w'\in W'$, $v\in V^{[\alpha]}$ and 
$i_{1}, \dots, i_{k+l}\in I$, $a\in A$, 
the series
$$\langle w', \phi_{W}^{i_{1};p}(z_{1})\cdots \phi_{W}^{i_{k};p}(z_{k})\psi^{a;p}(z)
\phi^{i_{k+1}}(z_{k+1})\cdots \phi^{i_{k+l}}(z_{k+l})v\rangle$$
is absolutely convergent in the region $|z_{1}|>\cdots>|z_{k}|>|z|>
|z_{k+1}|>\cdots >|z_{k+l}|>0$. Moreover, there exists a 
multivalued analytic function of the form
\begin{align}\label{k+l-prod-twted-twt-fd-1}
&\sum_{n_{1}, \dots, n_{k+l}, n=0}^{N}f_{n_{1}\cdots n_{k+l}n}(z_{1}, \dots,
z_{k+l}, z)\cdot\nn
&\quad\quad \cdot  
(z_{1}-z)^{-\alpha_{i_{1}}}\cdots (z_{k+l}-z)^{-\alpha_{i_{k+l}}}z^{-\alpha}
(\log (z_{1}-z))^{n_{1}}\cdots (\log (z_{k+l}-z))^{n_{k+l}}(\log z)^{n},
\end{align}
denoted by 
$$F(\langle w', \phi_{W}^{i_{1}}(z_{1})\cdots \phi_{W}^{i_{k}}(z_{k})\psi_{W}^{a}(z)
\phi^{i_{k+1}}(z_{k+1})\cdots \phi^{i_{k+l}}(z_{k+l})v\rangle),$$
where $N\in \N$ and $f_{n_{1}\cdots n_{k+l}n}(z_{1}, \dots,
z_{k+l}, z)$ for $n_{1}, \cdots, n_{k+l}, n=0, \dots, N$ are rational functions
of $z_{1}, \dots, z_{k+l}, z$ with the only possible poles $z_{i}=0$ for
$i=1, \dots, k+l$, $z=0$, $z_{i}-z_{j}=0$
for $i, j=1, \dots, k+l$, $i\ne j$, $z_{i}-z=0$ for $i=1, \dots, k+l$, 
such that its sum is equal to the branch
\begin{align}\label{k+l-prod-twted-twt-fd-2}
F&^{p}(\langle w', \phi_{W}^{i_{1}}(z_{1})\cdots \phi_{W}^{i_{k}}(z_{k})\psi_{W}^{a}(z)
\phi^{i_{k+1}}(z_{k+1})\cdots \phi^{i_{k+l}}(z_{k+l})v\rangle)\nn
&=\sum_{n_{1}, \dots, n_{k+l}, n=0}^{N}f_{n_{1}\cdots n_{k+l}n}(z_{1}, \dots,
z_{k+l})\cdot\nn
&\quad\quad \cdot  
e^{-\alpha_{i_{1}}l_{p}(z_{1}-z)}\cdots e^{-\alpha_{i_{k+l}}l_{p}(z_{k+l}-z)}
e^{-\alpha l_{p}(z)}
(l_{p}(z_{1}-z))^{n_{1}}\cdots (l_{p}(z_{k+l}-z))^{n_{k+l}}(l_{p}(z))^{n}
\end{align}
of (\ref{k+l-prod-twted-twt-fd-1}) in the region 
given by $|z_{1}|>\cdots>|z_{k}|>|z|>
|z_{k+1}|>\cdots >|z_{k+l}|>0$, $|\arg (z_{i}-z)-\arg z_{i}|<\frac{\pi}{2}$ for 
$i=1, \dots, k$ and $|\arg (z_{i}-z)-\arg z|<\frac{\pi}{2}$ for 
$i=k+1, \dots, k+l$. 
In addition,   the orders of the pole $z_{j}=0$ of the rational functions 
$f_{n_{1}\cdots n_{k+l}n}(z_{1}, \dots,z_{k+l}, z)$ 
have a lower bound independent of $v_{q}$ for $q\ne j$, $w$ and $w'$; 
the orders of the pole $z=0$ of the rational functions 
$f_{n_{1}\cdots n_{k+l}n}(z_{1}, \dots, z_{k+l}, z)$ 
have a lower bound independent of $v_{1}, \dots, v_{k+l}$ and $w'$; 
the orders of the pole $z_{j}=z_{m}$ of the rational functions 
$f_{n_{1}\cdots n_{k+l}n}(z_{1}, \dots, z_{k+l}, z)$ 
have a lower bound independent of $v_{q}$ for $q\ne j, m$, $v$, $w$ and $w'$;
the orders of the pole $z_{j}=z$ of the rational functions 
$f_{n_{1}\cdots n_{k+l}n}(z_{1}, \dots, z_{k+l}, z)$ 
have a lower bound independent of $v_{q}$ for $q\ne j$, $v$ and $w'$.

\item\label{property-16} For $v\in V$,  $w'\in W'$, $i\in I$, $a\in A$,
\begin{equation}\label{psi-w-i-comm}
F^{p}(\langle w', \phi_{W}^{i}(z_{1})\psi_{W}^{a}(z_{2})v\rangle)
=(-1)^{|\phi^{i}||\psi_{W}^{a}|} F^{p}(\langle w', \psi_{W}^{a}(z_{2})
\phi^{i}(z_{1})v\rangle).
\end{equation}

\end{enumerate}
In particular, when Properties 
\ref{property-1}--\ref{property-5} in 
Assumption \ref{basic-properties} hold,
Properties 
\ref{property-6} and \ref{property-7} in 
Assumption \ref{basic-properties} are equivalent to Properties 
\ref{property-13} and \ref{property-14} in Theorem \ref{phi-locality}
and Properties \ref{property-15} and \ref{property-16}.
\end{thm}
\pf
Assume that Properties \ref{property-15} and \ref{property-16} hold. 
For $w'\in W'$ and $v\in V$,
by Property \ref{property-5} in Assumption \ref{basic-properties},
\begin{align}\label{psi-gen-locality-1}
\langle w', \phi_{W}^{i;p}(z_{1})\psi_{W}^{a;p}(z_{2})v\rangle&=
\langle w', e^{-l_{p}(z_{1})\mathcal{N}_{g}}(\phi_{W}^{i})_{0}^{p}(z_{1})
e^{l_{p}(z_{1})\mathcal{N}_{g}}\psi_{W}^{a;p}(z_{2})v\rangle\nn
&=\sum_{r\in \N}\frac{(-1)^{r}}{r!}(l_{p}(z_{1}))^{r}
\langle w', [\overbrace{\mathcal{N}_{g}, \cdots, [\mathcal{N}_{g}}^{r},
(\phi_{W}^{i})_{0}^{p}(z_{1})]\cdots]\psi_{W}^{a;p}(z_{2})v\rangle\nn
&=\sum_{r\in \N}\frac{(-1)^{r}}{r!}
\langle w', (\phi_{W}^{\mathcal{N}_{g}^{r}(i)})_{0}^{p}(z_{1})\psi_{W}^{a;p}(z_{2})v\rangle
(l_{p}(z_{1}))^{r},
\end{align}
where, in our notation,  for $j\in I$,
$(\phi_{W}^{j})_{0}^{p}(z_{1})$
denotes the $p$-th branch of $(\phi_{W}^{j})_{0}(z_{1})$.
By Property \ref{property-15}, the left-hand side 
of (\ref{psi-gen-locality-1}) 
is absolutely convergent in the region given by $|z_{1}|>|z_{2}|>0$ and 
$|\arg (z_{1}-z_{2})-\arg z_{1}|<\frac{\pi}{2}$ to 
\begin{equation}\label{correl-prod-phi-psi}
\sum_{j, k, l, q=0}^{N}a_{qjkl}z_{1}^{-t}
e^{-(\alpha_{i}+m_{q})l_{p}(z_{1}-z_{2})}e^{n_{j}l_{p}(z_{2})}
(l_{p}(z_{1}-z_{2}))^{k}(l_{p}(z_{2}))^{l},
\end{equation}
where $t, m_{q}\in \Z$ for $q=0, \dots, N$  and $n_{j}\in \C$ for $j=0, \dots, N$. 
But the expansion of (\ref{correl-prod-phi-psi}) in the region given by $|z_{1}|>|z_{2}|>0$ 
and $|\arg (z_{1}-z_{2})-\arg z_{1}|<\frac{\pi}{2}$
as a power series in $l_{p}(z_{1})$ is 
\begin{align}\label{correl-prod-phi-psi-1}
&\sum_{k, j, l, q=0}^{N}a_{qjkl}z_{1}^{-t}
e^{-(\alpha_{i}+m_{q})l_{p}(z_{1}-z_{2})}e^{n_{j}l_{p}(z_{2})}
(l_{p}(z_{2}))^{l}\sum_{r=0}^{k}\binom{k}{r}
\left(\log\left(1-\frac{z_{2}}{z_{1}}\right)\right)^{k-r}
(l_{p}(z_{1}))^{r}\nn
&\quad=\sum_{r=1}^{N}\sum_{k=r}^{N}\sum_{j, l, q=0}^{N}a_{qjkl}z_{1}^{-t}
e^{-(\alpha_{i}+m_{q})l_{p}(z_{1}-z_{2})}e^{n_{j}l_{p}(z_{2})}
(l_{p}(z_{2}))^{l}\binom{k}{r}
\left(\log\left(1-\frac{z_{2}}{z_{1}}\right)\right)^{k-r}
(l_{p}(z_{1}))^{r},
\end{align}
where $\left(\log\left(1-\frac{z_{2}}{z_{1}}\right)\right)^{k-r}$ is understood as 
the branch obtained by using the expansion in nonnegative powers of $\frac{z_{2}}{z_{1}}$. 
Comparing (\ref{correl-prod-phi-psi}) and (\ref{correl-prod-phi-psi-1}) and using 
Proposition 2.1 in \cite{H-zhang-correction}, we obtain 
\begin{align}\label{correl-prod-phi-psi-2}
&\frac{(-1)^{r}}{r!}
\langle w', (\phi_{W}^{\mathcal{N}_{g}^{r}(i)})_{0}^{p}(z_{1})\psi_{W}^{a;p}(z_{2})v\rangle\nn
&\quad =\sum_{k=r}^{N}\sum_{j, l, q=0}^{N}a_{qjkl}z_{1}^{-t}
e^{-(\alpha_{i}+m_{q})l_{p}(z_{1}-z_{2})}e^{n_{j}l_{p}(z_{2})}
(l_{p}(z_{2}))^{l}\binom{k}{r}
\left(\log\left(1-\frac{z_{2}}{z_{1}}\right)\right)^{k-r}
\end{align}
for $r\in \N$ and for $r>N$, both sides of (\ref{correl-prod-phi-psi-2}) are $0$.
Then in the region given by $|z_{1}|>|z_{2}|>0$ 
and $|\arg (z_{1}-z_{2})-\arg z_{1}|<\frac{\pi}{2}$,
\begin{align*}
\langle w', &e^{l_{p}(z_{1}-z_{2})\mathcal{N}_{g}}
\phi_{W}^{i;p}(z_{1})e^{-l_{p}(z_{1}-z_{2})\mathcal{N}_{g}}
\psi_{W}^{a;p}(z_{2})v\rangle\nn
&=\langle w', e^{l_{p}(z_{1}-z_{2})\mathcal{N}_{g}}
e^{-l_{p}(z_{1})\mathcal{N}_{g}}(\phi_{W}^{i})_{0}^{p}(z_{1})
e^{l_{p}(z_{1})\mathcal{N}_{g}}e^{-l_{p}(z_{1}-z_{2})\mathcal{N}_{g}}
\psi_{W}^{a;p}(z_{2})v\rangle\nn
&=\langle w', e^{\left(\log\left(1-\frac{z_{2}}{z_{1}}\right)\right)\mathcal{N}_{g}}
(\phi_{W}^{i})_{0}^{p}(z_{1})
e^{-\left(\log\left(1-\frac{z_{2}}{z_{1}}\right)\right)\mathcal{N}_{g}}
\psi_{W}^{a;p}(z_{2})v\rangle\nn
&=\sum_{r\in \N}\frac{1}{r!}\left(\log\left(1-\frac{z_{2}}{z_{1}}\right)\right)^{r}
\langle w', [\overbrace{\mathcal{N}_{g}, \cdots, [\mathcal{N}_{g}}^{r},
(\phi_{W}^{i})_{0}^{p}(z_{1})]\cdots]\psi_{W}^{a;p}(z_{2})v\rangle\nn
&=\sum_{r\in \N}\frac{1}{r!}
\langle w', (\phi_{W}^{\mathcal{N}_{g}^{r}(i)})_{0}^{p}(z_{1})\psi_{W}^{a;p}(z_{2})v\rangle
\left(\log\left(1-\frac{z_{2}}{z_{1}}\right)\right)^{r}\nn
&=\sum_{r=1}^{N}\sum_{k=r}^{N}\sum_{j, l, q=0}^{N}(-1)^{r}a_{qjkl}z_{1}^{-t}
e^{-(\alpha_{i}+m_{q})l_{p}(z_{1}-z_{2})}e^{n_{j}l_{p}(z_{2})}
(l_{p}(z_{2}))^{l}\cdot\nn
&\quad\quad\quad\quad\quad\quad\quad\cdot 
\binom{k}{r}\left(\log\left(1-\frac{z_{2}}{z_{1}}\right)\right)^{k-r}
\left(\log\left(1-\frac{z_{2}}{z_{1}}\right)\right)^{r}\nn
&=\sum_{k=1}^{N}\sum_{j, l, q=0}^{N}a_{qjkl}z_{1}^{-t}
e^{-(\alpha_{i}+m_{q})l_{p}(z_{1}-z_{2})}e^{n_{j}l_{p}(z_{2})}
(l_{p}(z_{2}))^{l}\cdot\nn
&\quad\quad\quad\quad\quad\quad\quad\cdot 
\left(\sum_{r=1}^{k}\binom{k}{r}\left(\log\left(1-\frac{z_{2}}{z_{1}}\right)\right)^{k-r}
(-1)^{r}\left(\log\left(1-\frac{z_{2}}{z_{1}}\right)\right)^{r}\right)
\end{align*}
\begin{align}\label{psi-gen-locality-2}
&=\sum_{k=1}^{N}\sum_{j, l, q=0}^{N}a_{qjkl}z_{1}^{-t}
e^{-(\alpha_{i}+m_{q})l_{p}(z_{1}-z_{2})}e^{n_{j}l_{p}(z_{2})}
(l_{p}(z_{2}))^{l}\cdot\nn
&\quad\quad\quad\quad\quad\quad\quad\cdot 
\left(\left(\log\left(1-\frac{z_{2}}{z_{1}}\right)\right)
-\left(\log\left(1-\frac{z_{2}}{z_{1}}\right)\right)\right)^{k}\quad\quad\quad\quad\quad\nn
&=\sum_{j, l, q=0}^{N}a_{qj0l}z_{1}^{-t}
e^{-(\alpha_{i}+m_{q})l_{p}(z_{1}-z_{2})}e^{n_{j}l_{p}(z_{2})}
(l_{p}(z_{2}))^{l}.
\end{align}

On the other hand, in the region given by $|z_{2}|>|z_{1}|>0$ 
and $|\arg (z_{1}-z_{2})-\arg z_{2}|<\frac{\pi}{2}$,
\begin{align}\label{psi-gen-locality-3}
(-1)&^{|\phi^{i}||\psi^{a}|}
\langle w', \psi_{W}^{a;p}(z_{2})e^{l_{p}(z_{2}-z_{1})\mathcal{N}_{g}}
e^{\pi i\mathcal{N}_{g}}\phi^{i}(z_{1})e^{-\pi i\mathcal{N}_{g}}
e^{-l_{p}(z_{2}-z_{1})\mathcal{N}_{g}}v\rangle\nn
&=(-1)^{|\phi^{i}||\psi^{a}|}
\langle w', \psi_{W}^{a;p}(z_{2})e^{(l_{p}(z_{2}-z_{1})+\pi i)\mathcal{N}_{g}}
\phi^{i}(z_{1})e^{-(l_{p}(z_{2}-z_{1})+\pi i)\mathcal{N}_{g}}v\rangle\nn
&=\sum_{k\in \N}\frac{1}{k!}(l_{p}(z_{2}-z_{1})+\pi i)^{k}(-1)^{|\phi^{i}||\psi^{a}|}
\langle w', \psi_{W}^{a;p}(z_{2})
[\overbrace{\mathcal{N}_{g}, \cdots, [\mathcal{N}_{g}}^{k}, \phi^{i;p}(z_{1})]\cdots]
v\rangle\nn
&=\sum_{k\in \N}\frac{1}{k!}(l_{p}(z_{2}-z_{1})+\pi i)^{k}(-1)^{|\phi^{i}||\psi^{a}|}
\langle w', \psi_{W}^{a;p}(z_{2})
\phi^{\mathcal{N}_{g}^{k}(i);p}(z_{1})
v\rangle.
\end{align}
By Property \ref{property-5} in Assumption \ref{basic-properties} and (\ref{correl-prod-phi-psi-2}),
in the region given by $|z_{1}|>|z_{2}|>0$ 
and $|\arg (z_{1}-z_{2})-\arg z_{1}|<\frac{\pi}{2}$, we have
\begin{align}\label{psi-gen-locality-4}
\langle w',& \phi_{W}^{\mathcal{N}_{g}^{k}(i);p}(z_{1})\psi_{W}^{a;p}(z_{2})v\rangle\nn
&=\langle w', e^{-l_{p}(z_{1})\mathcal{N}_{g}}
(\phi_{W}^{\mathcal{N}_{g}^{k}(i)})^{p}_{0}(z_{1})
e^{-l_{p}(z_{1})\mathcal{N}_{g}}\psi_{W}^{a;p}(z_{2})v\rangle\nn
&=\sum_{n\in \N}\frac{(-1)^{n}}{n!}(l_{p}(z_{1}))^{n}
\langle w', [\overbrace{\mathcal{N}_{g}, \cdots
[\mathcal{N}_{g}}^{n}, (\phi_{W}^{\mathcal{N}_{g}^{k}(i)})^{p}_{0}(z_{1})
]\cdots]\psi_{W}^{a;p}(z_{2})v\rangle\nn
&=\sum_{n\in \N}\frac{(-1)^{n}}{n!}(l_{p}(z_{1}))^{n}\langle w',
(\phi_{W}^{\mathcal{N}_{g}^{k+n}(i)})^{p}_{0}(z_{1})\psi_{W}^{a;p}(z_{2})v\rangle\nn
&=\sum_{n\in \N}\frac{(-1)^{k}(k+n)!}{n!}(l_{p}(z_{1}))^{n}
\sum_{s=k+n}^{N}\sum_{j, l, q=0}^{N}a_{qjsl}z_{1}^{-t}
e^{-(\alpha_{i}+m_{q})l_{p}(z_{1}-z_{2})}e^{n_{j}l_{p}(z_{2})}
(l_{p}(z_{2}))^{l}\cdot\nn
&\quad\quad\quad\quad\cdot \binom{s}{k+n}
\left(\log\left(1-\frac{z_{2}}{z_{1}}\right)\right)^{s-(k+n)}\nn
&=\sum_{n\in \N}\sum_{s=k+n}^{N}\sum_{j, l, q=0}^{N}\frac{(-1)^{k}(k+n)!}{n!}
a_{qjsl}z_{1}^{-t}
e^{-(\alpha_{i}+m_{q})l_{p}(z_{1}-z_{2})}e^{n_{j}l_{p}(z_{2})}
(l_{p}(z_{2}))^{l}\cdot\quad\nn
&\quad\quad\quad\quad\cdot \binom{s}{k+n}
(l_{p}(z_{1}))^{n}\left(\log\left(1-\frac{z_{2}}{z_{1}}\right)\right)^{s-(k+n)}.
\end{align}
By Property \ref{property-16},
\begin{equation}\label{psi-gen-locality-5}
(-1)^{|\phi^{i}||\psi^{a}|}
\langle w',\psi_{W}^{a;p}(z_{2}) \phi_{W}^{\mathcal{N}_{g}^{k}(i);p}(z_{1})v\rangle
\end{equation}
and the left-hand side of (\ref{psi-gen-locality-4}) converges absolutely to the same branch 
of a multivalued analytic function in different regions. Then we see 
from (\ref{psi-gen-locality-4}) that 
in the region given by $|z_{2}|>|z_{1}|>0$ 
and $|\arg (z_{1}-z_{2})-\arg z_{2}|<\frac{\pi}{2}$,
(\ref{psi-gen-locality-5}) must converge absolutely to 
\begin{align}
\sum_{n\in \N}&\sum_{s=k+n}^{N}\sum_{j, l, q=0}^{N}\frac{(-1)^{k}(k+n)!}{n!}
a_{qjsl}z_{1}^{-t}
e^{-(\alpha_{i}+m_{q})l_{p}(z_{1}-z_{2})}e^{n_{j}l_{p}(z_{2})}
(l_{p}(z_{2}))^{l}\cdot\nn
&\quad\quad\quad\quad\cdot \binom{s}{k+n}
(l_{p}(z_{1}))^{n}\left(\log\left(1-\frac{z_{2}}{z_{1}}\right)\right)^{s-(k+n)}.
\end{align}
Thus the right-hand side of (\ref{psi-gen-locality-3}) is equal to 
\begin{align}\label{psi-gen-locality-6}
\sum_{k\in \N}&\sum_{n\in \N}\sum_{s=k+n}^{N}\sum_{j, l, q=0}^{N}\frac{(-1)^{k}(k+n)!}{k!n!}
a_{qjsl}z_{1}^{-t}
e^{-(\alpha_{i}+m_{q})l_{p}(z_{1}-z_{2})}e^{n_{j}l_{p}(z_{2})}
(l_{p}(z_{2}))^{l}\cdot\nn
&\quad\quad\quad\cdot \binom{s}{k+n}(l_{p}(z_{2}-z_{1})+\pi i)^{k}
(l_{p}(z_{1}))^{n}\left(\log\left(1-\frac{z_{2}}{z_{1}}\right)\right)^{s-(k+n)}\nn
&=\sum_{j, s, l, q=0}^{N}\sum_{r=0}^{s}
a_{qjsl}z_{1}^{-t}
e^{-(\alpha_{i}+m_{q})l_{p}(z_{1}-z_{2})}e^{n_{j}l_{p}(z_{2})}
(l_{p}(z_{2}))^{l}\cdot\nn
&\quad\quad\quad \cdot \binom{s}{r}
\left(\sum_{k\in \N}\binom{r}{k}(-1)^{k}(l_{p}(z_{2}-z_{1})+\pi i)^{k}
(l_{p}(z_{1}))^{s-k}\right)\left(\log\left(1-\frac{z_{2}}{z_{1}}\right)\right)^{s-r}\nn
&=\sum_{j, s, l, q=0}^{N}
a_{qjsl}z_{1}^{-t}
e^{-(\alpha_{i}+m_{q})l_{p}(z_{1}-z_{2})}e^{n_{j}l_{p}(z_{2})}
(l_{p}(z_{2}))^{l}\cdot\nn
&\quad\quad \quad \cdot \sum_{r=0}^{s}\binom{s}{r}
(l_{p}(z_{1})-(l_{p}(z_{2}-z_{1})+\pi i))^{s}
\left(\log\left(1-\frac{z_{2}}{z_{1}}\right)\right)^{s-r}\nn
&=\sum_{j, s, l, q=0}^{N}
a_{qjsl}z_{1}^{-t}
e^{-(\alpha_{i}+m_{q})l_{p}(z_{1}-z_{2})}e^{n_{j}l_{p}(z_{2})}
(l_{p}(z_{2}))^{l}\cdot\nn
&\quad\quad\quad \cdot \left(l_{p}(z_{1})-(l_{p}(z_{2}-z_{1})+\pi i)+
\log\left(1-\frac{z_{2}}{z_{1}}\right)\right)^{s}\nn
&=\sum_{j, l, q=0}^{N}
a_{qj0l}z_{1}^{-t}
e^{-(\alpha_{i}+m_{q})l_{p}(z_{1}-z_{2})}e^{n_{j}l_{p}(z_{2})}
(l_{p}(z_{2}))^{l}.
\end{align}

From (\ref{psi-gen-locality-2}) and the calculations from (\ref{psi-gen-locality-3}) 
to (\ref{psi-gen-locality-6}), we see that if we choose $M_{ia}\in \Z_{+}$ 
to be larger than $m_{q}$ for $q=0, \dots, N$, we have 
\begin{align*}
&e^{(\alpha_{i}+M_{ia})l_{p}(z_{1}-z_{2})} 
\langle w', e^{l_{p}(z_{1}-z_{2})\mathcal{N}_{g}}
\phi_{W}^{i;p}(z_{1})e^{-l_{p}(z_{1}-z_{2})\mathcal{N}_{g}}
\psi_{W}^{a;p}(z_{2})v\rangle\nn
&\quad =\sum_{j, l, q=0}^{N}
a_{qj0l}z_{1}^{-t}
(z_{1}-z_{2})^{M_{ia}-m_{q}}e^{n_{j}l_{p}(z_{2})}
(l_{p}(z_{2}))^{l}
\end{align*}
\begin{align}\label{psi-gen-locality-7}
&\quad\;\; =e^{(\alpha_{i}+M_{ia})l_{p}(z_{2}-z_{1})} e^{\pi i(\alpha_{i}+M_{ia})}
(-1)^{|\phi^{i}||\psi^{a}|}\cdot\nn
&\;\;\quad\quad\quad \cdot
\langle w', \psi_{W}^{a;p}(z_{2})e^{l_{p}(z_{2}-z_{1})\mathcal{N}_{g}}
e^{\pi i\mathcal{N}_{g}}\phi^{i}(z_{1})e^{-\pi i\mathcal{N}_{g}}
e^{-l_{p}(z_{2}-z_{1})\mathcal{N}_{g}}v\rangle
\end{align}
The formula (\ref{psi-gen-locality-7}) for all $w'\in W'$ and $v\in V$
is equivalent to the formal identity 
(\ref{psi-gen-weak-comm}).

Next we prove that Property \ref{property-7} in 
Assumption \ref{basic-properties} implies  Properties \ref{property-15} and 
\ref{property-16}.
By (b) in Data \ref{data}
and Property \ref{property-12} in Proposition  \ref{properties},
there exist $N_{1}, N_{2}\in \Z$ and $K\in \N$ such that 
$$\langle w', \phi_{W}^{i}(x)w\rangle=\sum_{n=N_{1}}^{N_{2}}\sum_{k=1}^{K}
\langle w', (\phi_{W}^{i})_{\alpha^{i}+n, k}w\rangle x^{-(\alpha^{i}+n)-1}(\log x)^{k}$$
for $w\in W$ and $w'\in W'$. Then  for $w'\in W'$ and $v\in V$,
\begin{align*}
(x_{1}-x_{2})^{\alpha_{i}+M_{ia}} &\langle w', 
(x_{1}-x_{2})^{\mathcal{N}_{g}}\phi_{W}^{i}(x_{1})(x_{1}-x_{2})^{-\mathcal{N}_{g}}
\psi_{W}^{a}(x_{2})v\rangle\nn
& \in 
\C((x_{1}^{-1}, x_{2}))[\log x_{2}],\nn
(x_{2}-x_{1})^{\alpha_{i}+M_{ia}} e^{\pi i(\alpha_{i}+M_{ia})}(-1)^{|\phi^{i}||\psi^{a}|}
& \langle w', \psi_{W}^{a}(x_{2})
(x_{2}-x_1)^{\mathcal{N}_{g}}e^{\pi i\mathcal{N}_{g}} \phi^{i}(x_{1})
e^{-\pi i\mathcal{N}_{g}}(x_{2}-x_1)^{-\mathcal{N}_{g}}v\rangle \nn
& \in
\C((x_{1}, x_{2}^{-1}))[\log x_{2}],
\end{align*}
where we use $\C((x_{1}^{-1}, x_{2}))$ ($\C((x_{1}, x_{2}^{-1}))$)
to denote the ring of Laurent series in $x_{1}$ and $x_{2}$ having only
finitely many terms in positive powers of $x_{1}$ and negative 
powers of $x_{2}$ (finitely many terms in negative powers of $x_{1}$ and positive 
powers of $x_{2}$). But by (\ref{psi-gen-weak-comm}) these two formal series are equal. 
So they must belong to 
$$\C[x_{1}, x_{1}^{-1}, x_{2}, x_{2}^{-1}][\log x_{1}, \log x_{2}].$$
This formal series can be written as
$$\sum_{j, l, q=0}^{N}
b_{qjl}x_{1}^{-t}
(x_{1}-x_{2})^{M_{ia}-m_{q}}x_{2}^{n_{j}}
(\log (x_{2}))^{l}$$
for $b_{qjl}\in \C$, $t\in \N$, $m_{q}\in \Z$ such that $m_{q}\le M_{ia}$, $n_{j}\in \C$
so that
\begin{align*}
&(x_{1}-x_{2})^{\alpha_{i}+M_{ia}} \langle w', 
(x_{1}-x_{2})^{\mathcal{N}_{g}}\phi_{W}^{i}(x_{1})(x_{1}-x_{2})^{-\mathcal{N}_{g}}
\psi_{W}^{a}(x_{2})v\rangle\nn
&\quad =\sum_{j, l, q=0}^{N}
b_{qjl}x_{1}^{-t}
(x_{1}-x_{2})^{M_{ia}-m_{q}}x_{2}^{n_{j}}
(\log (x_{2}))^{l}\nn
&\quad =(x_{2}-x_{1})^{\alpha_{i}+M_{ia}} e^{\pi i(\alpha_{i}+M_{ia})}(-1)^{|\phi^{i}||\psi^{a}|}\cdot\nn
&\quad\quad\quad \cdot
\langle w', \psi_{W}^{a}(x_{2})
(x_{2}-x_1)^{\mathcal{N}_{g}}e^{\pi i\mathcal{N}_{g}} \phi^{i}(x_{1})
e^{-\pi i\mathcal{N}_{g}}(x_{2}-x_1)^{-\mathcal{N}_{g}}v\rangle.
\end{align*}
Thus 
\begin{align}\label{psi-gen-locality-8}
&F^{p}(\langle w', 
(z_{1}-z_{2})^{\mathcal{N}_{g}}\phi_{W}^{i}(z_{1})(z_{1}-z_{2})^{-\mathcal{N}_{g}}
\psi_{W}^{a}(z_{2})v\rangle)\nn
&\quad =\sum_{j, l, q=0}^{N}
b_{qjl}z_{1}^{-t}
e^{-(\alpha_{i}+m_{q})l_{p}(z_{1}-z_{2})}e^{n_{j}l_{p}(z_{2})}
(l_{p}(z_{2}))^{l}\nn
&\quad =(-1)^{|\phi^{i}||\psi^{a}|}
F^{p}(\langle w', \psi_{W}^{a}(z_{2})
(z_{2}-z_1)^{\mathcal{N}_{g}}e^{\pi i\mathcal{N}_{g}} \phi^{i}(z_{1})
e^{-\pi i\mathcal{N}_{g}}(z_{2}-z_1)^{-\mathcal{N}_{g}}v\rangle).
\end{align}

Using Part (iii) of Property \ref{property-5}  in Assumption \ref{basic-properties}, we obtain
\begin{align}\label{psi-gen-locality-10}
\langle w&', \phi_{W}^{i;p}(z_{1})\psi_{W}^{a; p}(z_{2})v\rangle\nn
&=\langle w', e^{-l_{p}(z_{1}-z_{2})\mathcal{N}_{g}}e^{l_{p}(z_{1}-z_{2})\mathcal{N}_{g}}
\phi_{W}^{i; p}(z_{1})
e^{-l_{p}(z_{1}-z_{2})\mathcal{N}_{g}}e^{l_{p}(z_{1}-z_{2})\mathcal{N}_{g}}
\psi_{W}^{a;p}(z_{2})v\rangle\nn
&=\sum_{s\in \N}\frac{(-1)^{s}}{s!}(l_{p}(z_{1}-z_{2}))^{s}
\langle w', [\overbrace{\mathcal{N}_{g}, \cdots, [\mathcal{N}_{g}}^{s},
e^{l_{p}(z_{1}-z_{2})\mathcal{N}_{g}}
\phi_{W}^{i; p}(z_{1})
e^{-l_{p}(z_{1}-z_{2})\mathcal{N}_{g}}]\cdots]\psi_{W}^{a;p}(z_{2})v\rangle\nn
&=\sum_{s\in \N}\frac{(-1)^{s}}{s!}(l_{p}(z_{1}-z_{2}))^{s}
\langle w', 
e^{l_{p}(z_{1}-z_{2})\mathcal{N}_{g}}
[\overbrace{\mathcal{N}_{g}, \cdots, [\mathcal{N}_{g}}^{s}, \phi_{W}^{i; p}(z_{1})]\cdots]
e^{-l_{p}(z_{1}-z_{2})\mathcal{N}_{g}}\psi_{W}^{a;p}(z_{2})v\rangle\nn
&=\sum_{s=0}^{N}\frac{(-1)^{s}}{s!}(l_{p}(z_{1}-z_{2}))^{s}
\langle w', 
e^{l_{p}(z_{1}-z_{2})\mathcal{N}_{g}}
\phi_{W}^{\mathcal{N}_{g}^{s}(i); p}(z_{1})
e^{-l_{p}(z_{1}-z_{2})\mathcal{N}_{g}}\psi_{W}^{a;p}(z_{2})v\rangle.
\end{align}
By (\ref{psi-gen-locality-8}), the right-hand side of (\ref{psi-gen-locality-10})
is absolutely convergent in the region given by $|z_{1}|>|z_{2}|>0$ and 
$|\arg (z_{1}-z_{2})-\arg z_{1}|<\frac{\pi}{2}$ to 
\begin{equation}\label{psi-gen-locality-11}
\sum_{j, s, l, q=0}^{N}a_{qjsl}z_{1}^{-t}
e^{-(\alpha_{i}+m_{q})l_{p}(z_{1}-z_{2})}e^{n_{j}l_{p}(z_{2})}(l_{p}(z_{1}-z_{2}))^{s}
(l_{p}(z_{2}))^{l}
\end{equation}
for some $a_{qjsl}\in \C$. 

On the other hand,
by Part (iii) of Property \ref{property-5}  in Assumption \ref{basic-properties},
\begin{align*}
(-1)&^{|\phi^{i}||\psi^{a}|}
\langle w', \psi_{W}^{a; p}(z_{2})\phi_{W}^{i;p}(z_{1})v\rangle\nn
&=(-1)^{|\phi^{i}||\psi^{a}|}\langle w',
\psi_{W}^{a;p}(z_{2}) e^{-\pi i\mathcal{N}_{g}}e^{-l_{p}(z_{2}-z_{1})\mathcal{N}_{g}}
e^{l_{p}(z_{2}-z_{1})\mathcal{N}_{g}}e^{\pi i\mathcal{N}_{g}}\cdot\nn
&\quad\quad\quad\cdot 
\phi_{W}^{i; p}(z_{1})e^{-\pi i\mathcal{N}_{g}}
e^{-l_{p}(z_{2}-z_{1})\mathcal{N}_{g}}e^{l_{p}(z_{2}-z_{1})\mathcal{N}_{g}}
e^{\pi i\mathcal{N}_{g}}v\rangle\nn
&=\sum_{s\in \N}\frac{(-1)^{s}}{s!}(l_{p}(z_{2}-z_{1})+\pi i)^{s}(-1)^{|\phi^{i}||\psi^{a}|}
\cdot\nn
&\quad\quad\quad\cdot \langle w', \psi_{W}^{a;p}(z_{2})
[\overbrace{\mathcal{N}_{g}, \cdots, [\mathcal{N}_{g}}^{s},
e^{l_{p}(z_{2}-z_{1})\mathcal{N}_{g}}e^{\pi i\mathcal{N}_{g}}
\phi_{W}^{i; p}(z_{1})
e^{-\pi i\mathcal{N}_{g}}
e^{-l_{p}(z_{2}-z_{1})\mathcal{N}_{g}}]\cdots]v\rangle\nn
&=\sum_{s\in \N}\frac{(-1)^{s}}{s!}(l_{p}(z_{2}-z_{1})+\pi i)^{s}
(-1)^{|\phi^{i}||\psi^{a}|}\cdot\nn
&\quad\quad\quad\cdot \langle w', \psi_{W}^{a;p}(z_{2})
e^{l_{p}(z_{2}-z_{1})\mathcal{N}_{g}}e^{\pi i\mathcal{N}_{g}}
[\overbrace{\mathcal{N}_{g}, \cdots, [\mathcal{N}_{g}}^{s},\phi_{W}^{i; p}(z_{1})]\cdots]
e^{-\pi i\mathcal{N}_{g}}
e^{-l_{p}(z_{2}-z_{1})\mathcal{N}_{g}}v\rangle
\end{align*}
\begin{align}\label{psi-gen-locality-12}
&\;\quad\quad =\sum_{s=0}^{N}\frac{(-1)^{s}}{s!}(l_{p}(z_{2}-z_{1})+\pi i)^{s}
\cdot\nn
&\;\quad\quad\quad\quad\quad\cdot (-1)^{|\phi^{i}||\psi^{a}|}\langle w', \psi_{W}^{a;p}(z_{2})
e^{l_{p}(z_{2}-z_{1})\mathcal{N}_{g}}e^{\pi i\mathcal{N}_{g}}
\phi_{W}^{\mathcal{N}_{g}^{s}(i); p}(z_{1})
e^{-\pi i\mathcal{N}_{g}}
e^{-l_{p}(z_{2}-z_{1})\mathcal{N}_{g}}v\rangle.
\end{align}
By (\ref{psi-gen-locality-8}), 
$$ (-1)^{|\phi^{i}||\psi^{a}|}\langle w', \psi_{W}^{a;p}(z_{2})
e^{l_{p}(z_{2}-z_{1})\mathcal{N}_{g}}e^{\pi i\mathcal{N}_{g}}
\phi_{W}^{\mathcal{N}_{g}^{s}(i); p}(z_{1})
e^{-\pi i\mathcal{N}_{g}}
e^{-l_{p}(z_{2}-z_{1})\mathcal{N}_{g}}v\rangle$$
is absolutely convergent in the region given by $|z_{2}|>|z_{1}|>0$ and 
$|\arg (z_{1}-z_{2})-\arg z_{2}|<\frac{\pi}{2}$ to 
$$F^{p}(\langle w', 
(z_{1}-z_{2})^{\mathcal{N}_{g}}\phi_{W}^{\mathcal{N}_{g}^{s}(i); p}(z_{1})
(z_{1}-z_{2})^{-\mathcal{N}_{g}}
\psi_{W}^{a}(z_{2})v\rangle).$$
Thus 
the right-hand side of (\ref{psi-gen-locality-12})
is absolutely convergent in the region given by $|z_{2}|>|z_{1}|>0$ and 
$|\arg (z_{1}-z_{2})-\arg z_{2}|<\frac{\pi}{2}$ to (\ref{psi-gen-locality-11}).
In particular, we have proved Property 6 in the case $k=1$ and $l=0$ and 
also (\ref{psi-w-i-comm}). By Proposition \ref{phi-locality}, (\ref{phi-w-i-comm})
also holds. Thus Property \ref{property-16} holds. 

We still need to prove Property \ref{property-15}
for general $k$ and $l$. 
For $i, j\in I$, 
\begin{align}\label{psi-gen-locality-13}
(x&_{1}-x)^{\mathcal{N}_{g}}
\phi_{W}^{i}(x_{1})(x_{1}-x)^{-\mathcal{N}_{g}}(x_{2}-x)^{\mathcal{N}_{g}}
\phi_{W}^{j}(x_{2})(x_{2}-x)^{-\mathcal{N}_{g}}\nn
&=\sum_{r, s\in \N}\frac{(-1)^{r+s}}{r!s!}(\log (x_{1}-x))^{r}(\log (x_{2}-x))^{s}
\phi_{W}^{\mathcal{N}_{g}^{r}(i)}(x_{1})\phi_{W}^{\mathcal{N}_{g}^{s}(i)}(x_{2})
\end{align}
Since the right-hand side of (\ref{psi-gen-locality-13}) is in fact a finite sum, 
there exists $N_{ij}\in \Z_{+}$ such that  
\begin{align}\label{psi-gen-locality-14}
(x&_{1}-x_{2})^{N_{ij}}(x_{1}-x)^{\mathcal{N}_{g}}
\phi_{W}^{i}(x_{1})(x_{1}-x)^{-\mathcal{N}_{g}}(x_{2}-x)^{\mathcal{N}_{g}}
\phi_{W}^{j}(x_{2})(x_{2}-x)^{-\mathcal{N}_{g}}\nn
&=\sum_{r, s\in \N}\frac{(-1)^{r+s}}{r!s!}(\log (x_{1}-x))^{r}(\log (x_{2}-x))^{s}
(x_{1}-x_{2})^{N_{ij}}
\phi_{W}^{\mathcal{N}_{g}^{r}(i)}(x_{1})\phi_{W}^{\mathcal{N}_{g}^{s}(j)}(x_{2})\nn
&=\sum_{r, s\in \N}\frac{(-1)^{r+s}}{r!s!}(\log (x_{1}-x))^{r}(\log (x_{2}-x))^{s}
\cdot\nn
&\quad\quad\quad\quad\cdot 
(-1)^{|\phi^{\mathcal{N}_{g}^{r}(i)}||\phi^{\mathcal{N}_{g}^{s}(i)}|}
(x_{1}-x_{2})^{N_{ij}}
\phi_{W}^{\mathcal{N}_{g}^{s}(j)}(x_{2})\phi_{W}^{\mathcal{N}_{g}^{r}(i)}(x_{1})\nn
&=(x_{1}-x_{2})^{N_{ij}}(x_{2}-x)^{\mathcal{N}_{g}}
\phi_{W}^{j}(x_{2})(x_{2}-x)^{-\mathcal{N}_{g}}(x_{1}-x)^{\mathcal{N}_{g}}
\phi_{W}^{i}(x_{1})(x_{1}-x)^{-\mathcal{N}_{g}}.
\end{align}

For general $k, l\in \N$, consider the series
\begin{align}\label{psi-gen-locality-15}
&\prod_{q=1}^{k+l}(x_{q}-x)^{\alpha_{i_{q}}+M_{i_{q}a}}
\prod_{m=1}^{l}(x-x_{k+m})^{\alpha_{i_{k+m}}+M_{i_{k+m}a}}
e^{\pi i(\alpha_{i_{k+m}}+M_{i_{k+m}a})}\cdot\nn
&\quad\cdot 
\langle w', (x_{1}-x)^{\mathcal{N}_{g}}
\phi_{W}^{i_{1}}(x_{1})(x_{1}-x)^{-\mathcal{N}_{g}}
\cdots (x_{k}-x)^{\mathcal{N}_{g}}\phi_{W}^{i_{k}}(x_{k})
(x_{k}-x)^{-\mathcal{N}_{g}}\cdot\nn
&\quad\quad\quad\quad\cdot
\psi_{W}^{a}(x)(x-x_{k+1})^{\mathcal{N}_{g}}e^{\pi i\mathcal{N}_{g}}
\phi_{W}^{i_{k+1}}(x_{k+1})e^{-\pi i\mathcal{N}_{g}}(x-x_{k+1})^{-\mathcal{N}_{g}}
\cdot\nn
&\quad\quad\quad\quad
\cdots (x-x_{k+l})^{\mathcal{N}_{g}}e^{\pi i\mathcal{N}_{g}}\phi_{W}^{i_{k+l}}(x_{k+l})
e^{-\pi i\mathcal{N}_{g}}(x-x_{k+l})^{-\mathcal{N}_{g}}v\rangle,
\end{align}
where $v\in V$ and $w'\in W'$. Using (\ref{psi-gen-weak-comm}),
we see that (\ref{psi-gen-locality-15}) is equal to 
\begin{align*}
&\prod_{q=1}^{k+l}(x_{q}-x)^{\alpha_{i_{q}}+M_{i_{q}a}}
\prod_{m=1}^{l}(x_{k+m}-x)^{\alpha_{i_{k+m}}+M_{i_{k+m}a}}\cdot\nn
&\quad\cdot 
\langle w', (x_{1}-x)^{\mathcal{N}_{g}}
\phi_{W}^{i_{1}}(x_{1})(x_{1}-x)^{-\mathcal{N}_{g}}
\cdots (x_{k}-x)^{\mathcal{N}_{g}}\phi_{W}^{i_{k}}(x_{k})\cdot
\end{align*}
\begin{align}\label{psi-gen-locality-16}
&\quad\quad\quad\quad\cdot
(x_{k}-x)^{-\mathcal{N}_{g}}
(x_{k+1}-x)^{\mathcal{N}_{g}}
\phi_{W}^{i_{k+1}}(x_{k+1})(x_{k+1}-x)^{-\mathcal{N}_{g}}\cdot\nn
&\quad\quad\quad\quad
\cdots (x_{k+l}-x)^{\mathcal{N}_{g}}\phi_{W}^{i_{k+l}}(x_{k+l})
(x_{k+l}-x)^{-\mathcal{N}_{g}}\psi_{W}^{a}(x)v\rangle,
\end{align}
By (c) in Data \ref{data}, 
(\ref{psi-gen-locality-16}) has only finitely many negative integral powers in $x$
and finitely many nonnegative integral power terms in  $\log x$. 
Using (\ref{psi-gen-weak-comm}) again, we see that (\ref{psi-gen-locality-15}) is also equal to 
\begin{align}\label{psi-gen-locality-17}
&\prod_{q=1}^{k+l}(x-x_{q})^{\alpha_{i_{q}}+M_{i_{q}a}}
e^{\pi i(\alpha_{i_{q}}+M_{i_{q}a})}
\prod_{m=1}^{l}(x-x_{k+m})^{\alpha_{i_{k+m}}+M_{i_{k+m}a}}
e^{\pi i(\alpha_{i_{k+m}}+M_{i_{k+m}a})}\cdot\nn
&\quad\cdot 
\langle w', \psi_{W}^{a}(x)(x-x_{1})^{\mathcal{N}_{g}}e^{\pi i\mathcal{N}_{g}}
\phi_{W}^{i_{1}}(x_{1})e^{-\pi i\mathcal{N}_{g}}(x-x_{1})^{-\mathcal{N}_{g}}\cdot\nn
&\quad\quad\quad\quad
\cdots (x-x_{k})^{\mathcal{N}_{g}}e^{\pi i\mathcal{N}_{g}}\phi_{W}^{i_{k}}(x_{k})
e^{-\pi i\mathcal{N}_{g}}(x-x_{k})^{-\mathcal{N}_{g}}\cdot\nn
&\quad\quad\quad\quad\cdot
(x-x_{k+1})^{\mathcal{N}_{g}}e^{\pi i\mathcal{N}_{g}}
\phi_{W}^{i_{k+1}}(x_{k+1})e^{-\pi i\mathcal{N}_{g}}(x-x_{k+1})^{-\mathcal{N}_{g}}
\cdot\nn
&\quad\quad\quad\quad
\cdots (x-x_{k+l})^{\mathcal{N}_{g}}e^{\pi i\mathcal{N}_{g}}\phi_{W}^{i_{k+l}}(x_{k+l})
e^{-\pi i\mathcal{N}_{g}}(x-x_{k+l})^{-\mathcal{N}_{g}}v\rangle,
\end{align}
By (c) in Data \ref{data} and Property \ref{property-12} in Proposition \ref{properties},
 (\ref{psi-gen-locality-17})  has only finitely 
many positive integral power terms 
in $x$ and finitely many nonnegative integral power terms in  $\log x$. Thus 
(\ref{psi-gen-locality-15}), (\ref{psi-gen-locality-16}) and (\ref{psi-gen-locality-17}) 
are a Laurent polynomial in $x$ and a polynomial in $\log x$. In particular, 
\begin{align}\label{psi-gen-locality-18}
\langle w',& (x_{1}-x)^{\mathcal{N}_{g}}
\phi_{W}^{i_{1}}(x_{1})(x_{1}-x)^{-\mathcal{N}_{g}}
\cdots (x_{k}-x)^{\mathcal{N}_{g}}\phi_{W}^{i_{k}}(x_{k})\cdot\nn
&\quad\cdot
(x_{k}-x)^{-\mathcal{N}_{g}}
(x_{k+1}-x)^{\mathcal{N}_{g}}
\phi_{W}^{i_{k+1}}(x_{k+1})(x_{k+1}-x)^{-\mathcal{N}_{g}}\cdot\nn
&\quad
\cdots (x_{k+l}-x)^{\mathcal{N}_{g}}\phi_{W}^{i_{k+l}}(x_{k+l})
(x_{k+l}-x)^{-\mathcal{N}_{g}}\psi_{W}^{a}(x)v\rangle
\end{align}
is equal to this polynomial in $x$, $x^{-1}$ and $\log x$
 with series in $x_{1}, \dots, x_{k+l}$ as coefficients
multiplied by 
$$\prod_{q=1}^{k+l}(x_{q}-x)^{-(\alpha_{i_{q}}+M_{i_{q}a})}
\prod_{m=1}^{l}(x-x_{k+m})^{-(\alpha_{i_{k+m}}+M_{i_{k+m}a})}
e^{-\pi i(\alpha_{i_{k+m}}+M_{i_{k+m}a})}.$$
The coefficients of this polynomial in $x$, $x^{-1}$ and $\log x$ are given by 
the coefficients of (\ref{psi-gen-locality-16}) in powers of $x$ and $\log x$. 
These coefficients are finite sums of products of 
(finite) linear combinations of powers of  $x_{1}, \dots, x_{k+l}$, 
$\log x_{1}, \dots, \log x_{k+l}$ and series of the form 
$$\langle w',
\phi_{W}^{\mathcal{N}_{g}^{j_{1}}(i_{1})}(x_{1})\cdots
\phi_{W}^{\mathcal{N}_{g}^{j_{k+l}}(i_{k+l})}(x_{k+l})
(\psi_{W}^{a})_{r, s}v\rangle$$
for $j_{1}, \dots, j_{k+l}\in \N$, $r\in \C$ and $s\in \N$. By Property \ref{property-13} in 
Theorem \ref{phi-locality}, we see that 
these coefficients with $x_{q}^{n}$ and $\log x_{q}$ substituted by 
$e^{nl_{p}(z_{q})}$ and  $l_{p}(z_{q})$, respectively, 
are absolutely convergent in the region $|z_{1}|>\cdots >|z_{k+l}|>0$
to analytic functions of the form in Property \ref{property-13} in 
Theorem \ref{phi-locality} with 
$k$ there replaced by $k+l$. Thus (\ref{psi-gen-locality-18})
with $x_{q}^{n}$, $x^{n}$,
$\log x_{q}$ and $\log x$ substituted by $e^{nl_{p}(z_{q})}$, $e^{n\l_{p}(z)}$,
$l_{p}(z_{q})$ and  $l_{p}(z)$, respectively,
is absolutely convergent in the region $|z_{1}|>\cdots>|z_{k}|>|z|>
|z_{k+1}|>\cdots >|z_{k+l}|>0$ 
 to an analytic function of the form (\ref{k+l-prod-twted-twt-fd-2}),
except that there are no factors of the forms $(l_{p}(z_{1}-z))^{n_{1}}, \dots, 
(l_{p}(z_{k+l}-z))^{n_{k+l}}$. 
But 
\begin{equation}\label{psi-gen-locality-19}
\langle w',
\phi_{W}^{i_{1};p}(z_{1})\cdots \phi_{W}^{i_{k};p}(z_{k}) \psi_{W}^{a;p}(z)
\phi_{W}^{i_{k+1};p}(z_{k+1}) \cdots \phi_{W}^{i_{k+l};p}(z_{k+l})\rangle
\end{equation}
is a linear combination of series of the form 
\begin{align*}
&\langle w', e^{l_{p}(z_{1}-z)\mathcal{N}_{g}}
\phi_{W}^{i_{1};p}(z_{1})e^{-l_{p}(z_{1}-z)\mathcal{N}_{g}}
\cdots e^{l_{p}(z_{k}-z)\mathcal{N}_{g}}\phi_{W}^{i_{k};p}(z_{k})
e^{l_{p}(z_{k}-z)\mathcal{N}_{g}}\cdot\nn
&\quad\quad\quad\quad\cdot
\psi_{W}^{a;p}(z)e^{l_{p}(z-z_{k+1}\mathcal{N}_{g}}e^{\pi i\mathcal{N}_{g}}
\phi_{W}^{i_{k+1};p}(z_{k+1})e^{-\pi i\mathcal{N}_{g}}e^{-l_{p}(z-z_{k+1}\mathcal{N}_{g}}
\cdot\nn
&\quad\quad\quad\quad
\cdots e^{l_{p}(z-z_{k+l}\mathcal{N}_{g}}e^{\pi i\mathcal{N}_{g}}
\phi_{W}^{i_{k+l};p}(z_{k+l})
e^{-\pi i\mathcal{N}_{g}}e^{-l_{p}(z-z_{k+l}\mathcal{N}_{g}}v\rangle
\end{align*}
with nonnegative powers of $l_{p}(z_{1}-z), \dots, l_{p}(z_{k+l}-z)$ as coefficients. 
Thus (\ref{psi-gen-locality-19}) is absolutely convergent in the region 
given by $|z_{1}|>\cdots>|z_{k}|>|z|>
|z_{k+1}|>\cdots >|z_{k+l}|>0$, $|\arg (z_{i}-z)-\arg z_{i}|<\frac{\pi}{2}$ for 
$i=1, \dots, k$ and $|\arg (z_{i}-z)-\arg z|<\frac{\pi}{2}$ for 
$i=k+1, \dots, k+l$ to an analytic function of the form (\ref{k+l-prod-twted-twt-fd-2}). 
The remaining parts of Property 6 in Proposition \ref{basic-properties} follows
immediately from the proof above.
\epfv

We have the following most general commutativity which follows immediately from Property 
\ref{property-15} in Theorem \ref{psi-gen-locality}, Property 
\ref{property-14} in Theorem \ref{phi-locality} and Property 
\ref{property-16} in Theorem \ref{psi-gen-locality}:

\begin{cor}
Assume that Properties \ref{property-1}--\ref{property-7} in Assumption \ref{basic-properties}
hold. 
Then for $v\in V$, $w'\in V'$ and $\sigma\in S_{k}$,
$$
F^{p}(\langle w', \varphi_{1}(z_{1})\cdots \varphi_{k}(z_{k})v\rangle)
=\pm F^{p}(\langle w', \varphi_{\sigma(1)}(z_{\sigma(1)})\cdots
\varphi_{\sigma(k)}(z_{\sigma(k)})v\rangle),
$$
where one of $\varphi_{l}$ for $l=1, \dots, k$ is $\psi^{a;p}$ for some $a\in A$ and 
the others are in $\{\phi^{i;p}\;|\;i\in I\}$ and 
the sign $\pm$ is uniquely determined by $\sigma$
and $|\varphi_{1}|, \dots, |\varphi_{k}|$.
\end{cor}

\renewcommand{\theequation}{\thesection.\arabic{equation}}
\renewcommand{\thethm}{\thesection.\arabic{thm}}
\setcounter{equation}{0}
\setcounter{thm}{0}
\section{A construction theorem}

In this section, we construct a $g$-twisted $V$-module from 
the data in Data \ref{data}
satisfying Assumption \ref{basic-properties}.

First we need to define a twisted vertex operator map. 
Since we shall define the branches labeled by $p\in \Z$ of 
the twisted vertex operator map instead of defining the formal variable twisted 
vertex operator map, we need to show that these branches   labeled by $p\in \Z$
determine a formal variable map uniquely. 
So we first prove the following lemma:

\begin{lemma}\label{expansion}
Let $\phi$ be a multivalued analytic  map with preferred branch $\phi^{0}$
from $\C^{\times}$
to $\hom(W, \overline{W})$ and  $\phi^{p}$ for $p\in \Z$ are  labeled branches
of $\phi$.  If there exists $\wt \phi\in \frac{1}{2}\Z$ such that
$$[L^{g}_{W}(0), \phi^{p}(z)]=z\frac{d}{dz}\phi^{p}(z)+(\wt \phi)\phi^{p}(z),$$
then we have an expansion
\begin{equation}\label{expansion-0}
\phi^{p}(z)=\sum_{n\in \C}\sum_{k\in \N}
\phi_{n, k}e^{(-n-1)l_{p}(z)}(l_{p}(z))^{k}
\end{equation}
where  $\phi_{n, k}\in \hom(W, W)$ for $n\in \C$ and $k\in \N$
is homogeneous of weight $\wt \phi-n-1$. In particular, 
we obtain a formal series
$$\phi(x)=\sum_{n\in \C}\sum_{k\in \N}\phi_{n, k}x^{-n-1}(\log x)^{k}\in W\{x\}[\log x].$$
Moreover, for $w\in W$, there exist $N_{w}\in \R$ and $K_{n, w}\in \N$ for $n\in \C$ such that 
\begin{equation}\label{expansion-0.3}
\phi^{p}(z)w=\sum_{\Re{(n)}\ge N_{w}}\sum_{k=0}^{K_{n, w}}
\phi_{n, k}we^{(-n-1)l_{p}(z)}(l_{p}(z))^{k}.
\end{equation}
and for $w'\in w'$, there exist $N_{w'}\in \R$ such that 
\begin{equation}\label{expansion-0.7}
\langle w', \phi(z)\cdot\rangle=\sum_{\Re{(n)}\le N_{w'}}\sum_{k\in \N}
\langle w', \phi_{n, k}\cdot\rangle e^{(-n-1)l_{p}(z)}(l_{p}(z))^{k}.
\end{equation}
\end{lemma}
\pf
First by Property \ref{property-1} in Assumption \ref{basic-properties},  
$L_{W}(0)$ can be decomposed as the sum of its semisimple 
part $L_{W}(0)_{S}$ and its nilpotent part $L_{W}(0)_{N}$.
From the commutator formula for $L_{W}(0)$ and $\phi^{p}(z)$, we see that for
$c\in \C$,
$$e^{cL^{g}_{W}(0)}\phi^{p}(z)e^{-cL^{g}_{W}(0)}=e^{c(\swt \phi)}\phi^{p}(e^{c}z).$$
In particular, taking $c=-l_{p}(z)$, we obtain
\begin{align*}
\phi^{p}(z)&=e^{-l_{p}(z)(\swt \phi)}e^{l_{p}(z)L^{g}_{W}(0)}\phi^{p}(1)e^{-l_{p}(z)L^{g}_{W}(0)}\nn
&=e^{-l_{p}(z)(\swt \phi)}e^{l_{p}(z)L^{g}_{W}(0)_{S}}e^{l_{p}(z)L^{g}_{W}(0)_{N}}\phi^{p}(1)
e^{-l_{p}(z)L^{g}_{W}(0)_{N}}e^{-l_{p}(z)L^{g}_{W}(0)_{S}}.
\end{align*}
For $n\in \C$ and $w\in W_{[m]}$,  
$$e^{l_{p}(z)L^{g}_{W}(0)_{N}}\pi_{(\swt \phi)-n-1+m}\phi^{p}(1)
e^{-l_{p}(z)L^{g}_{W}(0)_{N}}w$$
is in fact a polynomial
in $l_{p}(z)$ with coefficients in $W_{[(\swt \phi)-n-1+m]}$, where for $r\in \C$, we
use $\pi_{r}$ to denote the projection from $W$ to $W_{[r]}$. Let $K_{n, w}$ be the degree of 
this polynomial in $l_{p}(z)$ and let $\phi_{n, k}w\in W$ for $k\in \N$
be the coefficient of the $k$-th power of $l_{p}(z)$ in  
$$e^{l_{p}(z)L^{g}_{W}(0)_{N}}\pi_{(\swt \phi)-n-1+m}\phi^{p}(1)
e^{-l_{p}(z)L^{g}_{W}(0)_{N}}w.$$ 
We obtain $\phi_{n, k}\in \hom(W, W)$ of weight $\wt \phi-n-1$ such that
$$\sum_{k=0}^{K_{n, w}}
\phi_{n, k}w(l_{p}(z))^{k}=e^{l_{p}(z)L^{g}_{W}(0)_{N}}\pi_{(\swt \phi)-n-1+m}\phi^{p}(1)
e^{-l_{p}(z)L^{g}_{W}(0)_{N}}w.$$
Then  for $n\in \Z$ and $w\in W$,
\begin{align*}
e&^{-l_{p}(z)(\swt \phi)}e^{l_{p}(z)L^{g}_{W}(0)_{S}}\sum_{k=0}^{K_{n, w}}\phi_{n, k}
e^{-l_{p}(z)L^{g}_{W}(0)_{S}}w (l_{p}(z))^{k}\nn
&\quad  =e^{l_{p}(z)L^{g}_{W}(0)_{S}}e^{l_{p}(z)L^{g}_{W}(0)_{N}}\pi_{(\swt \phi)-n-1+m}\phi^{p}(1)
e^{-l_{p}(z)L^{g}_{W}(0)_{N}}e^{-l_{p}(z)L^{g}_{W}(0)_{S}}w .
\end{align*}
Thus for $w\in W$,
\begin{align}\label{expansion-1}
\phi^{p}(z)w&=e^{-l_{p}(z)(\swt \phi)}e^{l_{p}(z)L^{g}_{W}(0)_{S}}e^{l_{p}(z)L^{g}_{W}(0)_{N}}\phi^{p}(1)
e^{-l_{p}(z)L^{g}_{W}(0)_{N}}e^{-l_{p}(z)L^{g}_{W}(0)_{S}}w\nn
&=\sum_{n\in \C}e^{-l_{p}(z)(\swt \phi)}e^{l_{p}(z)L^{g}_{W}(0)_{S}}e^{l_{p}(z)L^{g}_{W}(0)_{N}}\pi_{(\swt \phi)-n-1+m}\phi^{p}(1)
e^{-l_{p}(z)L^{g}_{W}(0)_{N}}e^{-l_{p}(z)L^{g}_{W}(0)_{S}}w\nn
&=\sum_{n\in \C}\sum_{k=0}^{K_{n, w}}
e^{-l_{p}(z)(\swt \phi)}e^{-l_{p}(z)(\swt \phi)} e^{l_{p}(z)L^{g}_{W}(0)_{S}}\phi_{n, k}
e^{-l_{p}(z)L^{g}_{W}(0)_{S}}w (l_{p}(z))^{k}\nn
&=\sum_{n\in \Z}\sum_{k=0}^{K_{n, w}}\phi_{n, k}w e^{(-n-1)l_{p}(z)} (l_{p}(z))^{k}.
\end{align}
Since $W$ is lower bounded with respect to the weights and the weight of $\phi_{n, k}$
is $\wt \phi-n-1$, we obtain (\ref{expansion-0.3}) from (\ref{expansion-1}) and we also have 
(\ref{expansion-0.7}).
Since nonhomogeneous elements of $W$ are finite sums of homogeneous elements of
$W$, (\ref{expansion-0}), (\ref{expansion-0.3})   and (\ref{expansion-0.7})
also holds for general $w\in W$. 
\epfv

The vertex operator map
we want to define is a linear map
\begin{align}\label{tw-vom}
Y^{g}_{W}: V\otimes W &\to W\{x\}[\log x],\nn
u\otimes w&\mapsto Y^{g}_{W}(u, x)w.
\end{align}
Such a map gives a multivalued analytic map (denoted using the same notation)
\begin{align*}
Y^{g}_{W}: \C^{\times}&\to \hom(V\otimes W, \overline{W}),\nn
z&\mapsto Y^{g}_{W}(\cdot, z)\cdot: u\otimes w\mapsto Y^{g}_{W}(u, z)w
\end{align*}
with labeled branches 
\begin{align*}
(Y^{g}_{W})^{p}: \C^{\times}&\to \hom(V\otimes W, \overline{W}),\nn
z&\mapsto (Y^{g}_{W})^{p}(\cdot, z)\cdot: u\otimes w\mapsto (Y^{g}_{W})^{p}(u, z)w
\end{align*}
for $p\in \Z$. Conversely, by Lemma \ref{expansion}
such a multivalued analytic map with
labeled branches  also determines a linear map of the form (\ref{tw-vom}).
Thus to define a twisted vertex operator map, we need only 
define $(Y^{g}_{W})^{p}$.

We first give the motivation of our definition. The idea is in fact the same as 
in \cite{H-2-const}. 
We define
$(Y^{g}_{W})^{p}(\phi_{-1}^{i}\one, z)w=\phi^{i;p}_{W}(z)w$ for $p\in \Z$, 
$i\in I$ and $w\in W$.
The vertex operator map should satisfy the duality property.
In particular, we should have
\begin{equation}\label{moti-tw-vom}
F^{p}(\langle w', Y^{g}_{W}(\phi^{i_{1}}(\xi_{1})\cdots \phi^{i_{k}}(\xi_{k})\one, z)w\rangle)
=F^{p}(\langle w', \phi_{W}^{i_{1}}(\xi_{1}+z)\cdots \phi_{W}^{i_{k}}(\xi_{k}+z)w\rangle)
\end{equation}
for $i_{1}, \dots, i_{k}\in I$, $w\in W$ and $w'\in W'$.
Note that $\phi^{i_{1}}(\xi_{1}), \dots, \phi^{i_{k}}(\xi_{k})$ are single-valued
analytic functions in $\xi_{1}, \dots, \xi_{k}$, respectively. Also by Property \ref{property-13}
in Theorem \ref{phi-locality}, the right-hand side of 
(\ref{moti-tw-vom}) is a single-valued analytic function of $\xi_{1}, \dots, \xi_{k}$
when $\xi_{i}+z\ne 0$ for $i=1, \dots, k$, $\xi_{i}\ne \xi_{j}$ 
for $i, j=1, \dots, k$ and $i\ne j$.

Motivated by (\ref{moti-tw-vom}), we define the vertex operator map as follows:
For $w'\in W'$, $w\in W$, $i_{1}, \dots, i_{k}\in I$, $m_{1}, \dots, m_{k}\in \Z$,
we define $(Y^{g}_{W})^{p}$ by
\begin{align}\label{vo}
\langle w', &(Y^{g}_{W})^{p}(\phi^{i_{1}}_{m_{1}}\cdots \phi^{i_{k}}_{m_{k}}\one, z)
w\rangle\nn
&=\res_{\xi_{1}=0}\cdots\res_{\xi_{k}=0}
\xi_{1}^{m_{1}}\cdots\xi_{k}^{m_{k}}
F^{p}(\langle w', \phi_{W}^{i_{1}}(\xi_{1}+z)\cdots \phi_{W}^{i_{k}}(\xi_{k}+z)
w\rangle).
\end{align}

Since there might be relations among elements of the form
$\phi^{i_{1}}_{m_{1}}\cdots \phi^{i_{k}}_{m_{k}}\one$,
we first have to show that the definition above indeed gives a well-defined map
from $\C^{\times}$ to $\hom(V\otimes W, \overline{W})$. Let $\phi^{0}$ be the map
from $\C^{\times}$ to $\hom(V, \overline{V})$ given by $\phi^{0}(z)=1_{V}$.
Let $\wt \phi^{0}=0$. Then Properties \ref{property-1}--\ref{property-7}
 in Assumption \ref{basic-properties} 
and Properties \ref{property-8}--\ref{property-12} in Proposition \ref{properties} still hold for
$\phi^{i}$, $i\in \tilde{I}=I\cup \{0\}$.
Then any relation among
such elements can always be written as
$$\sum_{\mu=1}^{M}\lambda_{\mu}\phi^{i^{\mu}_{1}}_{m^{\mu}_{1}}
\cdots \phi^{i^{\mu}_{k}}_{m^{\mu}_{k}}\one=0$$
for some $k\in \Z_{+}$, 
$ i^{\mu}_{j}\in \tilde{I}$ and $m^{\mu}_{j}\in \Z$ for
$\mu=1, \dots, M$, $j=1, \dots, k$, where $\phi^{i^{\mu}_{1}}_{m^{\mu}_{1}}
\cdots \phi^{i^{\mu}_{k}}_{m^{\mu}_{k}}\one$ for $\mu=1, \dots, M$ either
all belong to $V^{0}$ or all belong to $V^{1}$, that is,
$|\phi_{m_{1}^{\mu}}^{i_{1}^{\mu}}|+\cdots 
+|\phi_{m_{k}^{\mu}}^{i_{k}^{\mu}}|$ for 
$\mu=1, \dots, M$ are either all even or are all odd. In particular,
the parities of $|\phi_{m_{1}^{\mu}}^{i_{1}^{\mu}}|+\cdots 
+|\phi_{m_{k}^{\mu}}^{i_{k}^{\mu}}|$ are independent of $\mu$.
Since the parity  $|\phi_{m_{r}^{\mu}}^{i_{r}^{\mu}}|$ for 
$\mu=1, \dots, M$, $r=1, \dots, k$ are equal to the parity 
$|\phi^{i_{r}^{\mu}}|$, we see that the parities of 
$|\phi^{i_{1}^{\mu}}|+\cdots 
+|\phi^{i_{k}^{\mu}}|$ are independent of $\mu$.

\begin{lemma}\label{well-defined}
If $$\sum_{\mu=1}^{M}\lambda_{\mu}\phi^{i^{\mu}_{1}}_{m^{\mu}_{1}}
\cdots \phi^{i^{\mu}_{k}}_{m^{\mu}_{k}}\one=0,$$
then
\begin{equation}\label{well-defined-0}
\sum_{\mu=1}^{M}\lambda_{\mu}\res_{\xi_{1}=0}\cdots\res_{\xi_{k}=0}
\xi_{1}^{m^{\mu}_{1}}\cdots\xi_{k}^{m^{\mu}_{k}}F^{p}(\langle w', 
\phi_{W}^{i_{1}^{\mu}}(\xi_{1}+z)\cdots \phi_{W}^{i_{k}^{\mu}}(\xi_{k}+z)
w\rangle)=0
\end{equation}
for $w\in W$ and $w'\in W'$.
\end{lemma}
\pf
Since $V$ is generated by $\phi^{i}(x)$ for $i\in I$, 
by Property \ref{property-4} in Assumption \ref{basic-properties}, we can take
$$w=(\phi_{W})^{j_{1}}_{n_{1}, q_{1}}\cdots (\phi_{W})^{j_{l}}_{n_{l}, q_{l}}
(\psi_{W}^{a})_{n, q}
\phi^{j_{l+1}}_{n_{l+1}}
\cdots \phi^{j_{m}}_{n_{m}}\one.$$
Since this element
 is a coefficient of 
$$\phi_{W}^{j_{1}}(\zeta_{1})\cdots \phi_{W}^{j_{l}}(\zeta_{l})\psi_{W}^{a}(\zeta)
\phi^{j_{l+1}}(\zeta_{l+1})
\cdots \phi^{j_{m}}(\zeta_{m})\one,$$
we first prove 
\begin{align}\label{well-defined-1}
\sum_{\mu=1}^{M}&\lambda_{\mu}\res_{\xi_{1}=0}\cdots\res_{\xi_{k}=0}
\xi_{1}^{m^{\mu}_{1}}\cdots\xi_{k}^{m^{\mu}_{k}}\cdot\nn
& \cdot F^{p}(\langle w', 
\phi_{W}^{i_{1}^{\mu}}(\xi_{1}+z)\cdots \phi_{W}^{i_{k}^{\mu}}(\xi_{k}+z)
\phi_{W}^{j_{1}}(\zeta_{1})\cdots \phi_{W}^{j_{l}}(\zeta_{l})\psi_{W}^{a}(\zeta)
\phi^{j_{l+1}}(\zeta_{l+1})
\cdots \phi^{j_{m}}(\zeta_{m})\one\rangle)
=0.
\end{align}
We have
\begin{align*}
&\res_{\xi_{1}=0}\cdots\res_{\xi_{k}=0}
\xi_{1}^{m^{\mu}_{1}}\cdots\xi_{k}^{m^{\mu}_{k}}\cdot\nn
&\quad \quad \cdot F^{p}(\langle w', 
\phi_{W}^{i_{1}^{\mu}}(\xi_{1}+z)\cdots \phi_{W}^{i_{k}^{\mu}}(\xi_{k}+z)
\phi_{W}^{j_{1}}(\zeta_{1})\cdots \phi_{W}^{j_{l}}(\zeta_{l})\psi_{W}^{a}(\zeta)
\phi^{j_{l+1}}(\zeta_{l+1})
\cdots \phi^{j_{m}}(\zeta_{m})\one\rangle)\nn
&\quad=\prod_{r=1}^{k}\prod_{s=1}^{m}(-1)^{|\phi^{i^{\mu}_{r}}||\phi^{j_{s}}|} 
\prod_{r=1}^{k}(-1)^{|\phi^{i^{\mu}_{r}}||\psi_{W}^{a}|}\res_{\xi_{1}=0}\cdots\res_{\xi_{k}=0}
\xi_{1}^{m^{\mu}_{1}}\cdots\xi_{k}^{m^{\mu}_{k}}\cdot\nn
&\quad\quad\cdot F^{p}(\langle w',
\phi_{W}^{j_{1}}(\zeta_{1})\cdots \phi_{W}^{j_{l}}(\zeta_{l})\psi_{W}^{a}(\zeta)
\phi^{j_{l+1}}(\zeta_{l+1})
\cdots \phi^{j_{m}}(\zeta_{m})
\phi^{i^{\mu}_{1}}(\xi_{1}+z)\cdots \phi^{i^{\mu}_{k}}(\xi_{k}+z)\one\rangle)\nn
&\quad=\prod_{s=1}^{m}(-1)^{\left(\left|\phi^{i_{1}^{\mu}}\right|+\cdots 
+\left|\phi^{i_{k}^{\mu}}\right|\right)|\phi^{j_{s}}|} 
(-1)^{\left(\left|\phi^{i_{1}^{\mu}}\right|+\cdots 
+\left|\phi^{i_{k}^{\mu}}\right|\right)|\psi_{W}^{a}|}\res_{\xi_{1}=0}\cdots\res_{\xi_{k}=0}
\xi_{1}^{m^{\mu}_{1}}\cdots\xi_{k}^{m^{\mu}_{k}}\cdot\nn
&\quad\quad\cdot F^{p}(\langle e^{zL_{W}(-1)'}w',
\phi_{W}^{j_{1}}(\zeta_{1}-z)\cdots \phi_{W}^{j_{l}}(\zeta_{l}-z)
\psi_{W}^{a}(\zeta-z)\cdot\nn
&\quad\quad\quad\quad\quad\quad\quad\quad\quad\quad
\cdot \phi^{j_{l+1}}(\zeta_{l+1}-z)
\cdots \phi^{j_{m}}(\zeta_{m}-z)\phi^{i^{\mu}_{1}}(\xi_{1})\cdots \phi^{i^{\mu}_{k}}(\xi_{k})\one\rangle)\nn
&\quad=\prod_{s=1}^{m}(-1)^{\left(\left|\phi^{i_{1}^{\mu}}\right|+\cdots 
+\left|\phi^{i_{k}^{\mu}}\right|\right)|\phi^{j_{s}}|} 
(-1)^{\left(\left|\phi^{i_{1}^{\mu}}\right|+\cdots 
+\left|\phi^{i_{k}^{\mu}}\right|\right)|\psi_{W}^{a}|}
\cdot\nn
&\quad\quad \cdot  F^{p}(\langle e^{zL_{W}(-1)'}w',
\phi_{W}^{j_{1}}(\zeta_{1}-z)\cdots \phi_{W}^{j_{l}}(\zeta_{l}-z)\psi_{W}^{a}(\zeta-z)\cdot\nn
&\quad\quad\quad\quad\quad\quad\quad\quad\quad\quad
\cdot 
\phi^{j_{l+1}}(\zeta_{l+1}-z)
\cdots \phi^{j_{m}}(\zeta_{m}-z)
\phi^{i^{p}_{1}}_{m^{\mu}_{1}}\cdots \phi^{i^{p}_{k}}_{m^{\mu}_{k}}\one\rangle).
\end{align*}
Recalling that  $|\phi^{i_{1}^{\mu}}|+\cdots 
+|\phi^{i_{k}^{\mu}}|$ are independent of $\mu$, we obtain
\begin{align*}
&\sum_{\mu=1}^{M}\lambda_{\mu}\res_{\xi_{1}=0}\cdots\res_{\xi_{k}=0}
\xi_{1}^{m^{\mu}_{1}}\cdots\xi_{k}^{m^{\mu}_{k}}\cdot\nn
&\;\quad \cdot F^{p}(\langle w', 
\phi_{W}^{i_{1}^{\mu}}(\xi_{1}+z)\cdots \phi_{W}^{i_{k}^{\mu}}(\xi_{k}+z)
\phi_{W}^{j_{1}}(\zeta_{1})\cdots \phi_{W}^{j_{l}}(\zeta_{l})\psi_{W}^{a}(\zeta)
\phi^{j_{l+1}}(\zeta_{l+1})
\cdots \phi^{j_{m}}(\zeta_{m})\one\rangle)\nn
&\;=\sum_{\mu=1}^{M}\prod_{s=1}^{m}(-1)^{\left(\left|\phi^{i_{1}^{\mu}}\right|+\cdots 
+\left|\phi^{i_{k}^{\mu}}\right|\right)|\phi^{j_{s}}|} 
(-1)^{\left(\left|\phi^{i_{1}^{\mu}}\right|+\cdots 
+\left|\phi^{i_{k}^{\mu}}\right|\right)|\psi_{W}^{a}|}
\cdot\nn
&\;\quad \cdot  F^{p}(\langle e^{zL_{W}(-1)'}w',
\phi_{W}^{j_{1}}(\zeta_{1}-z)\cdots \phi_{W}^{j_{l}}(\zeta_{l}-z)\psi_{W}^{a}(\zeta-z)\cdot\nn
&\quad\quad\quad\quad\quad\quad\quad\quad\quad\quad
\cdot 
\phi^{j_{l+1}}(\zeta_{l+1}-z)
\cdots \phi^{j_{m}}(\zeta_{m}-z)
\phi^{i^{p}_{1}}_{m^{\mu}_{1}}\cdots \phi^{i^{p}_{k}}_{m^{\mu}_{k}}\one\rangle)\nn
&\;=\prod_{s=1}^{m}(-1)^{\left(\left|\phi^{i_{1}^{\mu}}\right|+\cdots 
+\left|\phi^{i_{k}^{\mu}}\right|\right)|\phi^{j_{s}}|} 
(-1)^{\left(\left|\phi^{i_{1}^{\mu}}\right|+\cdots 
+\left|\phi^{i_{k}^{\mu}}\right|\right)|\psi_{W}^{a}|}
 \cdot
\end{align*}
\begin{align*}
&\quad \quad\cdot
F^{p}\Biggl(\Biggl\langle e^{zL_{W}(-1)'}w',
\phi_{W}^{j_{1}}(\zeta_{1}-z)\cdots \phi_{W}^{j_{l}}(\zeta_{l}-z)
\psi_{W}^{a}(\zeta-z)\cdot\nn
&\quad\quad\quad\quad\quad\quad\quad\quad\quad\quad\quad
\cdot 
\phi^{j_{l+1}}(\zeta_{l+1}-z)
\cdots \phi^{j_{m}}(\zeta_{m}-z)
\left(\sum_{\mu=1}^{M}\lambda_{\mu}\phi^{i^{\mu}_{1}}_{m^{\mu}_{1}}\cdots
\phi^{i^{\mu}_{k}}_{m^{\mu}_{k}}\one\right)\Biggr\rangle\Biggr)
\quad\quad\;\nn
&\quad=0,
\end{align*}
proving (\ref{well-defined-1}).

For any fixed $z, \xi_{1}, \dots, 
\xi_{l}, \xi$, the left-hand side of 
(\ref{well-defined-1}) can be expanded in the region $|z|, |\xi_{1}|, \dots, 
|\xi_{l}|, |\xi|>|\zeta_{l+1}|>\cdots>|\zeta_{m}|$ as 
\begin{align}\label{well-defined-2}
\sum_{n_{1}, \dots, n_{m}\in \Z}\sum_{\mu=1}^{M}&\lambda_{\mu}
\res_{\xi_{1}=0}\cdots\res_{\xi_{k}=0}
\xi_{1}^{m^{\mu}_{1}}\cdots\xi_{k}^{m^{\mu}_{k}}\cdot\nn
&\quad \cdot F^{p}(\langle w', 
\phi_{W}^{i_{1}^{\mu}}(\xi_{1}+z)\cdots \phi_{W}^{i_{k}^{\mu}}(\xi_{k}+z)
\phi_{W}^{j_{1}}(\zeta_{1})\cdots \phi_{W}^{j_{l}}(\zeta_{l})\psi_{W}^{a}(\xi)
\phi^{j_{l+1}}_{n_{l+1}}
\cdots \phi^{j_{m}}_{n_{m}}\one\rangle)\cdot\nn
&\quad \cdot\zeta_{l+1}^{-n_{l+1}-1}\cdots \zeta_{m}^{-n_{m}-1}.
\end{align}
By (\ref{well-defined-1}), the Laurent series  (\ref{well-defined-2}) in $\zeta_{l+1}, \dots, \zeta_{m}$
is also $0$ and thus its coefficients are all $0$. 
So we obtain 
\begin{align}\label{well-defined-2.5}
\sum_{\mu=1}^{M}&\lambda_{\mu}
\res_{\xi_{1}=0}\cdots\res_{\xi_{k}=0}
\xi_{1}^{m^{\mu}_{1}}\cdots\xi_{k}^{m^{\mu}_{k}}\cdot\nn
&\quad \cdot F^{p}(\langle w', 
\phi_{W}^{i_{1}^{\mu}}(\xi_{1}+z)\cdots \phi_{W}^{i_{k}^{\mu}}(\xi_{k}+z)
\phi_{W}^{j_{1}}(\zeta_{1})\cdots \phi_{W}^{j_{l}}(\zeta_{l})\psi_{W}^{a}(\xi)
\phi^{j_{l+1}}_{n_{l+1}}
\cdots \phi^{j_{m}}_{n_{m}}\one\rangle)\nn
&=0.
\end{align}
By Assumption \ref{algebra}, $\phi^{j_{l+1}}_{n_{l+1}}
\cdots \phi^{j_{m}}_{n_{m}}\one$ is a generalized eigenvector for $g$ with the 
eigenvalue $e^{2\pi i(\alpha^{j_{l+1}}+\cdots +\alpha^{j_{m}})}$. 
Then there exists $N, K\in \N$ such that 
the left-hand side of (\ref{well-defined-2.5}) can be expanded when $\xi$ is sufficiently small
but not $0$ as 
\begin{align}\label{well-defined-2.6}
&\sum_{q=0}^{K}\sum_{n\in \alpha^{j_{l+1}}+\cdots +\alpha^{j_{m}}+N-\N}
\sum_{\mu=1}^{M}\lambda_{\mu}
\res_{\xi_{1}=0}\cdots\res_{\xi_{k}=0}
\xi_{1}^{m^{\mu}_{1}}\cdots\xi_{k}^{m^{\mu}_{k}}\cdot\nn
&\quad \cdot F^{p}(\langle w', 
\phi_{W}^{i_{1}^{\mu}}(\xi_{1}+z)\cdots \phi_{W}^{i_{k}^{\mu}}(\xi_{k}+z)
\phi_{W}^{j_{1}}(\zeta_{1})\cdots \phi_{W}^{j_{l}}(\zeta_{l})(\psi_{W}^{a})_{n, q}
\phi^{j_{l+1}}_{n_{l+1}}
\cdots \phi^{j_{m}}_{n_{m}}\one\rangle)\xi^{-n-1}l_{p}(\xi)^{q}.
\end{align}
By (\ref{well-defined-2.5}) and the fact that $(-\alpha^{j_{l+1}}-\cdots -
\alpha^{j_{m}}-N+\N)\times \{1, \dots, K\}$
is a unique expansion set (see Proposition 2.1 in \cite{H-zhang-correction}),
the expansion coefficients of (\ref{well-defined-2.6}) must be $0$,
that is,
\begin{align}\label{well-defined-3}
\sum_{\mu=1}^{M}&\lambda_{\mu}
\res_{\xi_{1}=0}\cdots\res_{\xi_{k}=0}
\xi_{1}^{m^{\mu}_{1}}\cdots\xi_{k}^{m^{\mu}_{k}}\cdot\nn
&\quad \cdot F^{p}(\langle w', 
\phi_{W}^{i_{1}^{\mu}}(\xi_{1}+z)\cdots \phi_{W}^{i_{k}^{\mu}}(\xi_{k}+z)
\phi_{W}^{j_{1}}(\zeta_{1})\cdots \phi_{W}^{j_{l}}(\zeta_{l})(\psi_{W}^{a})_{n, q}
\phi^{j_{l+1}}_{n_{l+1}}
\cdots \phi^{j_{m}}_{n_{m}}\one\rangle)\nn
&=0.
\end{align}

By Data \ref{data}, there exists 
$N_{l}, K_{l}\in \N$ such that
the left-hand side of 
(\ref{well-defined-3}) can be expanded when $|\zeta_{l}|$ is sufficiently small but not $0$ as 
\begin{align}\label{well-defined-4}
&\sum_{q_{l}=0}^{K_{l}}\sum_{n_{l}\in \alpha^{j_{l}}+N_{l}-\N}
\sum_{\mu=1}^{M}\lambda_{\mu}
\res_{\xi_{1}=0}\cdots\res_{\xi_{k}=0}
\xi_{1}^{m^{\mu}_{1}}\cdots\xi_{k}^{m^{\mu}_{k}}\cdot\nn
& \quad \cdot F^{p}(\langle w', 
\phi_{W}^{i_{1}^{\mu}}(\xi_{1}+z)\cdots \phi_{W}^{i_{k}^{\mu}}(\xi_{k}+z)
\phi_{W}^{j_{1}}(\zeta_{1})\cdots \phi_{W}^{j_{l-1}}(\zeta_{l-1})
(\phi_{W}^{j_{l}})_{n_{l}, q_{l}}
(\psi_{W}^{a})_{n, q}
\phi^{j_{l+1}}_{n_{l+1}}
\cdots \phi^{j_{m}}_{n_{m}}\one\rangle)\cdot \nn
&  \quad\cdot e^{(-n_{l}-1)l_{p}(\zeta_{l})}(l_{p}(\zeta_{l}))^{q_{l}}.
\end{align}
By (\ref{well-defined-3}),
(\ref{well-defined-4}) and the fact that  
$(-\alpha^{j_{l}}-N_{l}+\N)\times \{1, \dots, K_{l}\}$
is a unique expansion set (again see Proposition 2.1 in \cite{H-zhang-correction}), 
the expansion coefficients of (\ref{well-defined-4})
must be $0$, that is,
\begin{align*}
\sum_{\mu=1}^{M}&\lambda_{\mu}
\res_{\xi_{1}=0}\cdots\res_{\xi_{k}=0}
\xi_{1}^{m^{\mu}_{1}}\cdots\xi_{k}^{m^{\mu}_{k}}\cdot\nn
& \quad \cdot 
F^{p}(\langle w', 
\phi_{W}^{i_{1}^{\mu}}(\xi_{1}+z)\cdots \phi_{W}^{i_{k}^{\mu}}(\xi_{k}+z)\cdot\nn
& \quad\quad\quad\quad\quad\quad\quad \cdot
\phi_{W}^{j_{1}}(\zeta_{1})\cdots \phi_{W}^{j_{l-1}}(\zeta_{l-1})
(\phi_{W}^{j_{l}})_{n_{l}, q_{l}}
(\psi_{W}^{a})_{n, q}
\phi^{j_{l+1}}_{n_{l+1}}
\cdots \phi^{j_{m}}_{n_{m}}\one\rangle)\nn
&\quad =0.
\end{align*}
Continuing this process repeatedly for $\phi_{W}^{j_{l-1}}(\zeta_{l-1}), \dots, 
\phi_{W}^{j_{1}}(\zeta_{1})$, we obtain (\ref{well-defined-0}).
\epfv

From this lemma, we see that $(Y^{g}_{W})^{p}$ and thus 
the vertex operator map $Y^{g}_{W}$ are well defined.

The following result is our construction theorem:

\begin{thm}\label{const-thm}
The pair $(W, Y^{g}_{W})$ is a lower-bounded generalized $g$-twisted
$V$-module generated 
by $(\psi_{W}^{a})_{n, k}v$ for $a\in A$, $n\in \alpha+\Z$,
$k\in \N$, $v\in V^{[\alpha]}$ and $\alpha\in P_{V}$.
Moreover, this is the unique lower-bounded generalized 
$g$-twisted $V$-module structure on $W$ generated by $(\psi_{W}^{a})_{n, k}v$ 
for $a\in A$, $n\in \alpha+\Z$,
$k\in \N$, $v\in V^{[\alpha]}$ and $\alpha\in P_{V}$
such that $Y_{W}(\phi^{i}_{-1}\one, z)=\phi_{W}^{i}(z)$
for $i\in I$.
\end{thm}
\pf
The proof of this theorem is similar to the proof of Theorem 3.5 in \cite{H-2-const} but is
more complicated because the twisted vertex operator map
is multivalued. We refer the reader to \cite{H-twist-vo} for the definition of 
lower-bounded generalized $g$-twisted $V$-module.

The identity property
follow from of the definition of $Y^{g}_{W}$.

Let $L^{g}_{W}(0)'$ be the adjoint operator of $L^{g}_{W}(0)$.
For $w'\in W'$, $w\in W$, $ i_{1}, \dots, i_{k}\in I$ and $n_{1}, \dots, n_{k}\in \Z$, $q_{1}, \dots, q_{k}\in \N$, $c\in \C$,
\begin{align*}
\langle&  w', e^{cL^{g}_{W}(0)}(Y^{g}_{W})^{p}(\phi^{i_{1}}_{n_{1}}\cdots \phi^{i_{k}}_{n_{k}}\one, z)e^{-cL^{g}_{W}(0)}w\rangle\nn
&=\langle  e^{cL^{g}_{W}(0)'}w', (Y^{g}_{W})^{p}(\phi^{i_{1}}_{n_{1}}\cdots \phi^{i_{k}}_{n_{k}}\one, z)e^{-cL^{g}_{W}(0)}w\rangle\nn
&=\res_{\xi_{1}=0}\cdots\res_{\xi_{k}=0}\xi_{1}^{n_{1}}\cdots \xi_{k}^{n_{k}}
F^{p}(\langle  e^{cL^{g}_{W}(0)'}w', \phi_{W}^{i_{1}}(\xi_{1}+z)\cdots \phi_{W}^{i_{k}}(\xi_{k}+z)e^{-cL^{g}_{W}(0)}w\rangle)\nn
&=\res_{\xi_{1}=0}\cdots\res_{\xi_{k}=0}\xi_{1}^{n_{1}}\cdots \xi_{k}^{n_{k}}
F^{p}(\langle  v', e^{cL^{g}_{W}(0)}\phi_{W}^{i_{1}}(\xi_{1}+z)\cdots \phi_{W}^{i_{k}}(\xi_{k}+z)e^{-cL^{g}_{W}(0)}w\rangle)\nn
&=\res_{\xi_{1}=0}\cdots\res_{\xi_{k}=0}\xi_{1}^{n_{1}}\cdots \xi_{k}^{n_{k}}a^{c(\swt \phi^{i_{1}}+\cdots
\swt \phi^{i_{k}})}
F^{p}(\langle  w', \phi_{W}^{i_{1}}(a\xi_{1}+az)\cdots \phi_{W}^{i_{k}}(a\xi_{k}+az)w\rangle)
\end{align*}
\begin{align*}
&=\res_{\zeta_{1}=0}\cdots\res_{\zeta_{k}=0}\zeta_{1}^{n_{1}}\cdots \zeta_{k}^{n_{k}}a^{\swt \phi^{i_{1}}+\cdots
\swt \phi^{i_{k}}-k-n_{1}-\cdots -n_{k}}\cdot\nn
&\quad\quad\quad\quad\quad\quad\quad\quad\quad\quad\quad\quad\quad\quad\cdot
F^{p}(\langle  w', \phi_{W}^{i_{1}}(\zeta_{1}+az)\cdots \phi_{W}^{i_{k}}(\zeta_{k}+az)w\rangle)\quad\quad\quad\quad\quad\nn
&=\langle  w', (Y^{g}_{W})^{p}(e^{cL_{V}(0)}\phi^{i_{1}}_{n_{1}}\cdots \phi^{i_{k}}_{n_{k}}\one , az)w\rangle).
\end{align*}
This formula is equivalent to the $L(0)$-commutator formula.

From Property \ref{property-2} in Assumption \ref{basic-properties}
and the definition of $(Y^{g}_{W})^{p}$, we obtain
the $L(-1)$-commutator formula
$$\frac{d}{dz}(Y^{g}_{W})^{p}(\phi^{i_{1}}_{n_{1}}\cdots \phi^{i_{k}}_{n_{k}}\one, z)
=[L^{g}_{W}(-1), (Y^{g}_{W})^{p}(\phi^{i_{1}}_{n_{1}}\cdots \phi^{i_{k}}_{n_{k}}\one, z)].$$

Let $\{e_{n}\}_{n\in \Z}$ be a homogeneous basis of $W$ and $\{e_{n}'\}_{n\in \Z}$ its dual basis in $W'$.
Then we have
\begin{align}\label{prod}
\langle w&', (Y^{g}_{W})^{p}(\phi^{i_{1}}_{n_{1}}\cdots \phi^{i_{k}}_{n_{k}}\one, z_{1})
(Y^{g}_{W})^{p}(\phi^{j_{1}}_{m_{1}}\cdots \phi^{j_{l}}_{m_{l}}\one, z_{2})w\rangle\nn
&=\sum_{n\in \Z}\langle w', (Y^{g}_{W})^{p}(\phi^{i_{1}}_{n_{1}}\cdots \phi^{i_{k}}_{n_{k}}\one, z_{1})e_{n}\rangle
\langle e'_{n}, (Y^{g}_{W})^{p}(\phi^{j_{1}}_{m_{1}}\cdots \phi^{j_{l}}_{m_{l}}\one, z_{2})w\rangle\nn
&=\sum_{n\in \Z}\res_{\zeta_{1}=0}\cdots\res_{\zeta_{k}=0}\zeta_{1}^{n_{1}}\cdots \zeta_{k}^{n_{k}}
\res_{\xi_{1}=0}\cdots\res_{\xi_{l}=0}\xi_{1}^{m_{1}}\cdots \xi_{l}^{m_{l}}\cdot\nn
&\quad\quad\;\cdot F^{p}(\langle w', \phi_{W}^{i_{1}}(\zeta_{1}+z_{1})\cdots
 \phi_{W}^{i_{k}}(\zeta_{k}+z_{1})e_{n}\rangle)
F^{p}(\langle e'_{n}, \phi_{W}^{j_{1}}(\xi_{1}+z_{2})\cdots \phi_{W}^{j_{l}}(\xi_{l}+z_{2})w\rangle)\nn
&=\res_{\zeta_{1}=0}\cdots\res_{\zeta_{k}=0}\zeta_{1}^{n_{1}}\cdots \zeta_{k}^{n_{k}}
\res_{\xi_{1}=0}\cdots\res_{\xi_{l}=0}\xi_{1}^{m_{1}}\cdots \xi_{l}^{m_{l}}\cdot\nn
&\quad\quad\;\cdot \sum_{n\in \Z}F^{p}(\langle w', \phi_{W}^{i_{1}}(\zeta_{1}+z_{1})\cdots
 \phi_{W}^{i_{k}}(\zeta_{k}+z_{1})e_{n}\rangle)
F^{p}(\langle e'_{n}, \phi_{W}^{j_{1}}(\xi_{1}+z_{2})\cdots \phi_{W}^{j_{l}}(\xi_{l}+z_{2})w\rangle).\nn
\end{align}

By Property \ref{property-13} in Theorem \ref{phi-locality}, when $|z_{1}|>\cdots >|z_{k+l}|>0$,
\begin{align}\label{prod-1}
\sum_{n\in \Z}&F^{p}(\langle w', \phi_{W}^{i_{1}}(z_{1})\cdots
 \phi_{W}^{i_{k}}(z_{k})e_{n}\rangle)
F^{p}(\langle e'_{n}, \phi_{W}^{j_{1}}(z_{k+1})\cdots \phi_{W}^{j_{l}}(z_{k+l})w\rangle)\nn
&=\sum_{n\in \Z}\langle w', \phi_{W}^{i_{1}}(z_{1})\cdots
 \phi_{W}^{i_{k}}(z_{k})e_{n}\rangle
\langle e'_{n}, \phi_{W}^{j_{1}}(z_{k+1})\cdots \phi_{W}^{j_{l}}(z_{k+l})w\rangle\nn
&=\langle w', \phi_{W}^{i_{1}}(z_{1})\cdots
 \phi_{W}^{i_{k}}(z_{k})\phi_{W}^{j_{1}}(z_{k+1})\cdots \phi_{W}^{j_{l}}(z_{k+l})w\rangle
\end{align}
is absolutely convergent to the analytic function
\begin{equation}\label{prod-2}
F^{p}(\langle w', \phi_{W}^{i_{1}}(z_{1})\cdots
 \phi_{W}^{i_{k}}(z_{k})\phi_{W}^{j_{1}}(z_{k+1})\cdots \phi_{W}^{j_{l}}(z_{k+l})w\rangle)
\end{equation}
in $z_{1}, \dots, z_{k+l}$. On the other hand, also by Property \ref{property-13} in 
Theorem \ref{phi-locality}, 
there is a unique expansion of this branch of a multivalued function
in the region $|z_{1}|, \dots, |z_{k}|>|z_{k+1}|, \dots, |z_{k+l}|>0$, $z_{i}\ne z_{j}$ for $i\ne j$,
$i, j=1, \dots, k$ and $i, j=k+1, \dots, k+l$ such that each term
is a product of two analytic functions of the same form, one in $z_{1}, \dots, z_{k}$ 
and the other in $z_{k+1}, \dots,
z_{k+l}$. Since the left-hand side of (\ref{prod-1}) is a series of the same form
and is absolutely convergent
in the region $|z_{1}|>\cdots >|z_{k+l}|>0$ to  (\ref{prod-2}), it must be absolutely convergent
in the larger region $|z_{1}|, \dots,|z_{k}|>|z_{k+1}|, \dots, |z_{k+l}|>0$,
$z_{i}\ne z_{j}$ for $i\ne j$,
$i, j=1, \dots, k$ and $i, j=k+1, \dots, k+l$ to (\ref{prod-2}).

Substituting $\zeta_{i}+z_{1}$ for $z_{i}$ for $i=1, \dots, k$ and $\xi_{j}+z_{2}$ for $z_{k+j}$
for $j=1, \dots, l$, we see that
$$ \sum_{n\in \Z}F^{p}(\langle w', \phi_{W}^{i_{1}}(\zeta_{1}+z_{1})\cdots
 \phi_{W}^{i_{k}}(\zeta_{k}+z_{1})e_{n}\rangle)
F^{p}(\langle e'_{n}, \phi_{W}^{j_{1}}(\xi_{1}+z_{2})\cdots \phi_{W}^{j_{l}}(\xi_{l}+z_{2})w\rangle)$$
is absolutely convergent to
$$F^{p}(\langle w', \phi_{W}^{i_{1}}(\zeta_{1}+z_{1})\cdots \phi_{W}^{i_{k}}(\zeta_{k}+z_{1})
\phi_{W}^{j_{1}}(\xi_{1}+z_{2})\cdots \phi_{W}^{j_{l}}(\xi_{l}+z_{2})w\rangle)$$
when $|\zeta_{1}+z_{1}|, \dots, |\zeta_{k}+z_{1}|>|\xi_{1}+z_{2}|, \dots, |\xi_{l}+z_{2}|>0$,
$\zeta_{i}\ne \zeta_{j}$  for $i, j=1, \dots, k$ and $\xi_{i}\ne \xi_{j}$ for $i, j=1, \dots, l$. When
$|z_{1}|>|z_{2}|>0$, we can always find sufficiently small neighborhood of $0$ such that when
$\zeta_{1}, \dots, \zeta_{k}, \xi_{1}, \dots, \xi_{l}$ are in this neighborhood,
$|\zeta_{1}+z_{1}|, \dots, |\zeta_{k}+z_{1}|>|\xi_{1}+z_{2}|, \dots, |\xi_{l}+z_{2}|>0$
holds. Thus we see that when $|z_{1}|>|z_{2}|>0$, the right-hand side of (\ref{prod}) is absolutely convergent to
\begin{align}\label{prod-3}
\res_{\zeta_{1}=0}&\cdots\res_{\zeta_{k}=0}\zeta_{1}^{n_{1}}\cdots \zeta_{k}^{n_{k}}
\res_{\xi_{1}=0}\cdots\res_{\xi_{l}=0}\xi_{1}^{m_{1}}\cdots \xi_{l}^{m_{l}}\cdot\nn
&\quad\quad\quad\cdot F^{p}(\langle w', \phi_{W}^{i_{1}}(\zeta_{1}+z_{1})\cdots \phi_{W}^{i_{k}}(\zeta_{k}+z_{1})
\phi_{W}^{j_{1}}(\xi_{1}+z_{2})\cdots \phi_{W}^{j_{l}}(\xi_{l}+z_{2})w\rangle).
\end{align}
From the explicit expression of 
$$F^{p}(\langle w', \phi_{W}^{i_{1}}(\zeta_{1}+z_{1})\cdots \phi_{W}^{i_{k}}(\zeta_{k}+z_{1})
\phi_{W}^{j_{1}}(\xi_{1}+z_{2})\cdots \phi_{W}^{j_{l}}(\xi_{l}+z_{2})w\rangle)$$
(see Property \ref{property-13} in Theorem \ref{phi-locality}),
 it is clear that (\ref{prod-3}) is an analytic function in $z_{1}$ and $z_{2}$ of the form  
\begin{equation}\label{correl-fns-p-branch}
\sum_{i, j, k, l = 0}^N
a_{ijkl}e^{m_il_p(z_1)}e^{n_jl_p(z_2)}l_p(z_1)^kl_p(z_2)^l(z_1 -
z_2)^{-t}.
\end{equation}
In particular, the left-hand side of
(\ref{prod}), that is,
\begin{equation}\label{prod-4}
\langle w', (Y^{g}_{W})^{p}(\phi^{i_{1}}_{n_{1}}\cdots \phi^{i_{k}}_{n_{k}}\one, z_{1})
(Y^{g}_{W})^{p}(\phi^{j_{1}}_{m_{1}}\cdots \phi^{j_{l}}_{m_{l}}\one, z_{2})w\rangle,
\end{equation}
is absolutely convergent in the region $|z_{1}|>|z_{2}|>0$ to this analytic
function.

We have proved that the product of two vertex operators is convergent to an analytic 
function of the form (\ref{correl-fns-p-branch}), or equivalently, the corresponding branch
of a multivalued function with preferred branch of the form 
$$
f(z_1, z_2) = \sum_{i,
j, k, l = 0}^N a_{ijkl}z_1^{m_i}z_2^{n_j}({\rm log}z_1)^k({\rm
log}z_2)^l(z_1 - z_2)^{-t}.
$$
We are ready to
prove the commutativity.
The calculation above also shows that
\begin{equation}\label{prod-5}
\langle w',
(Y^{g}_{W})^{p}(\phi^{j_{1}}_{m_{1}}\cdots \phi^{j_{l}}_{m_{l}}\one, z_{2})
(Y^{g}_{W})^{p}(\phi^{i_{1}}_{n_{1}}\cdots \phi^{i_{k}}_{n_{k}}\one, z_{1})w\rangle
\end{equation}
is absolutely convergent to the rational function
\begin{align}\label{prod-6}
\res_{\xi_{1}=0}&\cdots\res_{\xi_{l}=0}\xi_{1}^{m_{1}}\cdots \xi_{l}^{m_{l}}
\res_{\zeta_{1}=0}\cdots\res_{\zeta_{k}=0}\zeta_{1}^{n_{1}}\cdots \zeta_{k}^{n_{k}}\cdot\nn
&\quad\quad\quad\quad \cdot F^{p}(\langle w',
\phi_{W}^{j_{1}}(\xi_{1}+z_{2})\cdots \phi_{W}^{j_{l}}(\xi_{l}+z_{2})
\phi_{W}^{i_{1}}(\zeta_{1}+z_{1})\cdots \phi_{W}^{i_{k}}(\zeta_{k}+z_{1})w\rangle),
\end{align}
in the regions $|z_{2}|>|z_{1}|>0$, respectively. By Property \ref{property-14}  in Theorem
\ref{phi-locality}, the analytic functions
(\ref{prod-3}) and (\ref{prod-6}) multiplied by
$$\prod_{r=1}^{k}\prod_{s=1}^{l}(-1)^{|\phi^{i_{r}}||\phi^{j_{s}}|}=(-1)^{|\phi^{i_{1}}_{n_{1}}\cdots \phi^{i_{k}}_{n_{k}}\one|
|\phi^{j_{1}}_{m_{1}}\cdots \phi^{j_{l}}_{m_{l}}\one|}$$
are equal.
Thus (\ref{prod-4}) and (\ref{prod-5}) multiplied by the sign
$(-1)^{|\phi^{i_{1}}_{n_{1}}\cdots \phi^{i_{k}}_{n_{k}}\one|
|\phi^{j_{1}}_{m_{1}}\cdots \phi^{j_{l}}_{m_{l}}\one|}$
are
absolutely convergent in the regions $|z_{1}|>|z_{2}|>0$ and $|z_{2}|>|z_{1}|>0$,
respectively, to a common analytic function of the form (\ref{correl-fns-p-branch}).

We now prove the associativity. For $i_{1}, \dots, i_{k}, j_{1}, \dots, j_{l}\in I$, $m_{1}, \dots, m_{l}\in \Z$,
$v\in V$ and $v'\in V'$,
using the expansion of $\phi^{i_{1}}(\xi_{1}), \dots, \phi^{i_{k}}(\xi_{k})$ and the
definition of $(Y^{g}_{W})^{p}$, we have
\begin{align}\label{phi-assoc}
\langle& w', (Y^{g}_{W})^{p}(\phi^{i_{1}}(z_{1})\cdots \phi^{i_{k}}(z_{k})
\phi^{j_{1}}_{m_{1}}\cdots \phi^{j_{l}}_{m_{l}}\one, z)w\rangle\nn
&=\sum_{p_{1}, \dots, p_{k}\in \Z}\langle w', (Y^{g}_{W})^{p}(\phi^{i_{1}}_{p_{1}}\cdots \phi^{i_{k}}_{p_{k}}
\phi^{j_{1}}_{m_{1}}\cdots \phi^{j_{l}}_{m_{l}}\one, z)w\rangle
z_{1}^{-p_{1}-1}\cdots z_{k}^{-p_{k}-1}\nn
&=\sum_{p_{1}, \dots, p_{k}\in \Z}\res_{\zeta_{1}=0}\cdots\res_{\zeta_{k}=0}\zeta_{1}^{p_{1}}\cdots \zeta_{k}^{p_{k}}
\res_{\xi_{1}=0}\cdots\res_{\xi_{l}=0}\xi_{1}^{m_{1}}\cdots \xi_{l}^{m_{l}}\cdot \nn
&\quad\quad\quad\;\cdot
F^{p}(\langle w', \phi_{W}^{i_{1}}(\zeta_{1}+z)\cdots \phi_{W}^{i_{k}}(\zeta_{k}+z)\phi_{W}^{j_{1}}(\xi_{1}+z)\cdots \phi_{W}^{j_{l}}(\xi_{l}+z)w\rangle)
z_{1}^{-p_{1}-1}\cdots z_{k}^{-p_{k}-1}.
\end{align}
We now expand
$$F^{p}(\langle w', \phi_{W}^{i_{1}}(\zeta_{1}+z)\cdots \phi_{W}^{i_{k}}(\zeta_{k}+z)\phi_{W}^{j_{1}}(\xi_{1}+z)\cdots \phi_{W}^{j_{l}}(\xi_{l}+z)w\rangle)$$
 as a Laurent series
$$\sum_{l\in \Z}f_{l}(\zeta_{1}, \dots, \zeta_{k-1}, \xi_{1}, \dots, \xi_{l}, z) \zeta_{k}^{-l-1}$$ 
in $\zeta_{k}$
in the region $|z|, |\zeta_{1}|, \dots, |\zeta_{k-1}|>|\zeta_{k}|>|\xi_{1}|, \dots, |\xi_{l}|$, where
$f_{l}(\zeta_{1}, \dots, \zeta_{k-1},  \xi_{1}, \dots, \xi_{l}, z)$ are analytic functions in $\zeta_{1}, \dots, \zeta_{k-1}$,
 $\xi_{1}, \dots, \xi_{l}$ and $z$.
Then in the region that the Laurent series expansion holds, we have
\begin{align}\label{phi-assoc-1}
\sum_{p_{k}\in \Z}&\res_{\zeta_{k}=0}\zeta_{k}^{p_{k}}
\left(\sum_{l\in \Z}f_{l}(\zeta_{1}, \dots, \zeta_{k-1}, \xi_{1}, \dots, \xi_{l}, z) \zeta_{k}^{-l-1}\right)
 z_{k}^{-p_{k}-1}\nn
&=\sum_{p_{k}\in \Z}f_{p_{k}}(\zeta_{1}, \dots, \zeta_{k-1}, \xi_{1}, \dots, \xi_{l}, z)  z_{k}^{-p_{k}-1}\nn
&=F^{p}(\langle w', \phi_{W}^{i_{1}}(\zeta_{1}+z)\cdots \phi_{W}^{i_{k-1}}(\zeta_{k-1}+z) \phi_{W}^{i_{k}}(z_{k}+z)
\phi_{W}^{j_{1}}(\xi_{1}+z)\cdots \phi_{W}^{j_{l}}(\xi_{l}+z)w\rangle).
\end{align}
Repeating this step for the variables $\zeta_{k-1}, \dots, \zeta_{1}$, we see that the right-hand side of
(\ref{phi-assoc}) is equal to the expansion of
\begin{equation}\label{phi-assoc-2}
\res_{\xi_{1}=0}\cdots\res_{\xi_{l}=0}\xi_{1}^{m_{1}}\cdots \xi_{l}^{m_{l}}
F^{p}(\langle w', \phi_{W}^{i_{1}}(z_{1}+z)\cdots \phi_{W}^{i_{k}}(z_{k}+z)\phi_{W}^{j_{1}}(\xi_{1}+z)\cdots \phi_{W}^{j_{l}}(\xi_{l}+z)w\rangle)
\end{equation}
as a  Laurent series in $z_{1}\dots, z_{k}$ in the region $|z|>|z_{1}|>\cdots>|z_{k}|>0$.
Thus the left-hand side of (\ref{phi-assoc})
is absolutely convergent to (\ref{phi-assoc-2}) in the region for this Laurent series expansion,
that is,  in the region $|z|>|z_{1}|>\cdots>|z_{k}|>0$,
\begin{align}\label{phi-assoc-3}
\langle w', &(Y^{g}_{W})^{p}(\phi^{i_{1}}(z_{1})\cdots \phi^{i_{k}}(z_{k})
\phi^{j_{1}}_{m_{1}}\cdots \phi^{j_{l}}_{m_{l}}\one, z)w\rangle\nn
&=\res_{\xi_{1}=0}\cdots\res_{\xi_{l}=0}\xi_{1}^{m_{1}}\cdots \xi_{l}^{m_{l}}\cdot\nn
&\quad\quad\quad\quad\cdot
F^{p}(\langle w', \phi_{W}^{i_{1}}(z_{1}+z)\cdots \phi_{W}^{i_{k}}(z_{k}+z)\phi_{W}^{j_{1}}(\xi_{1}+z)\cdots \phi_{W}^{j_{l}}(\xi_{l}+z)w\rangle).
\end{align}

On the other hand, we have
\begin{align}\label{iter-0}
\langle w', &(Y^{g}_{W})^{p}(Y_{V}(\phi^{i_{1}}_{n_{1}}\cdots \phi^{i_{k}}_{n_{k}}\one, z_{1}-z_{2})
\phi^{j_{1}}_{m_{1}}\cdots \phi^{j_{l}}_{m_{l}}\one, z_{2})w\rangle\nn
&=\sum_{n\in \Z}\langle w', (Y^{g}_{W})^{p}(e_{n}, z_{2})w\rangle\langle e_{n}',
Y_{V}(\phi^{i_{1}}_{n_{1}}\cdots \phi^{i_{k}}_{n_{k}}\one, z_{1}-z_{2})
\phi^{j_{1}}_{m_{1}}\cdots \phi^{j_{l}}_{m_{l}}\one\rangle\nn
&=\sum_{n\in \Z}\langle w', (Y^{g}_{W})^{p}(e_{n}, z_{2})w\rangle
\res_{\zeta_{1}=0}\cdots\res_{\zeta_{k}=0}\zeta_{1}^{n_{1}}\cdots \zeta_{k}^{n_{k}}\cdot\nn
&\quad\quad\quad\quad\quad\quad\cdot R(\langle e_{n}',\phi^{i_{1}}(\zeta_{1}+z_{1}-z_{2})\cdots \phi^{i_{k}}(\zeta_{k}+z_{1}-z_{2})
\phi^{j_{1}}_{m_{1}}\cdots \phi^{j_{l}}_{m_{l}}\one\rangle),
\end{align}
where we have used the definition of $Y_{V}$ in \cite{H-2-const}.
But by (\ref{phi-assoc-3}), in the region $|z_{2}|>|\zeta_{1}+z_{1}-z_{2}|>\cdots>|\zeta_{k}+z_{1}-z_{2}|>0$, $|\arg (\zeta_{k}+z_{1})-\arg z_{2}|<\frac{\pi}{2}, \dots, 
|\arg (\zeta_{1}+z_{1})-\arg z_{2}|<\frac{\pi}{2}$, we have
\begin{align}\label{iter-1}
\sum_{n\in \Z}&\langle w',  (Y^{g}_{W})^{p}(e_{n}, z_{2})w\rangle
\langle e_{n}',\phi^{i_{1}}(\zeta_{1}+z_{1}-z_{2})\cdots \phi^{i_{k}}(\zeta_{k}+z_{1}-z_{2})
\phi^{j_{1}}_{m_{1}}\cdots \phi^{j_{l}}_{m_{l}}\one\rangle\nn
&=\langle w', (Y^{g}_{W})^{p}(\phi^{i_{1}}(\zeta_{1}+z_{1}-z_{2})\cdots \phi^{i_{k}}(\zeta_{k}+z_{1}-z_{2})
\phi^{j_{1}}_{m_{1}}\cdots \phi^{j_{l}}_{m_{l}}\one, z_{2})w\rangle\nn
&=\res_{\xi_{1}=0}\cdots\res_{\xi_{l}=0}\xi_{1}^{m_{1}}\cdots \xi_{l}^{m_{l}}\cdot\nn
&\quad\quad\quad\quad\quad\cdot
F^{p}(\langle w', \phi_{W}^{i_{1}}(\zeta_{1}+z_{1})\cdots \phi_{W}^{i_{k}}(\zeta_{k}+z_{1})\phi_{W}^{j_{1}}(\xi_{1}+z_{2})\cdots \phi_{W}^{j_{l}}(\xi_{l}+z_{2})w\rangle).
\end{align}
The right-hand side of (\ref{iter-1}) is an analytic function in $\zeta_{1},
\dots, \zeta_{k}$, $z_{1}$ and $z_{2}$ of the form
\begin{align*}
&\sum_{i_{1}, \dots, i_{k}, i, n_{1}, \dots, n_{k}, n=0}^{N}f_{i_{1}
\cdots i_{k}in_{1}\cdots n_{k}n}(\zeta_{1}+z_{1}, \dots,
\zeta_{k}+z_{1}, z_{2})\cdot\nn
&\quad\quad\quad\quad\quad\quad \cdot 
e^{r^{(1)}_{i_{1}}l_{p}(\zeta_{1}+z_{1})}\cdots e^{r^{(k)}_{i_{k}}l_{p}(\zeta_{k}+z_{1})}
e^{r_{i}l_{p}(z_{2})}
(l_{p}(\zeta_{1}+z_{1}))^{n_{1}}\cdots (l_{p}(\zeta_{k}+z_{1}))^{n_{k}}(l_{p}(z_{2}))^{n},
\end{align*}
where 
$f_{i_{1}\cdots i_{k}in_{1}\cdots n_{k}n}(\zeta_{1}+z_{1}, \dots, \zeta_{k}+z_{1}, z_{1}, z_{2})$
for $i_{1}, \dots, i_{k}, i, n_{1}, \dots, n_{k}, n=0, \dots, N$
are rational functions in $\zeta_{1}, \dots, \zeta_{k}, z_{1}, z_{2}$ with the only possible poles
$\zeta_{i}-\zeta_{j}=0$ for $i\ne j$ and $\zeta_{i}+z_{1}-z_{2}=0$.
There is a unique expansion of such an analytic function
in the region $|z_{2}|>|\zeta_{1}+z_{1}-z_{2}|, \dots, |\zeta_{k}+z_{1}-z_{2}|>0$,
$|\arg (\zeta_{k}+z_{1})-\arg z_{2}|<\frac{\pi}{2}, \dots, 
|\arg (\zeta_{1}+z_{1})-\arg z_{2}|<\frac{\pi}{2}$,
$\zeta_{i}\ne \zeta_{j}$ for $i\ne j$,
$i, j=1, \dots, k$, such that each term
is a product of two analytic functions, one being analytic in $z_{2}$  of the form 
$$\sum_{j, q=0}^{M}b_{jq}e^{s_{j}l_{p}(z_{2})}(l_{p}(z_{2}))^{q}$$
and the other being a rational function in $\zeta_{1}+z_{1}-z_{2},
\dots, \zeta_{k}+z_{1}-z_{2}$ and $z_{1}$ with the only possible poles $\zeta_{i}-\zeta_{j}=0$ 
for $i\ne j$ and $\zeta_{i}+z_{1}-z_{2}=0$. Since
$$\sum_{n\in \Z}\langle w', (Y^{g}_{W})^{p}(e_{n}, z_{2})w\rangle
R(\langle e_{n}',\phi^{i_{1}}(\zeta_{1}+z_{1}-z_{2})\cdots \phi^{i_{k}}(\zeta_{k}+z_{1}-z_{2})
\phi^{j_{1}}_{m_{1}}\cdots \phi^{j_{l}}_{m_{l}}\one\rangle)$$
is a series of the same form  and is equal to
the left-hand side of (\ref{iter-1})  in the region
$|z_{2}|>|\zeta_{1}+z_{1}-z_{2}|>\cdots>|\zeta_{k}+z_{1}-z_{2}|>0$,
it must be absolutely convergent to the right-hand side of (\ref{iter-1}) in the larger region
$|z_{2}|>|\zeta_{1}+z_{1}-z_{2}|, \dots, |\zeta_{k}+z_{1}-z_{2}|>0$, 
$|\arg (\zeta_{k}+z_{1})-\arg z_{2}|<\frac{\pi}{2}, \dots, 
|\arg (\zeta_{1}+z_{1})-\arg z_{2}|<\frac{\pi}{2}$,. Therefore we obtain
\begin{align}\label{iter-2}
\sum_{n\in \Z}&\langle w', (Y^{g}_{W})^{p}(e_{n}, z_{2})w\rangle
R(\langle e_{n}',\phi^{i_{1}}(\zeta_{1}+z_{1}-z_{2})\cdots \phi^{i_{k}}(\zeta_{k}+z_{1}-z_{2})
\phi^{j_{1}}_{m_{1}}\cdots \phi^{j_{l}}_{m_{l}}\one\rangle)\nn
&=\res_{\xi_{1}=0}\cdots\res_{\xi_{l}=0}\xi_{1}^{m_{1}}\cdots \xi_{l}^{m_{l}}\cdot\nn
&\quad\quad\quad\quad\quad\cdot
F^{p}(\langle w', \phi_{W}^{i_{1}}(\zeta_{1}+z_{1})\cdots \phi_{W}^{i_{k}}(\zeta_{k}+z_{1})
\phi_{W}^{j_{1}}(\xi_{1}+z_{2})\cdots \phi_{W}^{j_{l}}(\xi_{l}+z_{2})w\rangle).
\end{align}
in the region $|z_{2}|>|\zeta_{1}+z_{1}-z_{2}|, \dots, |\zeta_{k}+z_{1}-z_{2}|>0$.
Thus when $|z_{2}|>|z_{1}-z_{2}|>0$, the right-hand side of (\ref{iter-0}) is absolutely convergent to
\begin{align*}
\res_{\zeta_{1}=0}&\cdots\res_{\zeta_{k}=0}\zeta_{1}^{n_{1}}\cdots \zeta_{k}^{n_{k}}
\res_{\xi_{1}=0}\cdots\res_{\xi_{l}=0}\xi_{1}^{m_{1}}\cdots \xi_{l}^{m_{l}}\cdot\nn
&\quad\quad\quad\quad\cdot F^{p}(\langle w', \phi_{W}^{i_{1}}(\zeta_{1}+z_{1})\cdots 
\phi_{W}^{i_{k}}(\zeta_{k}+z_{1})\phi_{W}^{j_{1}}(\xi_{1}+z_{2})\cdots \phi_{W}^{j_{l}}(\xi_{l}+z_{2})w\rangle),
\end{align*}
which has been proved above to be equal to the left hand side of (\ref{prod}) in 
the region $|z_{1}|>|z_{2}|>0$. The associativity is proved.

To prove the uniqueness, we need only show that any twisted $V$-module structure
on $W$ 
must have the vertex operator map defined by (\ref{vo}). 
But this is clear from the motivation that we have discussed before
the definition (\ref{vo}) of the vertex operator map $Y^{g}_{W}$.
\epfv

We shall say that the twisted vertex
operator map $Y_{W}^{g}$ is generated by the twisted fields 
$\phi^{i}(x)$ for $i\in I$. The $g$-twisted $V$-module $(W, Y_{W}^{g})$ 
is in fact generated by the coefficients of 
$\psi_{W}^{a}(x)v$ for $a\in A$ and $v\in V$. 

\begin{rema}
{\rm In this paper, we formulate and prove all our results for lower-bounded 
generalized twisted modules mainly because the explicit construction
in the next section gives in general only such twisted modules. But Theorem \ref{const-thm}
can be used to construct all different classes of twisted modules. 
If homogeneous subspaces of $W$ are finite dimensional, we obtain 
a grading-restricted generalized $g$-twisted $V$-module. If in addition
$L_{W}(0)$ acts on $W$ semisimply,  we obtain a $g$-twisted $V$-module. 
In the special case that $g=1_{V}$, Theorem \ref{const-thm} can be used to 
construct lower-bounded generalized $V$-modules, grading-restricted generalized 
$V$-modules and $V$-modules.}
\end{rema}

\renewcommand{\theequation}{\thesection.\arabic{equation}}
\renewcommand{\thethm}{\thesection.\arabic{thm}}
\setcounter{equation}{0}
\setcounter{thm}{0}
\section{An explicit construction of lower-bounded generalized
twisted modules satisfying a universal property}

In this section, we give an explicit construction of lower-bounded generalized
$g$-twisted $V$-modules satisfying
a universal property. As a consequence, every lower-bounded generalized
$g$-twisted $V$-module is the quotient of such a universal 
 lower-bounded generalized
$g$-twisted $V$-module.

We still assume  in this section that $V$ and $g$ satisfy Assumption \ref{algebra}. But we 
do not assume that we have the space, fields and operators in Data \ref{data}. 
In particular, we do not assume that Assumption \ref{basic-properties} holds.

Let 
$$\hat{V}_{\phi}^{[g]}
=\coprod_{i\in I, k\in \N}\C\mathcal{N}_{g}^{k}\phi^{i}_{-1}\one
\otimes t^{\alpha^{i}}\C[t, t^{-1}]\oplus \C L_{0}\oplus \C L_{-1},$$
where $L_{0}$ and $L_{-1}$ are fixed abstract basis elements of a vector space 
$\C L_{0}\oplus \C L_{-1}$.  Let $T(\hat{V}_{\phi}^{[g]})$ be the tensor algebra of
$\hat{V}_{\phi}^{[g]}$ and let 
$$\phi_{\hat{V}_{\phi}^{[g]}}^{i}(x)
=\sum_{n\in \alpha^{i}+\Z}((x^{-\mathcal{N}_{g}}
\phi^{i}_{-1}\one)\otimes t^{n})x^{-n-1}
\in x^{-\alpha^{i}}\hat{V}_{\phi}^{[g]}[[x, x^{-1}]][\log x]$$
for $i\in I$.
Then $\phi_{\hat{V}_{\phi}^{[g]}}^{i}(x)$ for $i\in I$ can be viewed as formal series of operators
 on $T(\hat{V}_{\phi}^{[g]})$. Also $L_{0}$ and $L_{-1}$ can  be viewed as operators
on $T(\hat{V}_{\phi}^{[g]})$. We shall use $L_{\hat{V}_{\phi}^{[g]}}(0)$ and $L_{\hat{V}_{\phi}^{[g]}}(-1)$
to denote the operators corresponding to $L_{0}$ and $L_{-1}$, respectively. 

For $i, j\in I$, we can always find $M_{i, j}\in \Z_{+}$ such that 
$x_{0}^{M_{i, j}}Y_{V}(\phi^{i}_{-1}\one, x_0)\phi^{j}_{-1}\one$
is a power series in $x_{0}$. For each pair $i, j\in I$, we choose $M_{i, j}$ to be 
the smallest of such positive integers. 
Let $J(\hat{V}_{\phi}^{[g]})$ be the ideal of 
$T(\hat{V}_{\phi}^{[g]})$ generated by the coefficients
of the formal series
\begin{align*}
(x_{1}-x_{2})^{M_{ij}}\phi_{\hat{V}_{\phi}^{[g]}}^{i}(x_{1})
\phi_{\hat{V}_{\phi}^{[g]}}^{j}(x_{2})
&-(-1)^{|\phi^{i}||\phi^{j}|}(x_{1}-x_{2})^{M_{ij}}
\phi_{\hat{V}_{\phi}^{[g]}}^{j}(x_{2})
\phi_{\hat{V}_{\phi}^{[g]}}^{i}(x_{1}),\\
L_{0}\phi_{\hat{V}_{\phi}^{[g]}}^{i}(x)-\phi_{\hat{V}_{\phi}^{[g]}}^{i}(x)L_{0}
&-x\frac{d}{dx}\phi_{\hat{V}_{\phi}^{[g]}}^{i}(x)-(\wt \phi^{i})\phi_{\hat{V}_{\phi}^{[g]}}^{i}(x),\\
L_{-1}\phi_{\hat{V}_{\phi}^{[g]}}^{i}(x)-\phi_{\hat{V}_{\phi}^{[g]}}^{i}(x)L_{-1}
&-\frac{d}{dx}\phi_{\hat{V}_{\phi}^{[g]}}^{i}(x)
\end{align*}
for $i, j\in I$, where the tensor product symbol $\otimes$ is omitted.
 Let $U(\hat{V}_{\phi}^{[g]})
=T(\hat{V}_{\phi}^{[g]})/J(\hat{V}_{\phi}^{[g]})$. 
Then $\phi_{\hat{V}_{\phi}^{[g]}}^{i}(x)$, $L_{\hat{V}_{\phi}^{[g]}}(0)$ and $L_{\hat{V}_{\phi}^{[g]}}(-1)$
can be viewed as formal series of operators and operators 
on $U(\hat{V}_{\phi}^{[g]})$ satisfying the weak commutativity
\begin{equation}\label{phi-weak-commty-0}
(x_{1}-x_{2})^{M_{ij}}\phi_{\hat{V}_{\phi}^{[g]}}^{i}(x_{1})
\phi_{\hat{V}_{\phi}^{[g]}}^{j}(x_{2})
=(-1)^{|\phi^{i}||\phi^{j}|}(x_{1}-x_{2})^{M_{ij}}
\phi_{\hat{V}_{\phi}^{[g]}}^{j}(x_{2})
\phi_{\hat{V}_{\phi}^{[g]}}^{i}(x_{1}),
\end{equation}
the $L(0)$-commutator formula
\begin{equation}\label{phi-L(0)-commutator-0}
L_{\hat{V}_{\phi}^{[g]}}(0)\phi_{\hat{V}_{\phi}^{[g]}}^{i}(x)
-\phi_{\hat{V}_{\phi}^{[g]}}^{i}(x)L_{\hat{V}_{\phi}^{[g]}}(0)
=x\frac{d}{dx}\phi_{\hat{V}_{\phi}^{[g]}}^{i}(x)+(\wt \phi^{i})\phi_{\hat{V}_{\phi}^{[g]}}^{i}(x)
\end{equation}
and the $L(-1)$-commutator formula
\begin{equation}\label{phi-L(-1)-commutator-0}
L_{\hat{V}_{\phi}^{[g]}}(-1)\phi_{\hat{V}_{\phi}^{[g]}}^{i}(x)
-\phi_{\hat{V}_{\phi}^{[g]}}^{i}(x)L_{\hat{V}_{\phi}^{[g]}}(-1)
=\frac{d}{dx}\phi_{\hat{V}_{\phi}^{[g]}}^{i}(x).
\end{equation}

Let $M$ be a $\Z_{2}$-graded vector space (graded by $\Z_{2}$-fermion numbers). Assume that $g$ acts on $M$ and there is an operator 
$L_{M}(0)$ on $M$. If $M$ is finite dimensional, then there exist operators $\mathcal{L}_{g}$,
$\mathcal{S}_{g}$, $\mathcal{N}_{g}$ such that on $M$, $g=e^{2\pi i\mathcal{L}_{g}}$
and $\mathcal{S}_{g}$ and $\mathcal{N}_{g}$ are the semisimple and nilpotent, respectively,
parts of $\mathcal{L}_{g}$. In this case, $M$ is also a direct sum of generalized 
eigenspaces for the operator $L_{M}(0)$ and $L_{M}(0)$ can be decomposed as 
the sum of its semisimple part $L_{M}(0)_{S}$ and nilpotent part $L_{M}(0)_{N}$. 
Moreover, the real parts of the eigenvalues of 
$L_{M}(0)$ has a lower bound. In the case that $M$ is infinite dimensional, we 
assume that all of these properties for $g$ and $L_{M}(0)$ hold. 
We call the 
eigenvalue of a generalized eigenvector 
$w\in M$ for $L_{M}(0)$ the {\it weight} of $w$ and denote 
it by $\wt w$. 
Let $\{w^{a}\}_{a\in A}$ be a basis of $M$ consisting of vectors homogeneous 
in weights, $\Z_{2}$-fermion numbers and 
$g$-weights (eigenvalues of $g$) such that for $a\in A$, 
either $L_{M}(0)_{N}w^{a}=0$ or there exists $L_{M}(0)_{N}(a)\in A$ 
such that $L_{M}(0)_{N}w^{a}=w^{L_{M}(0)_{N}(a)}$. For simplicity, when 
$L_{M}(0)_{N}w^{a}=0$, we shall use $w^{L_{M}(0)_{N}(a)}$ to denote $0$. 
Then for $a\in A$, we always have $L_{M}(0)_{N}w^{a}=w^{L_{M}(0)_{N}(a)}$.
For $a\in A$, let $\alpha^{a}\in \C$ such that $\Re(\alpha^{a})
\in [0, 1)$ and $e^{2\pi i\alpha^{a}}$ is the eigenvalue of $g$ for the generalized 
eigenvector $w^{a}$.

Let
$$\widetilde{M}^{[g]}=\coprod_{\alpha\in P_{V}}
U(\hat{V}_{\phi}^{[g]})\otimes (M\otimes 
t^{\alpha}\C[t, t^{-1}])\otimes V^{[\alpha]}.$$
Then $\widetilde{M}^{[g]}$ is a left $U(\hat{V}_{\phi}^{[g]})$-module.
In particular, $\phi_{\hat{V}_{\phi}^{[g]}}^{i}(x)$, $L_{\hat{V}_{\phi}^{[g]}}(0)$ 
and $L_{\hat{V}_{\phi}^{[g]}}(-1)$
act on $\widetilde{M}^{[g]}$  such that 
(\ref{phi-weak-commty-0}), (\ref{phi-L(0)-commutator-0}) and 
(\ref{phi-L(-1)-commutator-0}) hold for these operators. We shall denote 
their actions on $\widetilde{M}^{[g]}$ by 
$\phi_{\widetilde{M}^{[g]}}^{i}(x)$, $L_{\widetilde{M}^{[g]}}(0)$ 
and $L_{\widetilde{M}^{[g]}}(-1)$.
The actions of $g$, $e^{2\pi i\mathcal{S}_{g}}$ and 
$\mathcal{N}_{g}$ on $V$ and $M$ induce actions of $g$, 
$e^{2\pi i\mathcal{S}_{g}}$ and 
$\mathcal{N}_{g}$
on 
$\widetilde{M}^{[g]}$.

For $i\in I$, let  $K^{i}\in \N$ such that $\mathcal{N}_{g}^{K^{i}+1}\phi^{i}_{-1}\one=0$ and
we denote the actions of the elements 
$\frac{(-1)^{k}}{k!}
(\mathcal{N}_{g}^{k}\phi^{i}_{-1}\one)\otimes t^{n}$ for $n\in \alpha+\Z$ and $k=0, \dots,
K^{i}$ of $U(\hat{V}_{\phi}^{[g]})$ on $\widetilde{M}^{[g]}_{\ell}$ 
by $(\phi^{i}_{\widetilde{M}^{[g]}})_{n, k}$. 
For $a\in A$, $\alpha\in P_{V}$
and $n\in \alpha+\Z$, $w^{a}\otimes t^{n}$ can be viewed as a linear map 
from $V^{[\alpha]}$ to $\widetilde{M}^{[g]}$. We extend this map 
to a map from $V$ to $\widetilde{M}^{[g]}$ by mapping 
$V^{[\alpha']}$ to $0$ for $\alpha'\ne \alpha$.
We shall denote this map by $(\psi_{\widetilde{M}^{[g]}}^{a})_{n, 0}$. 
In general, for $n\in \alpha+\Z$ and $k\in \N$, 
we denote the linear map 
$v\mapsto \frac{(-1)^{k}}{k!}(\psi_{\widetilde{M}^{[g]}}^{a})_{n, 0}
\mathcal{N}_{g}^{k}v$ 
from $V^{[\alpha]}$ to $\widetilde{M}^{[g]}$
by $(\psi_{\widetilde{M}^{[g]}}^{a})_{n, k}$ and extend it to a linear map 
from $V$ to $\widetilde{M}^{[g]}$ in the same way. 
Then $\widetilde{M}^{[g]}$ is spanned by elements of the form
\begin{equation}\label{element-form}
(\phi^{i_{1}}_{\widetilde{M}^{[g]}})_{n_{1}, k_{1}}
\cdots (\phi^{i_{l}}_{\widetilde{M}^{[g]}})_{n_{l}, k_{l}}(L_{\widetilde{M}^{[g]}}(m))^{q})(\psi_{\widetilde{M}^{[g]}}^{a})_{n, k}v,
\end{equation}
for $i_{1}, \dots, i_{l}\in I$,  
$n_{1}\in \alpha^{i_{1}}+\Z,
\dots, n_{l}\in \alpha^{i_{l}}+\Z$, $0\le k_{1}\le K^{i_{1}}, \dots, 0\le k_{l}\le K^{i_{l}}$, 
$m=0, -1$, $q\in \N$,
$a\in A$, $n \in \alpha+\Z$, $0\le k\le K$, $v\in V^{[\alpha]}$,
$\alpha\in P_{V}$, where $K\in \N$ satisfying $\mathcal{N}_{g}^{K+1}v=0$. 
We already know that $M$ is graded by eigenvalues of $L_{M}(0)$.  
For the element (\ref{element-form}) with homogeneous  $v$, we define its
weight to be 
$$\wt \phi^{i_{1}}-n_{1}-1+\cdots +\phi^{i_{l}}-n_{l}-1-m+\wt w^{a} -n-1+\wt v.$$
Then
$$\widetilde{M}^{[g]}=\coprod_{n\in \C}\widetilde{M}^{[g]}_{[n]},$$
where $\widetilde{M}^{[g]}_{[n]}$ is the subspace of $\widetilde{M}^{[g]}$ 
consisting of elements of weight $n$. 

For $i\in I$, 
we have the formal series of operators on $\widetilde{M}^{[g]}$
\begin{align}\label{defn-phi}
\phi_{\widetilde{M}^{[g]}}^{i}(x)&
=\sum_{k=0}^{K_{i}}\sum_{n\in \alpha^{i}+\Z}
(\phi_{\widetilde{M}^{[g]}}^{i})_{n, k}x^{-n-1}(\log x)^{k}\nn
&=\sum_{n\in \alpha^{i}+\Z}
((x^{-\mathcal{N}_{g}}\phi^{i}_{-1}\one)\otimes t^{n})x^{-n-1}\nn
&=(x^{-\mathcal{N}_{g}}\phi^{i}_{-1}\one)
 \otimes \left(\frac{t}{x}\right)^{\alpha^{i}}
x^{-1}\delta\left(\frac{t}{x}\right).
\end{align}
Recall that $\phi_{\widetilde{M}^{[g]}}^{i}(x)$ is in fact the action of 
$\phi_{\hat{V}_{\phi}^{[g]}}^{i}(x)$ on $\widetilde{M}^{[g]}$. 
 For $v\in V$, there exists $K^{v}\in \N$ such that $\mathcal{N}_{g}^{K^{v}+1}v=0$. 
For $a\in A$ and $v\in V^{[\alpha]}$, let
\begin{align}\label{defn-psi}
\psi_{\widetilde{M}^{[g]}}^{a}(x)v
&=\sum_{k=0}^{K^{v}}\sum_{n\in \alpha+\Z}
(\psi_{\widetilde{M}^{[g]}}^{a})_{n, k}v x^{-n-1}(\log x)^{k}\nn
&=\sum_{n\in \alpha+\Z}
(w^{a}\otimes t^{n}) x^{-\mathcal{N}_{g}}v x^{-n-1}\nn
&=\left(w^{a}\otimes  \left(\frac{t}{x}\right)^{\alpha}
x^{-1}\delta\left(\frac{t}{x}\right)\right) (x^{-\mathcal{N}_{g}}v).
\end{align}
It is  a series in $x$ with coefficients in $\widetilde{M}^{[g]}$.
Then $\psi_{\widetilde{M}^{[g]}}^{a}(x)$ is a formal series with coefficients 
in $\hom(V, \widetilde{M}^{[g]})$.

Let $B\in \R$ such that $B\le \Re(\wt w)$ for any generalized eigenvector $w\in M$  of
$L_{M}(0)$. Such $B$ exists because the real parts of the eigenvalues of 
$L_{M}(0)$ is lower bounded.
Let $J_{B}(\widetilde{M}^{[g]})$ be the $U(\hat{V}_{\phi}^{[g]})$-submodule
of $\widetilde{M}^{[g]}$ generated by elements of the following forms:
(i) $(\psi_{\widetilde{M}^{[g]}}^{a})_{n, 0}\one$ for  $a\in A$,
and $n\not\in -\N-1$; (ii) (\ref{element-form}) for $i_{1}, \dots, i_{l}\in I$,  
$n_{1}\in \alpha^{i_{1}}+\Z,
\dots, n_{l}\in \alpha^{i_{l}}+\Z$, $0\le k_{1}\le K^{i_{1}}, \dots, 0\le k_{l}\le K^{i_{l}}$, 
$m=0, -1$, $a\in A$, $n \in \alpha+\Z$, $0\le k\le K$, $v\in V^{[\alpha]}$,
$\alpha\in P_{V}$ such that 
$$\Re(\wt \phi^{i_{1}}-n_{1}-1+\cdots +\phi^{i_{l}}-n_{l}-1-m+\wt w^{a} -n-1+\wt v)<B.$$
Consider the quotient  $U(\hat{V}_{\phi}^{[g]})$-module
 $\widetilde{M}^{[g]}/J_{B}(\widetilde{M}^{[g]})$. Since 
$J_{B}(\widetilde{M}^{[g]})$ is spanned by homogeneous elements,
$\widetilde{M}^{[g]}/J_{B}(\widetilde{M}^{[g]})$ is also graded.
In addition, $\widetilde{M}^{[g]}/J_{B}(\widetilde{M}^{[g]})$ is lower bounded 
with respect to the weight grading with a lower bound $B$. We shall still use the same notations 
to denote the elements of this quotient and operators on this quotient. Then in this quotient
$(\psi_{\widetilde{M}^{[g]}}^{a})_{n, 0}\one=0$ for  $a\in A$,
and $n\not\in -\N-1$.

For $i\in I$, $a\in A$, $n\in \C$ and $k\in \N$,  by (\ref{phi-L(0)-commutator-0})
$$\wt ((\phi_{\widetilde{M}^{[g]}}^{i})_{n, k}(\psi_{\widetilde{M}^{[g]}}^{a})_{-1, 0}\one)
=\wt \phi^{i}-n-1+\wt w^{a}.$$
Since  $B$ is a lower bound of the real parts of the weights of $\widetilde{M}^{[g]}/J_{B}(\widetilde{M}^{[g]})$ and 
$$\phi_{\widetilde{M}^{[g]}}^{i}(x)\in x^{-\alpha^{i}}
\left(\widetilde{M}^{[g]}/J_{B}(\widetilde{M}^{[g]})\right)[[x, x^{-1}]],$$
we have 
$$(\phi_{\widetilde{M}^{[g]}}^{i})_{n, k}(\psi_{\widetilde{M}^{[g]}}^{a})_{-1, 0}\one=0$$
when $\wt (\phi_{\widetilde{M}^{[g]}}^{i})_{n, k}w^{a}=\alpha^{i}+m$ where 
$m\in \Z$ and $m>\wt \phi^{i}-1+\Re(\wt w^{a})-B-\Re(\alpha^{i})$. 
For $i\in I$ and $a\in A$, let $M_{i, a}\in \Z_{+}$ be the smallest of 
$m\in \Z$ such that $m>\wt \phi^{i}-1+\Re(\wt w^{a})-B-\Re(\alpha^{i})$. Then 
$$x^{\alpha^{i}+M_{i, a}}
\phi_{\widetilde{M}^{[g]}}^{i}(x)(\psi_{\widetilde{M}^{[g]}}^{a})_{-1, 0}\one$$
 is a powers series in $x$ with polynomials in $\log x$ as coefficients.  Also
for $a\in A$ and $v\in V$, 
$\psi_{\widetilde{M}^{[g]}}^{a}(x)v$ has only finitely many 
terms with negative real parts of powers of $x$. 
In particular, for $i\in I$, $a\in A$ and $v\in V^{[\alpha]}$,
\begin{align}\label{phi-psi-element}
(x_{1}&-x_{2})^{\alpha^{i}+M_{i, a}} (x_{1}-x_{2})^{\mathcal{N}_{g}}
\phi_{\widetilde{M}^{[g]}}^{i}(x_{1}) 
(x_{1}-x_{2})^{-\mathcal{N}_{g}}
\psi_{\widetilde{M}^{[g]}}^{a}(x_{2})v\nn
&-(-1)^{|u||w|} (-x_{2}+x_{1})^{\alpha^{i}+M_{i, a}}
\psi_{\widetilde{M}^{[g]}}^{a}(x_{2})
(-x_{2}+x_1)^{\mathcal{N}_{g}} \phi^{i}(x_{1}) 
(-x_{2}+x_1)^{-\mathcal{N}_{g}}v,
\end{align}
in $x_{1}$ and $x_{2}$ are well defined as a formal series with coefficients in 
$\widetilde{M}^{[g]}/J_{B}(\widetilde{M}^{[g]})$.
Let $J(\widetilde{M}^{[g]}/J_{B}(\widetilde{M}^{[g]}))$ 
be the $U(\hat{V}_{\phi}^{[g]})$-submodule
of $\widetilde{M}^{[g]}/J_{B}(\widetilde{M}^{[g]})$ generated by 
 the coefficients of the formal series (\ref{phi-psi-element}) for 
$i\in I$, $a\in A$ and $v\in V^{[\alpha]}$ and the coefficients of the formal series 
\begin{align*}
L_{\widetilde{M}^{[g]}}(0)\psi_{\widetilde{M}^{[g]}}^{a}(x)v-
\psi_{\widetilde{M}^{[g]}}^{a}(x)L_{V}(0)v
&-x\frac{d}{dx}\psi_{\widetilde{M}^{[g]}}^{a}(x)v-
(\wt w^{a})\psi_{\widetilde{M}^{[g]}}^{a}(x)v-\psi_{\widetilde{M}^{[g]}}^{L_{M}(0)_{N}a}(x)v,
\\
L_{\widetilde{M}^{[g]}}(-1)\psi_{\widetilde{M}^{[g]}}^{a}(x)v&-
\psi_{\widetilde{M}^{[g]}}^{a}(x)L_{V}(-1)v
-\frac{d}{dx}\psi_{\widetilde{M}^{[g]}}^{a}(x)v
\end{align*}
for $a\in A$ and $v\in V$.
We then have a quotient $U(\hat{V}_{\phi}^{[g]})$-module
$$\widehat{M}^{[g]}_{B}=(\widetilde{M}^{[g]}_{\ell}/J_{B}(\widetilde{M}^{[g]}))
/J(\widetilde{M}^{[g]}/J_{B}(\widetilde{M}^{[g]})).$$
Again, we shall use the same notations
for the elements of $\widetilde{M}^{[g]}$ 
to denote the corresponding elements of $\widehat{M}^{[g]}_{B}$. We shall use
$\phi_{\widehat{M}^{[g]}_{B}}^{i}(x)$, $\psi_{\widehat{M}^{[g]}_{B}}^{a}(x)$,
$L_{\widehat{M}^{[g]}_{B}}(0)$
and $L_{\widehat{M}^{[g]}_{B}}(-1)$ to denote the series of operators and the operators
on $\widehat{M}^{[g]}_{B}$ induced from the corresponding series of 
operators and operators  on $\widetilde{M}^{[g]}$. 
Since $J(\widetilde{M}^{[g]}/J_{B}(\widetilde{M}^{[g]}))$
is spanned by homogeneous elements, the quotient 
$\widehat{M}^{[g]}_{B}=(\widetilde{M}^{[g]}_{\ell}/J_{B}(\widetilde{M}^{[g]}))
/J(\widetilde{M}^{[g]}/J_{B}(\widetilde{M}^{[g]}))$
is also graded  and is lower bounded with respect to the weight grading with a lower bound $B$. Moreover, 
in $\widehat{M}^{[g]}_{B}$, we have
\begin{equation}
(x_{1}-x_{2})^{M_{ij}}\phi_{\widehat{M}^{[g]}_{B}}^{i}(x_{1})
\phi_{\widehat{M}^{[g]}_{B}}^{j}(x_{2})
=(-1)^{|\phi^{i}||\phi^{j}|}(x_{1}-x_{2})^{M_{ij}}
\phi_{\widehat{M}^{[g]}_{B}}^{j}(x_{2})
\phi_{\widehat{M}^{[g]}_{B}}^{i}(x_{1}),\label{phi-weak-commty}
\end{equation}
\begin{equation}
L_{\widehat{M}^{[g]}_{B}}(0)\phi_{\widehat{M}^{[g]}_{B}}^{i}(x)
-\phi_{\widehat{M}^{[g]}_{B}}^{i}(x)L_{\widehat{M}^{[g]}_{B}}(0)
=x\frac{d}{dx}\phi_{\widehat{M}^{[g]}_{B}}^{i}(x)+(\wt \phi^{i})
\phi_{\widehat{M}^{[g]}_{B}}^{i}(x), \label{phi-L(0)-commutator}
\end{equation}
\begin{equation}
L_{\widehat{M}^{[g]}_{B}}(-1)\phi_{\widehat{M}^{[g]}_{B}}^{i}(x)-\phi_{\hat{V}_{\phi}^{[g]}}^{i}(x)
L_{\widehat{M}^{[g]}_{B}}(-1)
=\frac{d}{dx}\phi_{\widehat{M}^{[g]}_{B}}^{i}(x), \label{phi-L(-1)-commutator}
\end{equation}
\begin{align}
(x_{1}&-x_{2})^{\alpha^{i}+M_{i, a}} (x_{1}-x_{2})^{\mathcal{N}_{g}}
\phi_{\widehat{M}^{[g]}_{B}}^{i}(x_{1}) 
(x_{1}-x_{2})^{-\mathcal{N}_{g}}
\psi_{\widehat{M}^{[g]}_{B}}^{a}(x_{2})v\nn
&=(-1)^{|u||w|} (-x_{2}+x_{1})^{\alpha^{i}+M_{i, a}}
\psi_{\widehat{M}^{[g]}_{B}}^{a}(x_{2})
(-x_{2}+x_1)^{\mathcal{N}_{g}} \phi^{i}(x_{1}) 
(-x_{2}+x_1)^{-\mathcal{N}_{g}}v, \label{phi-psi-commty}
\end{align}
\begin{equation}
L_{\widehat{M}^{[g]}_{B}}(0)\psi_{\widehat{M}^{[g]}_{B}}^{a}(x)v-
\psi_{\widehat{M}^{[g]}_{B}}^{a}(x)L_{V}(0)v
=x\frac{d}{dx}\psi_{\widehat{M}^{[g]}_{B}}^{a}(x)+
(\wt w^{a})\psi_{\widehat{M}^{[g]}_{B}}^{a}(x)
+\psi_{\widehat{M}^{[g]}_{B}}^{L_{M}(0)_{N}a}(x),
\label{psi-L(0)-commutator}
\end{equation}
\begin{equation}
L_{\widehat{M}^{[g]}_{B}}(-1)\psi_{\widehat{M}^{[g]}_{B}}^{a}(x)v-
\psi_{\widehat{M}^{[g]}_{B}}^{a}(x)L_{V}(-1)v
=\frac{d}{dx}\psi_{\widehat{M}^{[g]}_{B}}^{a}(x)v.\label{psi-L(-1)-commutator}
\end{equation}
By (\ref{phi-L(0)-commutator}), (\ref{phi-L(-1)-commutator}), (\ref{psi-L(0)-commutator})
and (\ref{psi-L(-1)-commutator}), we see that $\widehat{M}^{[g]}_{B}$ is 
spanned by elements of the form 
\begin{equation}\label{element-form-no-L}
(\phi_{\widehat{M}^{[g]}_{B}}^{i_{1}})_{n_{1}, k_{1}}
\cdots (\phi_{\widehat{M}^{[g]}_{B}}^{i_{l}})_{n_{l}, k_{l}}
(\psi_{\widehat{M}^{[g]}_{B}}^{a})_{n, k}v.
\end{equation}

We now have the following main 
result giving an explicit construction of lower-bounded generalized $g$-twisted $V$-modules:

\begin{thm}\label{explicit-const}
The twisted fields 
$\phi_{\widehat{M}^{[g]}_{B}}^{i}$ for $i\in I$
generate a twisted vertex operator map 
$$Y^{g}_{\widehat{M}^{[g]}_{B}}: V\otimes \widehat{M}^{[g]}_{B}\to \widehat{M}^{[g]}_{B}\{x\}[\log x]$$
such that $(\widehat{M}^{[g]}_{B}, Y^{g}_{\widehat{M}^{[g]}_{B}})$ is a 
lower-bounded generalized $g$-twisted $V$-module. Moreover, this is the unique generalized 
$g$-twisted $V$-module structure on $\widehat{M}^{[g]}_{B}$ generated 
by the coefficients of $(\psi_{\widehat{M}^{[g]}_{B}}^{a})(x)v$ for $a\in A$ and $v\in V$
such that $Y^{g}_{\widehat{M}^{[g]}_{B}}(\phi^{i}_{-1}\one, z)=\phi_{\widehat{M}^{[g]}_{B}}^{i}(z)$
for $i\in I$.
\end{thm}
\pf
The space $\widehat{M}^{[g]}_{B}$ is graded by weights, $\Z_{2}$-fermion numbers and 
is a direct sum of generalized eigenspaces of an action of $g$. By construction, 
$(\widehat{M}^{[g]}_{B})_{[n]}=0$ when $\Re(n)<B$. 
We already have the linear maps $\phi_{\widehat{M}^{[g]}_{B}}^{i}$, 
$\psi_{\widehat{M}^{[g]}_{B}}^{a}$, $L_{\widehat{M}^{[g]}_{B}}(0)$
and $L_{\widehat{M}^{[g]}_{B}}(-1)$. 
We need only verify Properties \ref{property-1}--\ref{property-7} in 
Assumption \ref{basic-properties}.

By (\ref{phi-L(0)-commutator}) and (\ref{psi-L(0)-commutator}),
Property \ref{property-1} in Assumption \ref{basic-properties} holds. 

By (\ref{phi-L(-1)-commutator}) and (\ref{psi-L(-1)-commutator}), Property \ref{property-2} in 
Assumption \ref{basic-properties} holds. 

By the definition of $\widehat{M}^{[g]}_{B}$, Property \ref{property-3} in 
Assumption \ref{basic-properties} holds. 

By the definition of $\widehat{M}^{[g]}_{B}$, Property \ref{property-4} in 
Assumption \ref{basic-properties} holds. 

For $i\in I$ and $p\in \Z$, by (\ref{defn-phi}) and the definition of the actions of $g$, $e^{2\pi i\mathcal{S}_{g}}$ and $\mathcal{N}_{g}$ on $\widehat{M}^{[g]}_{B}$, we have
\begin{align*}
g\phi_{\widehat{M}^{[g]}_{B}}^{i; p+1}(z)g^{-1}
&=e^{2\pi i\mathcal{S}_{g}}e^{2\pi i\mathcal{N}_{g}}
\left((e^{-l_{p+1}(z)\mathcal{N}_{g}}\phi^{i}_{-1}\one)
 \otimes \frac{t^{\alpha^{i}}}{e^{\alpha^{i} l_{p+1}(z)}}
x^{-1}\delta\left(\frac{t}{z}\right)\right)e^{-2\pi i\mathcal{N}_{g}}
e^{-2\pi i\mathcal{S}_{g}}\nn
&=(e^{-(l_{p+1}(z)-2\pi i)\mathcal{N}_{g}}e^{2\pi i\mathcal{S}_{g}}\phi^{i}_{-1}\one)
\otimes \frac{t^{\alpha^{i}}}{e^{\alpha^{i} l_{p+1}(z)}}
z^{-1}\delta\left(\frac{t}{z}\right)\nn
&=(e^{-(l_{p+1}(z)-2\pi i)\mathcal{N}_{g}}\phi^{i}_{-1}\one)
\otimes \frac{t^{\alpha^{i}}}{e^{\alpha^{i} (l_{p+1}(z)-2\pi i)}}
z^{-1}\delta\left(\frac{t}{z}\right)\nn
&=(e^{-l_{p}(z)\mathcal{N}_{g}}a)
\otimes \frac{t^{\alpha^{i}}}{e^{\alpha^{i} l_{p}(z)}}
z^{-1}\delta\left(\frac{t}{z}\right)\quad\quad\quad\quad\quad\quad\quad\quad\quad\quad\quad\nn
&=\phi_{\widehat{M}^{[g]}_{B}}^{i; p}(z).
\end{align*}
This is Part (i) of Property \ref{property-5} in 
Assumption \ref{basic-properties}. 
Part (ii) of Property \ref{property-5} in 
Assumption \ref{basic-properties} follows immediately 
from  (\ref{defn-phi}), (\ref{defn-psi})  and the definition of the action of 
$\mathcal{N}_{g}$ on $\widehat{M}^{[g]}_{B}$. 
Part (iii) of Property \ref{property-5} in 
Assumption \ref{basic-properties}
follows immediately from the definition of the actions of 
$e^{2\pi i\mathcal{S}_{g}}$ and $\mathcal{N}_{g}$ on $\widehat{M}^{[g]}_{B}$
and Assumption \ref{algebra}. For $a\in A$ and $n\in -\N-1$, since $w^{a}$ is a generalized 
eigenvector of $g$ with eigenvalue  
$e^{2\pi i\alpha^{a}}$, by the definition of the action of $g$ on $\widehat{M}^{[g]}_{B}$,
$(\psi^{a}_{W})_{n, 0}\one$ is a generalized 
eigenvector of $g$ with eigenvalue  
$e^{2\pi i\alpha^{a}}$. This is Part (iv) of Property \ref{property-5} in 
Assumption \ref{basic-properties}. 

 Property \ref{property-6}  in 
Assumption \ref{basic-properties} in our case is in fact (\ref{phi-weak-commty}).

Property \ref{property-7}  in 
Assumption \ref{basic-properties} in our case is in fact (\ref{phi-psi-commty}).

Since the space $\widehat{M}^{[g]}_{B}$ equipped with $\phi_{\widehat{M}^{[g]}_{B}}^{i}$, 
$\psi_{\widehat{M}^{[g]}_{B}}^{a}$, $L_{\widehat{M}^{[g]}_{B}}(0)$
and $L_{\widehat{M}^{[g]}_{B}}(-1)$ satisfies Properties \ref{property-1}--\ref{property-7}
in Assumption \ref{basic-properties}, by Theorem \ref{const-thm}, we have 
a unique generalized 
$g$-twisted $V$-module structure on $\widehat{M}^{[g]}_{B}$ generated 
by the coefficients of $(\psi_{\widehat{M}^{[g]}_{B}}^{a})(x)v$ for $a\in A$ and $v\in V$
such that $Y^{g}_{\widehat{M}^{[g]}_{B}}(\phi^{i}_{-1}\one, z)=\phi_{\widehat{M}^{[g]}_{B}}^{i}(z)$
for $i\in I$.
\epfv

Now we prove a universal property of the generalized $g$-twisted $V$-module 
$\widehat{M}^{[g]}_{B}$
constructed in Theorem \ref{explicit-const}.

\begin{thm}
Let $(W, Y^{g}_{W})$ be a lower-bounded generalized $g$-twisted $V$-module
and $M_{0}$ a $\Z_{2}$-graded subspace of $W$ invariant under the actions of 
$g$, $\mathcal{S}_{g}$, $\mathcal{N}_{g}$, $L_{W}(0)$, 
$L_{W}(0)_{S}$ and
$L_{W}(0)_{N}$. Let $B\in \R$ such that $W_{[n]}=0$ when $\Re(n)<B$. 
Assume that there is a linear map $f: M\to M_{0}$ 
preserving the $\Z_{2}$-fermion number grading and commuting with the actions of 
$g$, $\mathcal{S}_{g}$, $\mathcal{N}_{g}$, $L_{W}(0)$ $(L_{\widehat{M}^{[g]}_{B}}(0))$, 
$L_{W}(0)_{S}$ $(L_{\widehat{M}^{[g]}_{B}}(0)_{S})$ and
$L_{W}(0)_{N}$ $(L_{\widehat{M}^{[g]}_{B}}(0)_{N})$. Then there exists a unique module 
map $\tilde{f}: \widehat{M}^{[g]}_{B}\to W$ such that $\tilde{f}|_{M}=f$. 
If $f$ is surjective and $(W, Y^{g}_{W})$ is generated by the 
coefficients of $(Y^{g})_{WV}^{W}(w_{0}, x)v$ for $w_{0}\in M_{0}$ and $v\in V$, 
where $(Y^{g})_{WV}^{W}$ 
is the twist vertex operator map obtained from $Y_{W}^{g}$ (see \cite{H-twist-vo}), then 
$\tilde{f}$ is surjective. 
\end{thm}
\pf
Note that $\widehat{M}^{[g]}_{B}$ 
is spanned by elements of the form (\ref{element-form-no-L}). 
We define $\tilde{f}$ by 
\begin{align*}
\tilde{f}&((\phi_{\widehat{M}^{[g]}_{B}}^{i_{1}})_{n_{1}, k_{1}}
\cdots (\phi_{\widehat{M}^{[g]}_{B}}^{i_{l}})_{n_{l}, k_{l}}
(\psi_{\widehat{M}^{[g]}_{B}}^{a})_{n, k}v)\nn
&=(Y_{W}^{g})_{n_{1}, k_{1}}(\phi^{i_{1}}_{-1}\one)
\cdots (Y_{W}^{g})_{n_{l}, k_{l}}(\phi^{i_{l}}_{-1}\one)
((Y^{g})_{WV}^{W})_{n, k}(f(w^{a}))v
\end{align*}
for $i_{1}, \dots, i_{l}\in I$, $n_{1}\in \alpha^{i_{1}}+\Z, \dots, n_{l}\in \alpha^{i_{l}}+\Z$,
$k_{1}, \dots, k_{l}\in \N$, $a\in A$, $\alpha\in P_{V}$, $n\in \alpha+\Z$, $k\in \N$
and $v\in V^{[\alpha]}$, where for $i\in I$, $n\in \alpha^{i}+\Z$ and $k\in \N$, 
$(Y_{W}^{g})_{n, k}(\phi^{i}_{-1}\one)$ is the coefficient of $x^{-n-1}(\log x)^{k}$ in the series
$Y_{W}^{g}(\phi^{i}_{-1}\one, x)$
and for $a\in A$, $n\in \C$ and $k\in \N$,  
$((Y^{g})_{WV}^{W})_{n, k}(f(w^{a}))$ is the coefficient of $x^{-n-1}(\log x)^{k}$
in the series $(Y^{g})_{WV}^{W}(f(w^{a}), x)$. 

We first need to show that $\tilde{f}$ is well defined. By the construction of 
$\widehat{M}^{[g]}_{B}$, we see that the only relations among elements of the 
form $(\phi_{\widehat{M}^{[g]}_{B}}^{i_{1}})_{n_{1}, k_{1}}
\cdots (\phi_{\widehat{M}^{[g]}_{B}}^{i_{l}})_{n_{l}, k_{l}}
(\psi_{\widehat{M}^{[g]}_{B}}^{a})_{n, k}v$ are the following:
$(\psi_{\widehat{M}^{[g]}_{B}}^{a})_{n, 0}\one=0$ for $a\in A$, $n\in \N-1$; 
$(\phi_{\widehat{M}^{[g]}_{B}}^{i_{1}})_{n_{1}, k_{1}}
\cdots (\phi_{\widehat{M}^{[g]}_{B}}^{i_{l}})_{n_{l}, k_{l}}
(\psi_{\widehat{M}^{[g]}_{B}}^{a})_{n, k}v=0$ when the real part of 
its weight is less than $B$; the coefficients of
(\ref{phi-weak-commty})--(\ref{psi-L(-1)-commutator}); 
the relations among $v\in V$. These relations also hold for 
elements of the form 
$$(Y_{W}^{g})_{n_{1}, k_{1}}(\phi^{i_{1}}_{-1}\one)
\cdots (Y_{W}^{g})_{n_{l}, k_{l}}(\phi^{i_{l}}_{-1}\one)
((Y^{g})_{WV}^{W})_{n, k}(f(w^{a}))v$$
of $W$ because $W$ is a 
lower-bounded generalized $g$-twisted 
$V$-module such that $W_{[n]}=0$ when $\Re(n)<B$, because the choices 
of $M_{i,j}$ for $i, j\in I$ and $M_{i, a}$ for $i\in I$ and $a\in A$
depend only on $\phi^{i}$ and $\phi^{j}$ and on $\phi^{i}$, $\wt w^{a}=\wt (\psi_{\widehat{M}^{[g]}_{B}}^{a})_{-1, 0}\one$
and $B$, respectively, and  because $f$ commutes with all the operators on $M$ and $M_{0}$.
Thus $\tilde{f}$ is indeed well defined. 

By the definition of $\tilde{f}$, it is a module map and $\tilde{f}|_{M}=f$. 
Since $\widehat{M}^{[g]}_{B}$ is determined uniquely by $M$, $\tilde{f}$ is unique.

If $f$ is surjective and $(W, Y^{g}_{W})$ is generated by the 
coefficients of $(Y^{g})_{WV}^{W}(w_{0}, x)v$ for $w_{0}\in M_{0}$ 
and $v\in V$, then $(W, Y^{g}_{W})$
is in fact generated by $(Y^{g})_{WV}^{W}(f(w), x)v$ for $w\in M$ and $v\in V$.
Since $\tilde{f}$ is a module map, we obtain $\tilde{f}(\widehat{M}^{[g]}_{B})=W$.
\epfv

Finally we have the following immediate consequence:

\begin{cor}\label{quotient}
Let $(W, Y^{g}_{W})$ be a lower-bounded generalized $g$-twisted $V$-module generated by the 
coefficients of $(Y^{g})_{WV}^{W}(w, x)v$ for $w\in M$, where $(Y^{g})_{WV}^{W}$ 
is the twist vertex operator map obtained from $Y_{W}^{g}$ (see \cite{H-twist-vo})
and $M$ is a $\Z_{2}$-graded subspace of $W$  invariant under the actions of 
$g$, $\mathcal{S}_{g}$, $\mathcal{N}_{g}$, $L_{W}(0)$, 
$L_{W}(0)_{S}$ and
$L_{W}(0)_{N}$. Let $B\in \R$ such that $W_{[n]}=0$ when $\Re(n)<B$. 
Then there is a generalized $g$-twisted $V$-submodule $J$ of 
$\widehat{M}^{[g]}_{B}$ such that $W$ is equivalent as a 
lower-bounded generalized $g$-twisted $V$-module to the quotient 
module $\widehat{M}^{[g]}_{B}/J$.
\end{cor}

\noindent {\small \sc Department of Mathematics, Rutgers University,
110 Frelinghuysen Rd., Piscataway, NJ 08854-8019}

\noindent {\em E-mail address}: yzhuang@math.rutgers.edu

\end{document}